\documentclass[11pt]{article}

\usepackage{amsmath}
\usepackage{amsfonts}
\usepackage{amssymb}
\usepackage{amsrefs}
\usepackage{mathtools}
\usepackage{showkeys}

\usepackage[draft]{todonotes}   

\usepackage{cases}
\usepackage{xcolor}

\newcommand{\dt}{\partial_t}
\newcommand{\dv}{\mathrm{div}}
\newcommand{\dr}{\partial_r}
\newcommand{\dx}{\partial_x}

\newcommand{\subeqref}[2]{$\eqref{#1}_{#2}$}
\newcommand{\abs}[2]{\bigl| #1 \bigr|^{#2}}
\newcommand{\norm}[2]{\bigl\Arrowvert #1 \bigr\Arrowvert_{#2}}
\newcommand{\supnorm}{L_t^{\infty}L_x^\infty}
\newcommand{\stnorm}[2]{L_t^{#1}L_x^{#2}}

\newtheorem{df}{Definition}
\newtheorem{thm}{Theorem}[section]
\newtheorem{lm}{Lemma}
\newtheorem{cor}{Corollary}
\newtheorem{prop}{Proposition}
\newenvironment{pf}{\paragraph{Proof}}{\hfill$\square$\\}

\newtheorem{rmk}{Remark}

\numberwithin{equation}{section}
\allowdisplaybreaks

\title{A Model of Radiational Gaseous Stars}

\author{Xin Liu}

\begin{document}
\allowdisplaybreaks
\maketitle

\begin{abstract}
We introduce a model concerning radiational gaseous stars and establish the existence theory of stationary solutions to the free boundary problem of hydrostatic equations describing the radiative equilibrium. We also concern the local well-posedness of the moving boundary problem of the corresponding Navier-Stokes-Fourier-Poisson system and construct a prior estimates of strong and classical solutions. Our results explore the vacuum behaviour of density and temperature near the free boundary for the equilibrium and capture such degeneracy in the evolutionary problem.
\end{abstract}

\tableofcontents

\section{Introduction}

\subsection{Problems}
The evolutionary configuration of a radiational gaseous star can be described by the following hydrodynamic system,
\begin{equation}\label{eq:dynm}
	\begin{cases}
		\rho_t + \dv (\rho \vec u) = 0   & x \in \Omega(t), \\
		(\rho \vec u )_t + \dv (\rho \vec u \otimes \vec u ) + \nabla P + \rho \nabla \psi = \iota \dv \mathbb{S} & x \in \Omega(t), \\
		(\rho e)_t + \dv(\rho e \vec u) + \theta P_\theta \dv \vec u = \Delta \theta + \epsilon \rho + \iota \mathbb S : \nabla \vec u  & x \in \Omega(t), \\
		\Delta \psi = \rho & x \in \mathbb{R}^3,
	\end{cases}
\end{equation}
where $ \rho, \vec u, P, \psi, e, \mathbb{S} $ represent the mass density, the velocity field, the pressure potential, the gravitation potential, the specific inner energy and the viscous stress tensor. The gas is assumed to be ideal and Newtonian. That is
$$ P = K \rho \theta = \theta P_\theta, ~ \mathbb{S} = \mu \bigl( \nabla \vec u + \nabla \vec{u}^\top \bigr) + \lambda \mathbb{I}_3 \dv \vec u, $$ with constants $ K > 0, \mu, \lambda > 0 $. Also, the specific inner energy is taken in the form $ e = c_\nu \theta $ with the specific heat coefficient $ c_\nu > 0 $. $ \iota = 0, 1 $ correspond to the inviscid case and the viscous case respectively. $ \epsilon > 0 $ on the right of \subeqref{eq:dynm}{3} represents the uniform rate of generation of energy. $ \Omega(t) : = \lbrace x \in \mathbb R^3 | \rho(x,t) > 0 \rbrace $ is the moving occupied domain of the gas. \eqref{eq:dynm} is complemented with the following boundary condition
\begin{equation}\label{eq:bdrcdn}
	\iota \mathbb{S}\vec{n}\bigr|_{\Gamma(t)} = 0, ~ \theta\bigr|_{\Gamma(t)} = 0, \lim\limits_{|x|\rightarrow \infty} \psi = 0, \mathcal{V}(t) = \vec u\cdot\vec{n}\bigr|_{\Gamma(t)},
\end{equation}
where $ \Gamma(t) = \partial \Omega(t) $ is the boundary of the moving domain, $ \mathcal V(t) $ represents the moving velocity of $ \Gamma(t) $ and $ \vec{n} $ is the exterior normal direction on $ \Gamma(t) $. 

To study \eqref{eq:dynm}, a fundamental question is whether there exists a steady state to such a hydrodynamic system. In fact, the equilibrium is determined by the steady state of \eqref{eq:dynm}. That is,  
\begin{equation}\label{rdmd01}
	\begin{cases}
		\nabla (K\rho\theta) = - \rho \nabla \psi & x \in \Omega,\\
		- \Delta \theta = \epsilon \rho & x \in \Omega,\\
		\Delta \psi = \rho & x \in \mathbb{R}^3,  
	\end{cases}
\end{equation}
where $ \Omega = \lbrace x\in \mathbb R^3 | \rho(x) > 0 \rbrace $. \eqref{rdmd01} is complemented with the following condition
\begin{gather}
 \int_{\Omega} \rho \,dx  =  M > 0 , ~ \lim\limits_{|x|\rightarrow \infty} \psi = 0,  \label{totalmas}\\
 \theta  =  0 ~~~ \text{on} ~ \Gamma : = \partial \Omega  \label{ztmb}.
\end{gather}

Compared with the usual zero normal heat flux condition, $ \theta $ is assumed to vanish in both the evolutionary boundary condition \eqref{eq:bdrcdn} and the steady boundary condition \eqref{ztmb}. This indicates that the temperature on the surface of a star is very low and close to zero. Roughly speaking, this is interpreted as follows. On the surface of a star, radiation transfers the heat to the space in the form of light (through emission of photons), and thus the temperature of the gas is relatively very low. Therefore, the equations \subeqref{rdmd01}{2} and \eqref{ztmb} form a balance of the heat source and the heat loss through the heat conductivity and radiation. We refer to \cite{Chandrasekhar1958} for the validity of this boundary condition. See also \cite{Hansen2004} for more discussion on the models of radiational gaseous stars. 

There are several questions arising from the problems \eqref{eq:dynm} and \eqref{rdmd01} naturally. That is,
\begin{enumerate}
	\item is there any regular solution to \eqref{rdmd01} with compact support;
	\item supposed the steady state solution exists, what is the corresponding boundary behavior of $ \rho $ and $ \theta $;
	\item if \eqref{eq:dynm} is given with the initial data $ (\rho(x,0), \vec{u}(x,0), \theta(x,0), \Omega(0)) = ( \rho_0, \vec{u}_0, \theta_0,\Omega_0 ) $ with $ (\rho_0, \theta_0) $ satisfying the boundary behavior of the steady state solution, is the system \eqref{eq:dynm} locally well-posed?
\end{enumerate}
In this work, we will answer these questions. In fact, we will establish the necessary and sufficient condition for the existence of steady states \eqref{rdmd01} with compact support in term of the rate of generation of energy $ \epsilon $, and the local well-posedness of strong and classical solutions to \eqref{eq:dynm} in the spherically symmetric motions for the viscous flows ($ \iota = 1 $).

The difficulties to solve \eqref{rdmd01} are as follows. On one hand, the temperature $ \theta $ is determined by a Poisson equation \subeqref{rdmd01}{2} with a Dirichlet boundary condition \eqref{ztmb}. Thus, to solve $ \theta$, we need the information from the distribution of the density $ \rho $ which determines the heat source and the domain of the Dirichlet problem for the Poisson equation. On the other hand, the distribution of the density $ \rho $ is determined through a free boundary problem consisting of \subeqref{rdmd01}{1} and \subeqref{rdmd01}{3} with \eqref{totalmas} which involves $ \theta $. This strong coupling makes the classical variational approach less obvious. However, despite the complex interaction of the radiation, the pressure and the self-gravity, we find out that \eqref{rdmd01} can essentially be reduced to the well-known Lane-Emden equation. Such reduction will yield the existence of steady state solutions and the corresponding boundary behaviors of the density and the temperature. Indeed, if and only if $ 1/6 < \epsilon K < 1 $, there exist infinitely many steady states (unique up to scaling) and they are spherically symmetric. On the boundary $ \partial \Omega $, 
\begin{equation}\label{PVT}
	\rho ,  \theta \bigr|_{\partial\Omega} = 0, ~-\infty < \nabla_n \rho^{\frac{\epsilon K }{1-\epsilon K}} , \nabla_n \theta \bigr|_{\partial\Omega} \leq - C < 0.
\end{equation}

Given the initial velocity $ \vec{u}_0 = u_0(r)\frac{x}{r} $ the initial density and the initial temperature $ \rho_0 = \rho_0(r), \theta_0 = \theta_0(r) $ satisfying \eqref{PVT}, where $ r = |x| $, the spherically symmetric motion $ \vec{u} = u(r,t) \frac{x}{r} $ of \eqref{eq:dynm} with $ \iota = 1 $ is described by,
\begin{equation}\label{eq:sphericalmotion}
	\begin{cases}
		\dt(\rho r^2) + \dr(\rho r^2 u) = 0 & r \in (0,R(t)) ,\\
		\dt(\rho r^2 u) + \dr(\rho r^2 u^2) + r^2 \dr(K\rho \theta) & \\
		~~ = - \rho r^2 \dr \psi + (2\mu+\lambda) r^2 \dr \bigl( \dfrac{\dr(r^2 u)}{r^2} \bigr)  & r \in (0,R(t)) ,  \\
		\dt(c_\nu \rho r^2 \theta) + \dr(c_\nu \rho r^2 \theta u) - \dr(r^2\dr\theta)  = - K\rho \theta\dr(r^2 u) &  \\ ~~ + \epsilon \rho r^2 + 2\mu r^2 \bigl( (\dr u)^2 + 2\bigl(\dfrac{u}{r} \bigr)^2 \bigr) + \lambda r^2 \bigl( \dr u + 2 \dfrac{u}{r} \bigr)^2   & r \in (0,R(t)),  \\
		\dr\psi = \frac{1}{r^2}\int_0^r s^2 \rho(s)\,ds & r \in (0,R(t)),
	\end{cases}
\end{equation}
where $ R(t) > 0 $ is the radius of the moving domain $ \Omega(t) $. \eqref{eq:sphericalmotion} is complemented with the following boundary condition
\begin{equation}{\label{symtrboundary}}
\begin{gathered}
	u(0,t)=0, ~ \bigl(2\mu r^2 \dr u + \lambda\dr(r^2u)\bigr)(R(t),t) = 0,~ \theta(R(t),t) = 0, \\
	\dt R(t) = u(R(t),t), ~ R(0) = R_0.
\end{gathered}
\end{equation}

The corresponding problem for an isentropic flow has already been  studied in \cite{Jang2010}. Just as in the classical literatures, the coordinate singularity at the centre need extra work to overcome. While in \cite{Jang2010} the author separated the entire domain into the boundary and the interior subdomains and considered the problem around the centre in Eulerian coordinates to avoid the coordinate singularity, the authors in \cite{LuoXinZeng2014,LuoXinZeng2016,LuoXinZeng2015} developed a technique to recover the regularity at the centre without imposing conditions on the derivatives of velocity in the Lagrangian coordinate. Instead, they only assume that the central velocity vanishes. Such technique was also applied in \cite{Jang2014}. As a matter of fact, the key is to study the problem in the Lagrangian coordinate induced by the flow trajectory. Such a coordinate system enables a careful track of the flow and in turn resolves the coordinate singularity. This program involves less derivative estimates in order to recover the central regularity comparing to the strategy of \cite{Jang2010}. We show in this work that 
we can as well recover the regularity of the temperature without imposing any condition on the derivatives of $ \theta $ or itself at the centre. In particular, we show that the derivative of $ \theta $ at the centre vanishes for a classical solution (see Corollary \ref{cor:e1est}). 

On the other hand, just as in \cite{Jang2010} and other related literatures, the vacuum boundary \eqref{PVT} will make the problem complicated and the heat conductivity does not help to make it easier. In the previous work \cite{LuoXinZeng2014}, by passing to the Lagrangian coordinate, the degenerate density $ \rho $ is reduced to a coefficient of the equation. Unfortunately, such a structure no longer exists for the temperature $ \theta $. We will show that for a classical solution, \eqref{PVT} can be tracked for the temperature (see Lemma \ref{lm:highregularity001}).

\subsection{Reviews of Related Works}

Before diving into our problems, we will review some related works. In the past decades, the mathematical study of gaseous stars mainly focuses on the degenerate gaseous star (\cite{Luo2008a,Luo2009b}). That is to consider the following isentropic system
\begin{equation}\label{IsenGS}
	\begin{cases}
		\nabla P(\rho) = - \rho \nabla \psi & x \in \Omega, \\
		\Delta \psi = \rho & x \in \mathbb{R}^3.
	\end{cases}
\end{equation}
In this case, the variational structure associated with \eqref{IsenGS} plays important roles. For example, in \cite{AuchmutyBeals1971,Friedman2010,Li1991}, the authors considered the rotating gaseous stars problem (that is, on the left of \subeqref{IsenGS}{1} is an additional convection) with prescribed angular momentum or angular velocity. Under some assumptions on the equation of state, the prescribed angular momentum or angular velocity, the authors have shown the existence of solutions to the rotating gaseous star problem with compact support (see also \cite{Friedman1981, Luo2008a, Luo2009b}). For a non-rotating gaseous star, the variational structure associated with \eqref{IsenGS}, together with the Hardy-Littlewood-Sobolev inequality (see \cite{Lieb2001,Lieb1987}, also known as the monotonically decreasing rearranging inequality), shows that the solution is spherically symmetric. Thus such a problem is reduced to solve an ordinary differential equation (ODE). In particular, for $ P(\rho) = \rho^\gamma$ with $ \gamma > 1 $, such solutions are called the Lane-Emden solutions for gaseous stars. Meanwhile, in \cite{Lin1997}, the author studied the linear stability and instability of the Lane-Emden solutions. Roughly speaking, for $ \gamma > 6/5 $, there is a Lane-Emden solution to \eqref{IsenGS} with compact support (see \cite{Chandrasekhar1958}). Moreover, for $ 4/3 < \gamma < 2 $, there exists a unique ball-type solution and it is neutrally stable; for $ \gamma = 4/3 $, only for a special total mass $ M^* > 0 $ (known as the critical mass for radiation gaseous stars \cite{Fu1998}) there are infinitely many ball-type solutions, and they are marginally stable; for $ 6/5 < \gamma < 4/3 $, it admits ball-type, singular ball-type and singular ground state solutions and the ball-type solution is unique and unstable; for $ \gamma = 6/5 $, there is only a ground state solution with finite total mass; for $ 1 < \gamma < 6/5 $, no solution with finite total mass exists. For the ball-type solutions, the physical vacuum boundary holds, i.e.
\begin{equation}\label{PV}
	\rho\bigr|_{x\in \partial\Omega} = 0, ~ - \infty < \nabla_{\vec n} P'(\rho) \bigr|_{x\in \partial\Omega} \leq -C < 0,
\end{equation}
for some constant $ 0 < C < \infty $, where $ \vec n $ is the exterior normal vector on the surface $ \partial \Omega $ and $ c^2 := P'(\rho) = \frac{d}{d \rho} P(\rho) $ represents the square of sound speed. 

It shall be emphasised that the physical vacuum \eqref{PV} causes big challenges in the study of the evolutionary problem of \eqref{IsenGS}. That is to study the following Euler-Poisson system
\begin{equation}\label{EPE}
	\begin{cases}
		\dt \rho + \dv (\rho \vec u ) = 0 & x \in \Omega(t),\\
		\rho \dt \vec u + \rho \vec u \cdot \nabla \vec u + \nabla P = - \rho \nabla\psi & x \in \Omega(t),\\
		\Delta \psi = \rho & x \in \mathbb{R}^3.
	\end{cases}
\end{equation}
Indeed, the physical vacuum boundary implies that the sound speed $ c $ is only $ 1/2 $-H\"older continuous rather than Lipschitz continuous across the boundary. As pointed out in \cite{Liu1996}, this singularity makes the standard hyperbolic method fail. It is only recently that the local well-posedness theory is studied by Coutand, Lindblad and Shkoller \cite{Coutand2010,Coutand2011a,Coutand2012}, Jang and Masmoudi \cite{Jang2009b,Jang2015}, Gu and Lei \cite{Gu2012,Gu2015}, Luo, Xin and Zeng \cite{LuoXinZeng2014} in the setting of one spatial dimension, three spatial dimension and spherical symmetry with or without self-gravitation. See \cite{Jang2010} for a viscous flow. We refer to \cite{Liu1996,ZengHH2015a,ZengHH2015,ZengHH2014} for other discussions on the physical vacuum problem and the references therein for other vacuum problems. 

As for the global dynamic property of \eqref{EPE}, the nonlinear stability and instability of the equilibrium is partially understood only for spherically symmetric motions. In particular, the works from Jang and Tice \cite{Jang2008a,Jang2014,Jang2013a} show that for $ 6/5 \leq \gamma < 4/3 $, the Lane-Emden solutions are unstable, and additional viscosities can not reduce such instability. 
When it comes to the case $ 4/3 < \gamma < 2 $, the asymptotic stability theory is first studied in \cite{LuoXinZeng2016} in the viscous case. See \cite{LuoXinZeng2015} for the model with a degenerate viscosity. Meanwhile, the case $ \gamma = 4/3 $, which is also referred to as the radiation case (see \cite{Chandrasekhar1958,Hadzic2016}), is rather complicated. On one hand, only for a specific total mass $ M^*>0$, called the critical mass, the problem \eqref{IsenGS} admits solutions with compact support. On the other hand, by studying a class of self-similar solutions in \cite{Fu1998}, the authors showed the existence of expanding and collapsing solutions with variant total mass. Such phenomena are far from fully understood. Recently, in \cite{Hadzic2016}, the authors show the asymptotic stability of the expanding solutions
. Their analysis contains the exploration of a damping structure in the Lagrangian coordinate. For other discussions on the stability, property of rotating or other models of gaseous stars, see \cite{Deng2002,Rein2003a,Makino2015,Li1991,Chanillo1994,Luo2004,Luo2009b,Luo2008a,Luo2014,Wu2013,Chanillo2012}. 

\subsection{Lagrangian Formulation, Methodology and Results}

%
%
%
%

I will summarize the key observation and technique in this work. In order to solve the steady state problem \eqref{rdmd01}, the key observation is that from $ \eqref{rdmd01}_1 $, 
\begin{equation*}
	K \dfrac{\theta}{\rho} \nabla \rho + K \nabla \theta = - \nabla \psi,
\end{equation*}
where $ \rho, \theta $ both are positive interiorly and vanish on the boundary, $ \nabla \psi $ must be parallel to the normal direction of the surface $ \partial \Omega $. In fact, with some regularity assumption on the regular solution, I will show that for any regular solution, such a parallel relation holds. Moreover, it implies that the gravitation potential $ \psi $ is a constant along the boundary $ \partial \Omega $. Therefore, a Pohozaev type argument will imply that $ \epsilon \nabla \psi = - \nabla \theta $ inside $ \Omega $. 
Such a property together with $ \eqref{rdmd01}_1 $ yields that the temperature is in the form $ \theta = \theta(\rho) $. Hence, after substituting this relation into \eqref{rdmd01}, the problem is reduced to the classical isentropic gaseous star system with $ P = P(\rho) $. Therefore, we can achieve the existence theory and the spherically symmetric property of the solutions after applying the classical existence theory. 
In fact, the solution will satisfy \eqref{PVT} and 
\begin{equation}\label{intro:exb}
\begin{gathered}
	\abs{\dr^k \rho}{} \leq O(\rho^{\frac{1-(1+k)\epsilon K}{1-\epsilon K}}),  ~ 
	\abs{\dr^k \theta}{} \leq O(1) + O(\theta^{\frac{1-\epsilon K}{\epsilon K} - k + 2}),
\end{gathered}
\end{equation}
near the boundary, where $ k \in \mathbb Z^+ $ and $ \dr $ represents the derivative along the radial direction.

As for the evolutionary problem \eqref{eq:sphericalmotion}, we shall work in the Lagrangian coordinate.
In other words, denote $ r = r(x,t) $ defined by $ \frac{d r}{dt} = u(r(x,t),t) $ and $ r(x,0) = x $, where $ x \in (0, R_0) $. Then \eqref{eq:sphericalmotion} can be written as,
\begin{equation}\label{eq:LagrangianCoordinates}
	\begin{cases}
		\bigl(\dfrac{x}{r}\bigr)^2 \rho_0 \dt v + \bigl(K\dfrac{x^2\rho_0}{r^2r_x} \Theta \bigr)_x = - \dfrac{x^2\rho_0}{r^4} \int_0^x s^2 \rho_0(s)\,ds & \\ ~~~~~~~~~~ + (2\mu+\lambda)  \bigl( \dfrac{(r^2 v)_x}{r^2 r_x}\bigr)_x  & x \in (0, R_0),  \\
		c_\nu x^2 \rho_0 \dt \Theta + K \dfrac{x^2 \rho_0}{r^2r_x} \Theta (r^2 v)_x - \bigl(\dfrac{r^2}{r_x} \Theta_x \bigr)_x = \epsilon x^2 \rho_0 & \\ ~~~~~~~~~~~ + 2\mu r^2 r_x \bigl( \bigl(\dfrac{v_x}{r_x}\bigr)^2 + 2 \bigl( \dfrac{v}{r}\bigr)^2 \bigr) + \lambda r^2 r_x \bigl( \dfrac{v_x}{r_x} + 2 \dfrac{v}{r} \bigr)^2 & x \in (0, R_0),
	\end{cases}
\end{equation} 
where the Lagrangian unknowns are defined as
\begin{equation}
	v(x,t) = u(r(x,t),t) = \dt r(x,t), ~ \Theta(x,t) = \theta (r(x,t),t).
\end{equation}
\eqref{eq:LagrangianCoordinates} is complemented with the initial and boundary conditions
\begin{gather}
	\Theta(x,0) = \Theta_0(x) = \theta(x,0),~ v(x,0) = v_0(x) = u(x,0), ~ x \in (0,R_0), {\label{initialcndtnLG}}\\
	v(0,t)= 0, ~ \Theta(R_0, t) = 0, ~ \biggl\lbrack (2\mu+\lambda)\dfrac{v_x}{r_x} + 2 \lambda \dfrac{v}{r}\biggr\rbrack (R_0,t) = 0, ~ t \geq 0, \label{boundarycndtnLG}
\end{gather}
where $ \rho_0 > 0, \Theta_0 > 0 $ for $ x \in (0,R_0) $ and are imposed with the properties \eqref{PVT} and \eqref{intro:exb}. That is, for $ x \in (R_0/2,R_0] $
\begin{equation}\label{PV:lagrangian}
	\begin{gathered}
		\rho_0(R_0), \Theta_0(R_0) = 0 , ~ -\infty < \bigl( \rho_0^{\frac{\epsilon K}{1-\epsilon K}} \bigr)_x , \bigl( \Theta_0 \bigr)_x \leq -C < 0, \\
		\abs{\dx^k \rho_0}{} \leq O(\rho_0^{\frac{1-(1+k)\epsilon K}{1-\epsilon K}}), ~ \abs{\dx^k \Theta_0}{} \leq O(1) + O(\Theta_0^{\frac{1-\epsilon K}{\epsilon K} - k + 2}),
	\end{gathered}
\end{equation}
where $ k\in \mathbb{Z}^+ $.

To study the well-posedness of the problem \eqref{eq:LagrangianCoordinates}, which will yield the well-posedness of \eqref{eq:sphericalmotion}, we will focus on resolving the coordinate singularity. In the energy estimates, we use the Hardy's inequality and the rescaled Poincar\'e inequality to manipulate the coordinate singularity at the center($ x = 0 $) and the vacuum singularity on the boundary($ x = R_0 $, see Lemma \ref{lm:hardy}). Compared with the isentropic case in \cite{LuoXinZeng2016}, we need the estimate of the pressure potential, or equivalently, the estimate of $ \Theta $, in order to recover the spatial regularity of $ r, v $ at the coordinate center. This can be done by calculating the $ L^\infty_t L^2_x $ norm of
\begin{equation*}
	\dfrac{1}{r\sqrt{r_x}} \bigl( \dfrac{r^2}{r_x} \Theta_x\bigr)_x = \dfrac{\text{the rest terms in \subeqref{eq:LagrangianCoordinates}{2}}}{r\sqrt{r_x}}.
\end{equation*}
The denominator $ r\sqrt{r_x} $ here is specially designed to avoid extra remaining nonlinear term. This yields the $ L^\infty_tL^2_x $ estimates of $ \Theta_x $ and $ x \Theta_{xx} $. On the other hand, the quantity $ \mathcal G : = \log \frac{r^2r_x}{x^2} $ satisfies an equation of the form
\begin{equation*}
	\mathcal G_{xt} = f \mathcal G_{x} + g
\end{equation*}
for some functions $ f, g $. Thus by applying Gr\"onwall's inequality, it implies the bound of $ L^\infty_tL^2_x $ norms of $ \mathcal G_{x} $ and $\mathcal G_{xt}$. This will yield the bound of $ L^\infty_t L^2_x $ norm of $ r_{xx}, v_{xx}, \bigl( \frac{r}{x} \bigr)_x $ and $ \bigl(\frac{v}{x}\bigr)_x $. Such a structure is first used in \cite{LuoXinZeng2015, LuoXinZeng2016}. These steps will resolve the coordinate singularity. 

Besides, we obtain some new point-wise estimates without using the embedding theory, which will dramatically reduce the nonlinearity of the system and may be useful in future study. By applying the multiplier $ r^3 $ to \subeqref{eq:LagrangianCoordinates}{1}, we will be able to obtain the point-wise estimates of $ v_x $ and $ v/x $. In addition, after integrating \subeqref{eq:LagrangianCoordinates}{2} from the centre, we will obtain the point-wise bound of $ x^{1/2}\Theta_x $ (see Lemma \ref{lm:pointwise}). For a classical solution (higher order regularity), similar point-wise estimate will yield $ \Theta_x(0,t) = 0 $ (see Corollary \ref{cor:e1est}). Also, $ \Theta_{xt} $ is point-wisely bounded near the boundary for such a solution. Then the fundamental theory of calculus shows that $ - \infty < \Theta_x < - C $ (and $ C < \Theta / \sigma < \infty $) with some constant $ C>0 $ for a short time (see Lemma \ref{lm:highregularity001}). Thus the vacuum property \subeqref{PV:lagrangian}{1}(and hence \eqref{PVT}) is recovered. 

The rest of the paper will be organized as follows. In the next section, we will introduce the main theorems and some notations used in this work. In Section \ref{sec:steadysol}, we show the existence theory for the steady state problem \eqref{rdmd01} and prove Theorem \ref{thm:steadystatesol}. In Section \ref{sec:aprior}, we show the a prior estimates for the evolutionary problem \eqref{eq:LagrangianCoordinates}. In particular, the point-wise estimates are given in Section \ref{sec:pointwise}. The energy estimates for the strong solutions are given in Section \ref{sec:energy}. Section \ref{sec:elliptic} is devoted to the regularity estimates for the strong solution. In Section \ref{sec:classical}, we present the corresponding estimates for the classical solutions and it will finish the proof of Theorem \ref{thm:evolutionaryproblem-aprioriest}. In Section \ref{sec:well-posedness}, we employ the fixed point theory to derive the well-posedness theory (i.e. Theorem \ref{thm:evolutionaryproblem-solution}).

\section{Main Theorems and Preliminaries}

During this work, the following notations are adopted. For any quantities $ A, B $, by $ A \lesssim B $ it means that
\begin{equation}
	\exists C > 0, ~ A \leq C\cdot B,
\end{equation}
where $ C $ is a generic constant and will be different from lines to lines but independent of the solutions. Also, $ A \simeq B $ means $ A \lesssim B $ and $ B \lesssim A $. 
By $ C = C(\cdot) $, it means a positive constant $ C $ depends solely on the inputs. Similarly, by $ P = P(\cdot) $, it means a positive polynomial of the inputs which is increasing in each input.

A regular solution to the steady state problem \eqref{rdmd01} is defined as follows.
\begin{df}\label{def:regularsol}
	We call the triple $ \left(\rho, \theta, \Omega \right) $ a regular solution to the free boundary problem \eqref{rdmd01} with complementing conditions \eqref{totalmas} and \eqref{ztmb} if it satisfies,
	\begin{itemize}
		\item $ \Omega \subset \mathbb{R}^3 $ is a bounded, stared-shaped domain, and the origin $ O = (0,0,0) $ lies inside $ \Omega $;
		\item the boundary of $ \Omega $ is a $ C^2 $ closed surface in $ \mathbb{R}^3 $;
		\item $ \rho, ~ \theta $ vanish on the boundary of $ \Omega $, and are positive inside $ \Omega $;
		\item $ \theta \in C^2(\bar{\Omega}) $, and $ \rho \in C^1 (\Omega) \cap C(\bar \Omega) $.
	\end{itemize}
\end{df}
To consider the evolutionary problem \eqref{eq:LagrangianCoordinates}, we will work in the space-time $ (x,t) \in \lbrack 0, R_0 \rbrack \times \lbrack 0, T \rbrack $ for some constant $ T > 0 $ which will be determined later. Also, the Sobolev space $ X_tY_x $ is defined as
\begin{equation}
	X_tY_x : = X( (0, T), Y(0,R_0)),
\end{equation}
where $ X, Y $ represent some Sobolev spaces in the temporal variable and the spatial variable. The interior cut-off function $ \chi: \lbrack 0 , R_0 \rbrack \rightarrow \lbrack 0,1 \rbrack $ is defined as a smooth function satisfying 
\begin{equation}
	\chi = \begin{cases}
		1 & x\in \lbrack 0, R_0/4\rbrack, \\
		0 & x \in \lbrack R_0/2, R_0\rbrack,
	\end{cases}
\end{equation}
and $ - 8/R_0 \leq \chi' \leq 0 $.

Using these notation, a strong solution to \eqref{eq:LagrangianCoordinates} is a solution $ (r, v, \Theta) $ satisfying the following.
\begin{df}\label{def:strongsol} Given $ 0 < T <\infty $, 
	the energy functional for a strong solution to \eqref{eq:LagrangianCoordinates} is defined as
	\begin{equation}
	\begin{aligned}
		\mathcal E_1 = & \mathcal{E}_1(v,\Theta) := \sum_{i=0}^{1} \bigl( \norm{x\sqrt{\rho_0}\dt^i v}{\stnorm{\infty}{2}}^2 + \norm{x\sqrt{\rho_0}\dt^i \Theta}{L_t^\infty L_x^2}^2 \\
		&  + \norm{\sqrt{\chi\rho_0}\dt^i v}{\stnorm{\infty}{2}}^2 \bigr)  + \norm{xv_x}{\stnorm{\infty}{2}}^2 + \norm{v}{\stnorm{\infty}{2}}^2 + \norm{x\Theta_x}{\stnorm{\infty}{2}}^2 \\
		& + \norm{\sqrt{\chi}v_x}{\stnorm{\infty}{2}}^2  + \norm{\sqrt{\chi}\frac{v}{x}}{\stnorm{\infty}{2}}^2 + \sum_{i=0}^1 \bigl( \norm{x\dt^i v_x}{\stnorm{2}{2}}^2 + \norm{\dt^i v }{\stnorm{2}{2}}^2 \\
		& + \norm{x\dt^i \Theta_x}{\stnorm{2}{2}}^2 + \norm{\sqrt{\chi}\dt^i v_x}{\stnorm{2}{2}}^2 + \norm{\sqrt{\chi}\frac{\dt^iv}{x}}{\stnorm{2}{2}}^2 \bigr)  \\
		& + \norm{x\sqrt{\rho_0}v_t}{\stnorm{2}{2}}^2 + \norm{x\sqrt{\rho_0}\Theta_t}{\stnorm{2}{2}}^2 + \norm{\sqrt{\chi\rho_0}v_t}{\stnorm{2}{2}}^2 .
	\end{aligned}
	\end{equation}
	A strong solution is a triple $ (r, v, \Theta) = (r(x,t), v(x,t), \Theta(x,t)) $ satisfying the system \eqref{eq:LagrangianCoordinates} with \eqref{initialcndtnLG}, \eqref{boundarycndtnLG} and has the following regularity
	\begin{equation}\label{Regularity:Strong}
		\begin{gathered}
			r, v\in L_t^{\infty}H_x^{2}\cap L_t^2H_x^1 , ~ \Theta \in L_t^\infty H_x^1 , ~ x\Theta \in L_t^\infty H_x^2\cap L_t^2 H_x^1,\\ 
			v_t \in L_t^2H_x^1, ~ x\Theta_t \in L_t^2H_x^1
		\end{gathered}
	\end{equation}
	Also, the strong solution is in the functional space 
	\begin{equation}\label{functionalspace:strong}
	\begin{aligned}
		\mathfrak X = & \mathfrak X_T : = \biggl\lbrace (r,v,\Theta)\bigl| \mathcal E_1 < \infty ~ \text{and} ~  \norm{r_{xx}}{\stnorm{\infty}{2}}, \norm{\bigl(\dfrac{r}{x}\bigr)_x}{\stnorm{\infty}{2}}, \\ 
		& \norm{v_{xx}}{\stnorm{\infty}{2}},\norm{\bigl(\dfrac{v}{x}\bigr)_x}{\stnorm{\infty}{2}},
			\norm{\Theta_x}{\stnorm{\infty}{2}}, \norm{x\Theta_{xx}}{\stnorm{\infty}{2}} < \infty 
		\biggr\rbrace.
	\end{aligned}
	\end{equation}
\end{df}

Similarly, a classical solution to \eqref{eq:LagrangianCoordinates} is defined as follows.
\begin{df}\label{def:classicalsol}
Given $ 0 < T <\infty $, 
	the energy functional for a classical solution to \eqref{eq:LagrangianCoordinates} is defined as
	\begin{equation}
	\begin{aligned}
		\mathcal E_2 = & \mathcal{E}_2 ( v, \Theta) : = \mathcal{E}_1 +\bigl( \norm{x\sqrt{\rho_0} v_{tt}}{\stnorm{\infty}{2}}^2 + \norm{x\sqrt{\rho_0} \Theta_{tt}}{L_t^\infty L_x^2}^2 \\
		& + \norm{\sqrt{\chi\rho_0} v_{tt}}{\stnorm{\infty}{2}}^2 \bigr)   + \norm{xv_{xt}}{\stnorm{\infty}{2}}^2 + \norm{v_t}{\stnorm{\infty}{2}}^2 + \norm{x\Theta_{xt}}{\stnorm{\infty}{2}}^2 \\
		&   + \norm{\sqrt{\chi}v_{xt}}{\stnorm{\infty}{2}}^2 + \norm{\sqrt{\chi}\frac{v_t}{x}}{\stnorm{\infty}{2}}^2 +  \bigl( \norm{x v_{xtt}}{\stnorm{2}{2}}^2 + \norm{ v_{tt} }{\stnorm{2}{2}}^2 \\
		& + \norm{x \Theta_{xtt}}{\stnorm{2}{2}}^2 + \norm{\sqrt{\chi} v_{xtt}}{\stnorm{2}{2}}^2 + \norm{\sqrt{\chi}\frac{v_{tt}}{x}}{\stnorm{2}{2}}^2 \bigr)\\
		& + \norm{x\sqrt{\rho_0}v_{tt}}{\stnorm{2}{2}}^2  + \norm{x\sqrt{\rho_0}\Theta_{tt}}{\stnorm{2}{2}}^2 + \norm{\sqrt{\chi\rho_0}v_{tt}}{\stnorm{2}{2}}^2.
	\end{aligned}
	\end{equation}
	Then the classical solution is a triple $ (r, v, \Theta) = (r(x,t), v(x,t), \Theta(x,t)) $ in $ \mathfrak X $ satisfying \eqref{eq:LagrangianCoordinates} and has the following regularity
	\begin{equation}\label{Regularity:Classical}
		\begin{gathered}
			r, v\in L_t^{\infty}H_x^{3} \cap L_t^2 H_x^1 , ~ \Theta \in L_t^\infty H_x^2 , ~ x\Theta \in L_t^\infty H_x^3 \cap L_t^2 H_x^1, \\
			v_{t} \in L_t^\infty H_x^2 \cap L_t^2H_x^1 , ~ \Theta_t \in L_t^\infty H_x^1 , ~ x\Theta_t \in L_t^\infty H_x^2 \cap L_t^2 H_x^1, \\
			 v_{tt} \in L_t^2H_x^1, ~ x\Theta_{tt} \in L_t^2 H_x^1 .
		\end{gathered}
	\end{equation}
	Also, the strong solution is in the functional space 
	\begin{equation}\label{functionalspace:classical}\begin{aligned}
	\mathfrak Y = & \mathfrak Y_T : = \biggl\lbrace (r,v,\Theta)\in \mathfrak X \bigl| \mathcal E_2 < \infty ~ \text{and} ~  \norm{\Theta_{xt}}{\stnorm{\infty}{2}}, \norm{x\Theta_{xxt}}{\stnorm{\infty}{2}},\\
	 & \norm{v_{xxt}}{\stnorm{\infty}{2}},  \norm{\bigl(\dfrac{v_t}{x}\bigr)_x}{\stnorm{\infty}{2}} , \norm{r_{xxx}}{\stnorm{\infty}{2}}, \norm{\bigl(\dfrac{r}{x}\bigr)_{xx}}{\stnorm{\infty}{2}}, \\ & \norm{v_{xxx}}{\stnorm{\infty}{2}},
	 \norm{\bigl(\dfrac{v}{x}\bigr)_{xx}}{\stnorm{\infty}{2}},
	   \norm{\Theta_{xx}}{\stnorm{\infty}{2}}, \norm{x\Theta_{xxx}}{\stnorm{\infty}{2}} < \infty
		\biggr\rbrace .
	\end{aligned}	
	\end{equation}	
\end{df}
Also, we denote the following initial data for the energy functionals $ \mathcal E_1 $ and $ \mathcal E_2 $,
\begin{align}
	\mathcal E_1^0 = & \mathcal E_1^0 (v_0,\Theta_0) :=  \sum_{i=0}^{1} \bigl( \norm{x\sqrt{\rho_0}\dt^i v_0}{L^2}^2 + \norm{x\sqrt{\rho_0}\dt^i \Theta_0}{L^2}^2 {\nonumber} \\
	& + \norm{\sqrt{\chi\rho_0}\dt^i v_0}{L^2}^2 \bigr)   + \norm{xv_{0,x}}{L^2}^2 + \norm{v_0}{L^2}^2 + \norm{x\Theta_{0,x}}{L^2}^2 {\label{initialenergy-strong}} \\
	&  + \norm{\sqrt{\chi}v_{0,x}}{L^2}^2  + \norm{\sqrt{\chi}\frac{v_0}{x}}{L^2}^2, {\nonumber} \\
	\mathcal E_2^0 = & \mathcal E_2^0 (v_0,\Theta_0) :=  \mathcal E_1^0 + \bigl( \norm{x\sqrt{\rho_0}v_{0,tt}}{L^2}^2 + \norm{x\sqrt{\rho_0} \Theta_{0,tt}}{L^2}^2{\nonumber} \\
	&  + \norm{\sqrt{\chi\rho_0} v_{0,tt}}{L^2}^2 \bigr)  + \norm{xv_{0,xt}}{L^2}^2 + \norm{v_{0,t}}{L^2}^2 + \norm{x\Theta_{0,xt}}{L^2}^2 {\label{initialenergy-classical}} \\
	&  + \norm{\sqrt{\chi}v_{0,xt}}{L^2}^2 + \norm{\sqrt{\chi}\frac{v_{0,t}}{x}}{L^2}^2, {\nonumber}
\end{align}
where $ L^2 = L^2_x(0,R_0) $ and $ v_{0,t}, v_{0,tt}, \Theta_{0,t}, \Theta_{0,tt} $ are defined by the equation \eqref{eq:Lg000} and \eqref{eq:Lg010}. 

Now we are able to write down the main theorems in this work. For the steady state problem \eqref{rdmd01}, we have,
\begin{thm}\label{thm:steadystatesol}
	Given $ M > 0 $ in \eqref{totalmas}, if and only if $ 1/6 < \epsilon K < 1 $, there are infinitely many (self-similar) regular solutions to \eqref{rdmd01} as defined in Definition \ref{def:regularsol}. In addition, the solutions are spherically symmetric and satisfy \eqref{PVT}, \eqref{intro:exb} on the gas-vacuum interface. 
\end{thm}
For the evolutionary problem \eqref{eq:LagrangianCoordinates}, we have,
\begin{thm}\label{thm:evolutionaryproblem-aprioriest}
	Given initial data $ (r(x,0),v(x,0),\Theta(x,0)) = (x,v_0(x),\Theta_0(x)) $ for $ x \in (0, R_0) $, the smooth solution to \eqref{eq:LagrangianCoordinates} satisfies the following. 
	\begin{enumerate}
	\item If $ \mathcal E_1^0 < \infty $, there is a positive time $ T_* $ such that for $ T \leq T_* $, $ (r,v,\Theta) \in \mathfrak X $, where $ T_* \geq 1/ P_*(\mathcal E_1^0 + 1) $, for some positive polynomial $ P_* = P_*(\cdot) $;
	\item if $ \mathcal E_1^0 < \mathcal E_2^0 < \infty $, there is a positive time $ T_{**} $ such that for $ T \leq T_{**} $, $ (r,v,\Theta) \in \mathfrak Y $, where $ T_{**} \geq 1/P_{**}(\mathcal E_1^0+1,\mathcal E_2^0 + 1) $, for some positive polynomial $ P_{**} = P_{**}(\cdot) $.
	\end{enumerate}
\end{thm}
Also, we have the following existence theory of \eqref{eq:LagrangianCoordinates}(or equivalently \eqref{eq:dynm}).
\begin{thm}\label{thm:evolutionaryproblem-solution}
	Given initial data $ (r(x,0),v(x,0),\Theta(x,0)) = (x,v_0(x),\Theta_0(x)) $ for $ x \in (0, R_0) $, if $ \mathcal E_1^0 < \infty $, there is a constant $ 0 < \tilde T_* < \infty $ such that there is a strong solution (defined in Definition \ref{def:strongsol}) to \eqref{eq:LagrangianCoordinates} for $ 0 < T < \tilde T_* $. Also, if the initial data $ (r(x,0),v(x,0),\Theta(x,0)) = (x,v_0(x),\Theta_0(x)) $ satisfies $ \mathcal E_2^0 < \infty $, there is a constant  $ 0 < \tilde T_{**} < \infty $  such that there is a classical solution  (defined in Definition \ref{def:classicalsol}) to \eqref{eq:LagrangianCoordinates} for $ 0 < T < \tilde T_{**} $.
\end{thm}

In this work, we will make use of the following lemmas. The first lemma shows that the gradient of a function with positive value interiorly and zero value on the boundary is parallel to the exterior normal direction on the boundary. 

\begin{lm}[Parallel Relation]\label{lm:pllrltn}
	Consider the domain $ \Omega^+ = \mathbb T^2 \times \mathbb R^+ =  (-1,1)^2 \times (0, + \infty) $ and a function $ f $ satisfies:
	\begin{gather*}
		f \in C^1 ( \Omega^+ ) \text{ and } f > 0 \text{ in } \Omega^+; ~ f = 0 ~ \text{on} ~ \partial \Omega^+ = \left\lbrace{(x_1,x_2,x_3)|x_3 = 0}\right\rbrace, \\
		V  ~ \text{is uniformly continuous and bounded up to the boundary,}
	\end{gather*}
	where
	\begin{displaymath}
		V = (V_1,V_2,V_3)^\top:= x_3 \dfrac{\nabla f}{f}.
	\end{displaymath}
	Then $ V $ can be continuously extended to the boundary $\partial \Omega^+$ and is also denoted as $ V $. In addition, we have $ V \parallel (0,0,1)^\top $ and $ V \cdot (0,0,1)^\top \geq 0 $ on $ \partial \Omega^+ $.
\end{lm}
\begin{pf}
	After shifting, it suffices to show, in the sense of trace,
	\begin{align*}
		V_3(0,0,0) \geq 0, ~ \text{and} ~ V_i(0,0,0) = 0, ~ i=1,2.
	\end{align*}
	First, $ V_3 (0,0,0) \geq 0 $ is easy to verify. Suppose it was not true, i.e. $ V_3 (0,0,0) < 0 $. Then there would be a constant $\delta > 0 $ such $ V_3 (0,0,x_3) < 0 $ for $ x_3 \in (0,\delta) $, which would imply, $ \partial_{x_3} f (0,0,x_3) < 0 $. In particular, $ 0 < f(0,0,\delta/2) < f(0,0,0) = 0 $ which is impossible. 
	
	Next, we show $ V_1(0,0,0), V_2(0,0,0) = 0 $ through a contradictory argument. Otherwise, without loss of generality, suppose $ V_2(0,0,0) \geq c > 0 $ for some positive constant $ c $. In the following, we work in the plane $$ S : = \left\lbrace (x_1,x_2,x_3) \in \Omega^+| x_1 = 0 \right\rbrace. $$ We consider two cases, $ V_3(0,0,0) =0 $ and $ V_3(0,0,0) > 0 $. 
	
	If $ V_3(0,0,0) = 0 $, we have for some $ 0 < \epsilon <<1 $ and constant $ c > 0 $, there is a constant $ \delta >0 $ such that the following holds for $ (x_1,x_2,x_3) \in B_\delta(0) \cap \Omega^+ $,
	\begin{equation}\label{geomequ0000}
		\left|V_3(x_1,x_2,x_3)\right| < \epsilon, ~ V_2(x_1,x_2,x_3) > c-\epsilon > 0.  
	\end{equation}
	This will imply $ \partial_{x_2} f(0, x_2, x_3) > 0 $ in $B_\delta(0) \cap S $. Then $ f $ is strictly increasing in the $ x_2 $-variable inside such neighbourhood, and by the implicit function theorem, through any point $ p \in B_\delta(0) \cap S $, there is a neighbourhood of $ p $ such that the level set of $ f $ passing through $ p $ is well-defined. This will suggest that any level set of $ f $ does not have end-points inside $ B_\delta(0) $. Let $ \delta^* \leq \delta/2 $. Consider the level set of $ f $, $ \left\lbrace f = f(0,0,\delta^*)> 0 \right\rbrace $. Then the intersection of such level set and the plane $ S $ is a curve, parametrised by $$ \left\lbrace \tau(t) = (0,x_2(t),x_3(t)), t \in (-\omega_1,+\omega_2) \right\rbrace $$ for some $ \omega_1,\omega_2 > 0 $. Moreover, since $ \tau $ has no end-point inside $ B_\delta(0) \cap S $, it satisfies $ \left| \tau'(t) \right| > 0$ , $ \tau(0) = (0,0,\delta^*) $ and
	\begin{equation}\label{geomequ0001}
		x_2'(t) \cdot V_2(\tau(t)) + x_3'(t) \cdot V_3(\tau(t)) = 0,
	\end{equation}
	as long as $ \tau $ is still inside such neighbourhood, which follows from the fact $\tau $ has no end-points inside this area and the fact $ V $ is pointing toward the normal direction of $ \tau $. First, $ x_3'(0) \neq 0 $. Otherwise, from \eqref{geomequ0000} and \eqref{geomequ0001}, $ x_2'(0) = 0 $, which is impossible. Without lost of generality, let's assume $ x_3'(0) < 0 $. Denote $ T = \sup \left\lbrace t | \tau(s) \in B_\delta(0)~ \text{for} ~ 0 < s < t \right\rbrace $, called the escaping time of the curve. Due to the monotonicity of $ f $ in $ x_2 $-variable, the curve $ \tau $ will eventually extend outside $ B_\delta(0) $ through the boundary $ \partial B_\delta(0) \cap \left\lbrace x_3 > 0 \right\rbrace $, since otherwise $ \tau $ would be a closed curve and this would be contradictory to the monotonicity of $f$. In the meantime, we claim $ x_3'(t) < 0 $ for $ 0 < t < T $. Otherwise, there was a $ t^* \in (0,T ) $ such that $ x_3'(t^*)=0 $. Together with \eqref{geomequ0000} and \eqref{geomequ0001}, it would imply $ x_2'(t^*) =0 $  which is impossible. Therefore, the curve $ \tau $ will pass the boundary section $ \partial B_\delta(0) \cap \left\lbrace 0 < x_3 < \delta /2 \right\rbrace $. It follows $ \left| x_2(T) \right| \geq \delta /2 $. On the other hand, from \eqref{geomequ0001} $x_2'(t) = 0 $ if and only if $ V_3(\tau(t)) = 0 $. Hence $\lbrace t\in(0,T)|V_3(\tau(t))\neq 0 \rbrace = \lbrace t \in (0,T) | x'_2(t) \neq 0 \rbrace $. Direct calculation shows
	\begin{align*}
			& x_3(T) = x_3(0) - \int_0^T  \left| x_3'(t) \right| \,dt \leq \delta^* - \int_{ \lbrace t \in (0,T) | V_3(\tau(t)) \neq 0 \rbrace} \left| \dfrac{V_2(\tau(t))}{V_3(\tau(t))} x_2'(t) \right| \,dt \\
			& ~~~~ \leq \delta^* - \dfrac{c-\epsilon}{\epsilon}\int_{ \lbrace t \in (0,T) | x_2'(t) \neq 0 \rbrace}  \left| x_2'(t)\right| \,dt \leq \delta^* - \dfrac{c-\epsilon}{\epsilon} \left| \int_0^T x_2'(t) \,dt \right|\\
			& ~~~~ = \delta^* - \dfrac{c-\epsilon}{\epsilon} \left| x_2(T)\right| \leq \delta^* - \dfrac{c-\epsilon}{\epsilon} \delta/2 < 0,
	\end{align*}
	for $ \delta^* $ small enough. This would mean the level set $ \left\lbrace f= f(0,0,\delta^*) > 0 \right\rbrace $ intersects with the boundary $ \partial \Omega^+ $ which is impossible.
	
	If $ V_3(0,0,0) > 0 $, there is a constant $\delta > 0 $ such that the following holds. For $ (x_1,x_2,x_3) \in B_\delta(0)\cap\Omega^+ $ with some $ 0 < c < C < \infty $,
	\begin{equation}\label{geomequ0002}
		c < V_2(x_1,x_2,x_3), V_3(x_1,x_2,x_3) < C.
	\end{equation}
	Similar as before, $ f $ is strictly increasing in $ x_2,x_3 $-variables inside such neighbourhood and any level set of $ f $ will not have end-points inside $ B_\delta(0) $. Let $ \delta^* < \delta/2 $. Consider the level set of $ f $, $ \left\lbrace f= f(0,0,\delta^*)> 0 \right\rbrace $ . Then the intersection of such level set and the plane $ S $ is a curve, parametrised by $$ \left\lbrace \tau(t) = (0,x_2(t),x_3(t)), t \in (-\omega_1,+\omega_2)\right\rbrace $$ for some $\omega_1, \omega_2 > 0 $. Moreover, $ \tau $ satisfies \eqref{geomequ0001} and $ \tau(0) = (0,0,\delta^*) $, $\abs{\tau'(t)}{} > 0 $ as long as $ \tau $ is inside $ B_\delta(0) $. Then by applying the same argument as before, $ x_3'(0) \neq 0 $. Without loss of generality, we assume $ x_3'(0) < 0 $. Then, $ x_2'(0) > 0 $. In particular, $ x_2'(t) > 0 $ as long as $ x_3'(t) < 0 $. Denote $ T = \sup \left\lbrace t | \tau(s) \in B_\delta(0)~ \text{for} ~ 0 < s < t \right\rbrace $. Then $ x_2'(t) > 0 $ and $ x_3'(t) < 0 $ for $ 0 < t < T $. Therefore, the curve $ \tau $ will pass the boundary section $ \partial B_\delta(0) \cap \left\lbrace 0 < x_3 < \delta/2, x_2 > 0 \right\rbrace $. Therefore, $ x_2(T) \geq \delta /2 $. Direct calculation shows
	\begin{align*}
			& x_3( T ) = x_3(0) + \int_0^T x_3'(t) \,dt = \delta^* - \int_0^T \dfrac{V_2(\tau(t))}{V_3(\tau(t))}x_2'(t) \,dt \\
			& ~~~~ \leq \delta^* - \dfrac{c}{C} \int_0^T x_2'(t) \,dt = \delta^* - \dfrac{c}{C} x_2(T) \leq \delta^* - \dfrac{c}{C} \delta/2 < 0,
	\end{align*}
	for $ \delta^* $ small enough. This will reach a contradiction. 
\end{pf}

We briefly recall some results on the classical non-rotating gaseous star problem
\begin{equation}\label{rdprb01}
	\begin{cases}
		\nabla \bigl(\tilde{K} \rho^{\gamma}\bigr) = -\rho \nabla\psi  & x \in \Omega,\\
		\Delta \psi = \rho & x \in \mathbb{R}^3,
	\end{cases}
\end{equation}
for some constants $ \tilde K > 0, \gamma > 1 $.
The following proposition is from \cite{Lin1997,Chandrasekhar1958}. We omit the proof here and refer to \cite{AuchmutyBeals1971,Lieb2001,Lieb1987,Luo2004,LinSS1989,Fu1998,KUAN1996}.
\begin{prop}\label{exiisentro}
	For fixed $ \tilde{K} > 0 $, $ \gamma > 1$, we have the following existence theory for \eqref{rdprb01},
	\begin{enumerate}
		\item Any regular solution to \eqref{rdprb01} with \eqref{totalmas} is spherically symmetric;
		\item For $ \gamma > 4/3 $, there is a regular ball type solution to \eqref{rdprb01} for any $ M > 0 $ in \eqref{totalmas}; 
		\item For $ \gamma = 4/3 $, there are infinitely many regular ball-type solutions to \eqref{rdprb01} only for $ M = M_c $ with $ M_c $ being a critical mass depending on $ \tilde{K} $;
		\item For $ 6/5 < \gamma < 4/3 $, there are unique regular ball type, singular ball type, and singular ground-state type solutions respectively with the same total mass $ M>0 $;
		\item For $ \gamma = 6/5 $, there is a unique ground-state type solution with finite total mass for any $ M > 0 $;
		\item For $ 1 < \gamma < 6/5 $, there is no solution to \eqref{rdprb01} with finite total mass.
	\end{enumerate}
By a ball type solution it means a solution with spherical symmetry and compact support. By a ground-state type solution it means a solution with unbounded support. By a singular solution it means a solution with $ \rho(0) = \infty $. By a regular solution it means a solution with $ \rho \in C(\bar \Omega) \cap C^1(\Omega) $. In addition, the physical vacuum \eqref{PV} on the boundary holds for the ball type solutions.
\end{prop}
Moreover, in the case of $ \gamma = 4/3 $, we have the following properties of the solutions. The following two propositions are based on the homology invariant property of \eqref{rdprb01} (see \cite{Chandrasekhar1958}).

\begin{prop}\label{ch4/3}
	For a fixed $ \tilde{K} > 0 $, and $ \gamma = 4/3 $, the solutions to \eqref{rdprb01} with $ M = M_c $ in \eqref{totalmas} are parametrised by $ \rho_s(x) = s^3 \bar{\rho}(sx) $ and $ \Omega_s = B_{s^{-1}\bar R}(0) $ for $ s > 0 $, where $ ( \bar{\rho}, B_{\bar R}(0) ) $ is one of the regular ball-type solutions to \eqref{rdprb01} with \eqref{totalmas}. In particular, the total mass is independent of the central density $ \rho_s(0) $. 
\end{prop}

\begin{pf}
	For any solution $ (\rho,\Omega) $ to \eqref{rdprb01} with $ \gamma = 4/3 $, i.e., 
	\begin{equation}\label{4/300}
	\begin{cases}
	\nabla\bigl(\tilde{K} \rho^{4/3}\bigr) = - \rho \nabla \psi   ~~~~ & x \in \Omega,\\
	\Delta \psi = \rho &  x \in \mathbb{R}^3,\\
	\int_{\Omega} \rho\,dx = M_c,
	\end{cases}
	\end{equation}
	after dividing $ \eqref{4/300}_1 $ with $ \rho $ and taking divergence on both sides, it reaches the following Lane-Emden equation,
	\begin{equation}\label{4/3001}
	4\tilde{K} \Delta \rho^{1/3} + \rho = 0    ~~~~ x \in \Omega.
	\end{equation}
	From Proposition \ref{exiisentro}, $ \rho $ is spherically symmetric. Define $ u(r) $ as the solution to the following ODE,
	\begin{equation}\label{ODE001}
	4 \tilde{K} \biggl(\dfrac{d^2}{dr^2} u(r) + \dfrac{2}{r} \dfrac{d}{dr} u(r) \biggr)+ u^3(r) = 0, ~~~~ \text{with} ~ u(0) = \rho^{1/3}(0), ~ \partial_r u(0) = 0.
	\end{equation}
	Assume $ R $ is the first zero of $ u $, i.e. $ u(r)>0 $ for $ 0\leq r < R $, $ u(R) = 0 $. Then $ \rho(r) : = u^3(r) $ for $ r \in \left\lbrack 0, R \right\rbrack $ is the solution to \eqref{4/3001} (or equivalently \eqref{4/300}), where $ r $ denotes the radius $ \abs{x}{} $.  Similarly, define $ \bar u $ to be the solution to \eqref{ODE001} with $ \bar u(0) = \bar \rho^{1/3}(0), \dr \bar u(0) = 0 $. Let $ \bar R $ be the first zero of $ \bar u $. Then $ \bar \rho $ agrees with $ \bar u^3 $ for $ r \in \left\lbrack 0, \bar R\right\rbrack $. Meanwhile, define \begin{equation}
		g(r) = s\bar u(sr) = s \bar{\rho}^{1/3}(sr) ~~~~\text{with} ~ s=\dfrac{u(0)}{\bar{u}(0)}. 
	\end{equation}
	Simple calculation shows that $ g $ satisfies the same ODE \eqref{ODE001} with the same initial values as $ u $, and the first zero of $ g $ is $ s^{-1}\bar{R} $. The uniqueness theory of ODE then implies $ u \equiv g $, i.e.
	\begin{equation}
		\rho^{1/3}(r) \equiv s\bar\rho^{1/3}(sr), ~~~~ s R = \bar R.
	\end{equation}
	It is easy to verify
	\begin{equation*}
		\int_0^{R} \rho(s) s^2 \,ds = \int_0^{\bar R}\bar\rho(s) s^2\,ds.
	\end{equation*}
\end{pf}
%
%

The next proposition concerns the relation between $ \tilde{K} $ and $ M_c = M_c(\tilde{K}) $.

\begin{prop}\label{4/3}
	For $ \gamma = 4/3 $, the critical mass $ M_c $ depends only on the gas dynamic constant $ \tilde K $. Moreover, 
	$ M_c = M_c(\tilde{K}) $ is monotonically increasing with respect to $ \tilde{K} $. More precisely, there exists a constant $ M_1 > 0 $ such that $ M_c(\tilde{K}) = \tilde{K}^{3/2} M_1 $.
\end{prop}

\begin{pf}
	We already know from Proposition \ref{ch4/3}, the total (critical) mass $ M_c $ is independent of the central density. Without loss of generality, we start with the solution $ u_1 $ to the following ODE(i.e. \eqref{ODE001} with $ \tilde{K} = 1 $),
	\begin{equation}
	4 \biggl(\dfrac{d^2}{dr^2} u_1(r) + \dfrac{2}{r} \dfrac{d}{dr} u_1(r) \biggr)+ u_1^3(r) = 0, ~~~~ \text{with} ~ u_1(0) = 1, ~ \partial_r u_1(0) = 0,
	\end{equation}
	with the first zero $ R_1 $ of $ u_1 $. The total mass of $ u_1 $ is defined to be \begin{displaymath} M_1 = 4\pi \int_{0}^{R_1} u_1^3(r)r^2\,dr . \end{displaymath} For any $ s>0 $, define $ u_s(r) = u_1(s r) $. Then $ u_s $ satisfies
	\begin{displaymath}
	4 \biggl(\dfrac{d^2}{dr^2} u_s(r) + \dfrac{2}{r} \dfrac{d}{dr} u_s(r) \biggr)+ s^2 u_s^3(r) = 0, ~~~~ \text{with} ~ u_s(0) = 1, ~ \partial_r u_s(0) = 0,
	\end{displaymath}
	or
	\begin{displaymath}
		4 s^{-2}\biggl(\dfrac{d^2}{dr^2} u_s(r) + \dfrac{2}{r} \dfrac{d}{dr} u_s(r) \biggr)+ u_s^3(r) = 0,
	\end{displaymath}
	with the first zero $ R_s = s^{-1}R_1 $, and the total mass of $ u_s $ is given by \begin{displaymath} M_s = 4\pi \int_{0}^{R_s} u_s^3(r) r^2 \,dr = 4\pi \int_{0}^{s^{-1}R_1} u_1^3 (sr) r^2 \,dr = s^{-3} M_1. \end{displaymath} 
	For any solution $ (\rho, \Omega) $ to \eqref{4/300} with $ \tilde K > 0 $ and $ \rho(0) = 1 $, by choosing $ s = \tilde K^{-1/2} $, the uniqueness theory of ODE yields that $ \rho $ agrees with $ u_s^3 $ for $ 0 \leq r \leq R_s $ through a similar argument as in the proof of Proposition \ref{ch4/3}. In particular, the value of the critical mass is
	\begin{equation*}
		M_c = M_s = \bigl( \tilde K^{-1/2} \bigr)^{-3} M_1 = \tilde K^{3/2} M_1.
	\end{equation*}
%
%
\end{pf}


To overcome the difficulties caused by the coordinate singularity and the degeneracy of the density and the temperature in the problem \eqref{eq:LagrangianCoordinates}, the following form of Hardy's inequality and Poincar\'e inequality will be useful.

\begin{lm} \label{lm:hardy} Let $ k $ be a given real number, and let $ g $ be a function satisfying $ \int_0^1 s^k (g^2 + g'^2)\,dx < \infty $.
	\begin{enumerate}
		\item[a)] If $ k > 1 $, then we have
		\begin{equation*}
			\int_0^1 s^{k-2} g^2\,ds \leq C \int_{0}^{1} s^k (g^2 + g'^2)\,ds.
		\end{equation*}
		\item[b)] If $ k < 1 $, then $ g $ has a trace at $ x=0 $ and moreover
		\begin{equation*}
			\int_0^1 s^{k-2} (g- g(0))^2\,ds\leq C\int_0^1 s^k g'^2 \,ds.
		\end{equation*}
		\item[c)] In particular, if $ g(1) = 0 $, the following Poincar\'e type inequality holds
		\begin{equation*}
			\int_0^1 s^2 g^2 \,ds \leq C \int_0^1 s^2 g'^2 \, ds. 
		\end{equation*}
		\item[d)] Suppose $ g(s) $ is defined in $ s \in [0,\omega) $ and $ g(0) = 0 $, the standard Poicar\'e inequality reads,
		\begin{equation*}
			\int_0^\omega g^2 \,ds \leq C \omega^2 \int_0^\omega g'^2 \,ds.
		\end{equation*}
	\end{enumerate}
\end{lm}
\begin{pf}
	While d) is the rescaling of the standard Poincar\'e inequality 
	, we refer a) and b) to \cite{Jang2014} and mainly show c) in the following. Notice, by applying the Cauchy inequality and a),
	\begin{align*}
		& s^2 g^2 = \biggl( \int_s^1 ( \tau g' + g) \,d\tau \biggr)^2 \leq (1-s) C \biggl( \int_s^1 \tau^2 g'^2\,d\tau + \int_s^1 g^2 \,d\tau \biggr) \\
		& ~~ \leq (1-s) C \biggl( \int_0^1 \tau^2 g^2 \,d\tau + \int_0^1 \tau^2 g'^2 \,d\tau \biggr).
	\end{align*}
	Then
	\begin{align*}
		& \int_0^1 s^2 g^2 \,ds = \int_0^{s_0} s^2 g^2 \,ds +  \int_{s_0}^1 s^2 g^2 \,ds \\
		& ~~~~~~\leq (s_0 - 1/2 s_0^2) C \biggl( \int_0^1 \tau^2 g^2 \,d\tau + \int_0^1 \tau^2 g'^2 \,d\tau \biggr) + \int_{s_0}^1 g^2\,ds. 
	\end{align*}
	Then by choosing $ s_0 > 0 $ small enough, it holds
	\begin{align*}
		& \int_0^1 s^2g^2\,ds \leq C_{s_0} \int_0^1 s^2 g'^2\,ds + C_{s_0} \int_{s_0}^1 g^2 \,ds \\
		& ~~~~~~ \leq C_{s_0} \int_0^1 s^2 g'^2\,ds + C_{s_0}\int_{s_0}^1 g'^2\,ds \leq C_{s_0} \int_0^1 s^2 g'^2\,ds,
	\end{align*}
	where in the second inequality we have applied the Poincar\'e inequality. 
\end{pf}

In the following, we adopt the notations
\begin{equation}
	\int \cdot \,dt : = \int_0^T \cdot\,dt, ~ \int \cdot \,dx : = \int_0^{R_0} \cdot\,dx.
\end{equation}
Also, for the constant $ \delta \in (0,1) $, $ C_\delta \simeq 1/\delta > 1 $.

\section{The Steady States \eqref{rdmd01}}\label{sec:steadysol}

Before showing the non-existence and existence theory of the steady state problem \eqref{rdmd01}, we first show the parallel relation 
\begin{equation}\label{prop:parallelrelation} \nabla \psi \parallel \vec{n} ~ \text{on} ~ \partial \Omega \end{equation} 
for the regular solution defined in \eqref{def:regularsol}. 

From $ \eqref{rdmd01}_1 $, $ K \theta \nabla \rho + K \rho \nabla \theta = - \rho \nabla \psi $, which is equivalent to say, in $ \Omega $,
	\begin{equation}\label{nonextc}
		K \dfrac{\theta}{\rho} \nabla \rho + K \nabla \theta = - \nabla \psi.
	\end{equation}
On one hand, $ \nabla \theta $ and $ \nabla \psi $ are uniformly continuous and uniformly bounded up to the boundary $ \partial \Omega $, which follows from the standard elliptic estimates on the equations \subeqref{rdmd01}{2} and \subeqref{rdmd01}{3}. This implies
$$ V := \dfrac{\theta}{\rho}\nabla \rho $$
is uniformly continuous and uniformly bounded up to the boundary. Meanwhile, the strong maximum principle and the Hopf lemma imply that $ \theta \simeq d(x) $ where $ d(x) $ denotes the distance from $ x \in \Omega $ to the boundary and $ \nabla \theta \parallel \vec{n} $ on $ \partial \Omega $. By employing Lemma \ref{lm:pllrltn}, $ V $ is parallel to the normal direction and $ V \cdot \vec{n} \leq 0 $ on the boundary $ \partial \Omega $. Therefore, it holds that $ \nabla \psi = - K \nabla \theta - K V $ is parallel to $ \vec{n} $ on the boundary and
\begin{equation}\label{prop:snbndint}
	\int_{\partial \Omega} \dfrac{\theta}{\rho}\nabla \rho \cdot \vec{n} \,dx \leq 0, ~ \int_{\partial \Omega} \theta \nabla \rho \cdot \vec{n} \,dx \leq 0,
\end{equation}
where the integrants are defined by \eqref{nonextc} on the boundary, i.e.
\begin{equation}\label{prop:traceoftheboundary}
\begin{gathered}
	\dfrac{\theta}{\rho}\nabla \rho \cdot \vec{n} := V\cdot \vec{n} = - \nabla \theta \cdot \vec{n} - \dfrac{1}{K} \nabla \psi \cdot \vec{n},\\
	\theta \nabla \rho \cdot \vec{n} := \rho V \cdot\vec{n} = - \rho \nabla \theta \cdot \vec{n} - \dfrac{\rho}{K} \nabla \psi \cdot \vec{n}.
\end{gathered}\end{equation}

In the following sections, we will show that $ 1/6 < \epsilon K < 1 $ is a necessary and sufficient condition for the existence of steady state solutions to \eqref{rdmd01} and study the degeneracy of the density and the temperature near the boundary.

\subsection{Non-Existence for High and Low Rate of Generation of Energy}\label{sec:nonexist}
In this section, we show that $ 1/6 < \epsilon K < 1 $ is a necessary condition for the existence of regular solutions to \eqref{rdmd01}. That is to say for $ \epsilon K \geq 1 $ or $ 0 < \epsilon K \leq 1/6 $, there is no regular solution. We show this in the following two lemmas.

\begin{lm}
	There is no regular solution to \eqref{rdmd01} with $ \epsilon K \geq 1 $.
\end{lm}

\begin{pf}
	We prove this by contradiction. Suppose for some $ \epsilon $ satisfying $ \epsilon K \geq 1 $, there is a regular solution $ (\rho, \theta, \Omega) $ to \eqref{rdmd01} with $ \eqref{totalmas}, ~ \eqref{ztmb} $. By dividing the first equation $ \eqref{rdmd01}_1 $ by $ \rho $ and taking divergence in the resulting equation, one can derive
	\begin{equation}\label{nonex01}
	K \dv \biggl( \dfrac{\theta}{\rho} \nabla \rho \biggr) =\left(\epsilon K - 1 \right)\rho.
	\end{equation}
	Integration this equation in $ \Omega $ yields
	\begin{equation}\label{nonex03}
	K \int_{\partial \Omega} \dfrac{\theta}{\rho} \nabla \rho \cdot \vec{n}\,dS = \left(\epsilon K - 1 \right) \int_{\Omega} \rho\,dx = M \left(\epsilon K - 1\right),
	\end{equation}
	where $ \vec{n} $ denotes the exteriorly normal vector on the surface $ \partial \Omega $ and we have used \eqref{totalmas}. The integration on the left of \eqref{nonex03} is understood in the sense of trace (see \eqref{prop:traceoftheboundary}). Therefore, \eqref{prop:snbndint} implies
	\begin{equation*}
	\epsilon K - 1 \leq 0 ~~~ \text{or} ~~~ \epsilon K \leq 1.
	\end{equation*}
	Thus $ \epsilon K = 1 $. Hence, the right of \eqref{nonex01} vanishes. Multiply \eqref{nonex01} with $ \rho $ and integrate the resulting equation in $ \Omega $. It follows,
	\begin{equation}\label{nonex02}
	- K \int_{\Omega} \dfrac{\theta}{\rho} |\nabla \rho|^2 \, dx + K \int_{\partial \Omega} \theta \nabla \rho \cdot \vec{n}\,dS= 0.
	\end{equation} 
	However, for a regular solution,
	\begin{displaymath}
	\int_{\Omega} \dfrac{\theta}{\rho} |\nabla \rho|^2 \, dx > 0 ~~ \text{and} ~~ \int_{\partial \Omega} \theta \nabla \rho \cdot \vec{n} \, dS \leq 0, 
	\end{displaymath}
	(see \eqref{prop:snbndint}) and therefore,  
	\begin{displaymath}
		- K \int_{\Omega} \dfrac{\theta}{\rho} |\nabla \rho|^2 \, dx + K \int_{\partial \Omega} \theta \nabla \rho \cdot \vec{n}\,dS < 0,
	\end{displaymath}
	which contradicts \eqref{nonex02}.
\end{pf}

\begin{lm}
	There is no regular solution to \eqref{rdmd01} with $ 0 < \epsilon K \leq \dfrac{1}{6} $.
\end{lm}

\begin{pf}
	Again, we will prove this through a contradictory argument. Suppose that for some $ \epsilon $ satisfying $ 0 < \epsilon K \leq \frac{1}{6} $, there is a regular solution $ (\rho, \theta, \Omega) $ to \eqref{rdmd01}. We employ a Pohozaev type argument as in \cite{Deng2002}. Multiply $ \eqref{rdmd01}_{1} $ with $ x $ and integrate the resulting in $ \Omega $
	\begin{displaymath}
	\int_{\Omega} \nabla \left(K\rho \theta \right) \cdot x\,dx = - \int_{\Omega} \rho \nabla \psi \cdot x \,dx.
	\end{displaymath}
	Integration by parts on the left yields
	\begin{equation}\label{nonex001}
	- 3K \int_{\Omega}\rho \theta \,dx = - \int_{\Omega} \rho \nabla \psi \cdot x\,dx.
	\end{equation}
	Meanwhile, multiply $ \eqref{rdmd01}_3 $ with $ \nabla \psi \cdot x $, and integrate the resulting in $ \Omega $
	\begin{equation}\label{nonex002}
	\int_{\Omega} \Delta \psi \left(\nabla \psi \cdot x \right) \,dx = \int_{\Omega}\rho \left(\nabla \psi \cdot x\right)\,dx.
	\end{equation}
	Integration by parts on the left yields,
	\begin{align}\label{nonex003}
	& \int_{\Omega} \Delta \psi \left(\nabla \psi \cdot x \right) \,dx {\nonumber} \\
	& ~~~~~~~~ = - \int_{\Omega} |\nabla \psi|^2\,dx - \dfrac{1}{2} \int_{\Omega} x \cdot \nabla \left|\nabla \psi\right|^2 \,dx + \int_{\partial \Omega} \left(\nabla \psi \cdot \vec{n}\right) \left(\nabla \psi \cdot x\right) \,dS {\nonumber} \\
	& ~~~~~~~~ =  \dfrac{1}{2} \int_{\Omega} |\nabla \psi|^2\,dx + \dfrac{1}{2} \int_{\partial \Omega} \left| \nabla \psi \right|^2 x\cdot \vec{n}\,dS,
	\end{align}
	where we have used the fact
	\begin{displaymath}
		\int_{\partial \Omega} \left(\nabla \psi \cdot \vec{n}\right) \left(\nabla \psi \cdot x\right) \,dS = \int_{\partial \Omega} |\nabla \psi|^2 x \cdot\vec{n} \,dS
	\end{displaymath}
	due to $ \nabla \psi \parallel \vec{n} $ on $ \partial \Omega $ (see \eqref{prop:parallelrelation}). 
	\eqref{nonex001}, \eqref{nonex002}, \eqref{nonex003} yield the identity
	\begin{equation}\label{nonex006}
	- 3K \int_{\Omega} \rho \theta \,dx = - \dfrac{1}{2} \int_{\Omega}\left|\nabla\psi\right|^2\,dx - \dfrac{1}{2}\int_{\partial \Omega} \left|\nabla\psi\right|^2 x\cdot\vec{n}\,dS.
	\end{equation}
	
	On the other hand, multiply $ \eqref{rdmd01}_3 $ with $ \psi $ and integrate the resulting in $ \Omega $. After integration by parts, it follows
	\begin{equation}\label{nonex004}
	- \int_{\Omega} \left|\nabla \psi\right|^2 \,dx + \int_{\partial \Omega} \psi \left(\nabla\psi \cdot \vec{n}\right)\,dS = \int_{\Omega} \rho \psi \,dx.
	\end{equation}
	Multiply $ \eqref{rdmd01}_2 $ with $ \psi $ and $ \eqref{rdmd01}_3 $ with $ \theta $ and integrate the resulting expressions in $ \Omega $. After integration by parts,
	\begin{gather*}
	\int_{\Omega} \nabla \theta \cdot \nabla \psi \,dx - \int_{\partial \Omega} \psi \left(\nabla \theta \cdot \vec{n}\right) \,dS = \epsilon \int_{\Omega} \rho \psi \,dx, \\
	- \int_{\Omega} \nabla \psi \cdot \nabla \theta \,dx=  \int_{\Omega} \rho \theta \,dx,
	\end{gather*}
	which together with \eqref{nonex004} yields
	\begin{equation}\label{nonex005}
	\begin{aligned}
	& \int_{\Omega} \rho \theta \,dx =  - \int_{\partial \Omega} \psi \left(\nabla \theta \cdot \vec{n}\right)\,dS - \epsilon\int_{\Omega}\rho \psi \,dx \\
	& ~~~~~~~ = - \int_{\partial \Omega} \psi \left(\nabla \theta \cdot \vec{n}\right)\,dS + \epsilon \int_{\Omega}\left|\nabla \psi\right|^2\,dx -  \epsilon \int_{\partial \Omega} \psi \left(\nabla\psi \cdot \vec{n}\right)\,dS.
	\end{aligned}
	\end{equation}
	Therefore, $ \eqref{nonex005} \times 3K + \eqref{nonex006} $ implies 
	\begin{equation}\label{nonex007}
	\begin{aligned}
	& \dfrac{1-6\epsilon K}{2}\int_{\Omega} \left|\nabla\psi\right|^2\,dx  = - \dfrac{1}{2} \int_{\partial\Omega} \left|\nabla\psi\right|^2 x\cdot\vec{n}\,dS\\
	& ~~~~~~~~~~~~~~~~ - 3K \biggl(\int_{\partial\Omega}\psi \bigl(\nabla\theta\cdot\vec{n}\bigr)\,dS + \epsilon \int_{\partial\Omega}\psi \bigl(\nabla\psi\cdot\vec{n}\bigr)\,dS\biggr).
	\end{aligned}
	\end{equation}
	Notice, since $ \theta = 0 $ on $ \partial \Omega $,
	\begin{equation}\label{nonex009}
	\begin{aligned}
	& \int_{\partial\Omega}\psi \left(\nabla\theta\cdot\vec{n}\right)\,dS + \epsilon \int_{\partial\Omega}\psi \left(\nabla\psi\cdot\vec{n}\right)\,dS = \dfrac{1}{\epsilon} \int_{\partial \Omega} \left(\theta + \epsilon\psi \right) \nabla \left(\theta + \epsilon \psi\right)\cdot \vec{n}\,dS\\
	& ~~~~~~ = \dfrac{1}{2\epsilon} \int_{\partial \Omega} \nabla \left(\theta + \epsilon \psi\right)^2 \cdot \vec{n} \,dS.
	\end{aligned}
	\end{equation}
	Meanwhile, from $ \eqref{rdmd01}_2 $ and $ \eqref{rdmd01}_3 $,
	\begin{align*}
	&\Delta \left(\theta + \epsilon \psi \right)^2 = 2 \left| \nabla \left( \theta + \epsilon \psi \right) \right|^2 + 2\left(\theta + \epsilon \psi\right) \Delta \left(\theta + \epsilon \psi\right)\\
	& ~~~~~~~= 2 \left| \nabla \left( \theta + \epsilon \psi \right) \right|^2 + 2\left(\theta + \epsilon \psi\right) \left( -\epsilon \rho + \epsilon \rho\right) = 2 \left| \nabla \left( \theta + \epsilon \psi \right) \right|^2 \geq 0.
	\end{align*}
	Integrating this expression in $ \Omega $, after integration by parts, implies
	\begin{equation}\label{nonex008}
	\int_{\partial \Omega} \nabla \left(\theta + \epsilon \psi \right)^2 \cdot \vec{n}\,dS = 2 \int_{\Omega} \left|\nabla\left(\theta + \epsilon \psi\right)\right|^2\,dx \geq 0.
	\end{equation}
	\eqref{nonex007}, \eqref{nonex009}, \eqref{nonex008} then yield
	\begin{displaymath}
	0 \leq \dfrac{1-6\epsilon K}{2}\int_{\Omega} \left|\nabla\psi\right|^2\,dx  \leq  - \dfrac{1}{2}\int_{\partial\Omega} \left|\nabla\psi\right|^2 x\cdot\vec{n}\,dS \leq 0,
	\end{displaymath}
	due to the contradictory assumption $ 0 < \epsilon K \leq \frac{1}{6} $ and the fact that $ \Omega $ is star-shaped. Therefore, $\nabla \psi \equiv 0 $ in $ \bar{\Omega} $. This is impossible for a regular solution. Otherwise, after integrating $ \eqref{rdmd01}_3 $ over $ \Omega $ it follows,
	\begin{displaymath} 0 = \int_{\partial \Omega} \nabla \psi \cdot \vec{n}\, dS = \int_{\Omega} \rho \,dx = M > 0.\end{displaymath} 
	Thus it reaches a contradiction. 
\end{pf}

\subsection{Existence of Steady State Solutions}\label{sec:existence}
In this section, we will construct the regular solution to \eqref{rdmd01} when $ 1/6 < \epsilon K < 1 $ by reducing the problem to the classical Lane-Emden equation. But before we start, it should be noticed that even though it looks like the non-existence theory of the Lane-Emden equation implies the non-existence theory in Section \ref{sec:nonexist}, it actually does not. In fact, in the following reduction process, it has already assumed $ \epsilon K < 1 $(see Lemma \ref{equvlc01}). Therefore the reduction does not work for the non-existence theory and the proof in Section \ref{sec:nonexist} is not redundant. 

We start with the following property of the regular solution to \eqref{rdmd01}. 
\begin{lm}\label{improp}
	Any regular solution to \eqref{rdmd01} will satisfy
	\begin{equation}\label{idtt01}
	\nabla \left(\theta + \epsilon \psi \right) = 0~~~~~~~ \text{in} ~ \Omega.
	\end{equation}
\end{lm}

\begin{pf}
	From $ \eqref{rdmd01}_2 $ and $ \eqref{rdmd01}_3 $, 
	\begin{equation}\label{harmonicid01}
	\Delta \left( \theta + \epsilon \psi \right) = -\epsilon \rho + \epsilon \rho = 0 ~~~~~~ x \in \Omega.
	\end{equation}
	Again, we employ a Pohozaev type argument. Multiply \eqref{harmonicid01} with $ \nabla \left(\theta + \epsilon \psi \right) \cdot x $, and integrate the resulting in $ \Omega $. Integration by parts then yields
	\begin{align*}
	& 0 = \int_{\Omega}\Delta \left(\theta + \epsilon \psi \right) \nabla \left(\theta + \epsilon \psi 
	\right) \cdot x \,dx\\ 
	& ~~~~ = - \int_{\Omega} \left|\nabla\left(\theta+\epsilon\psi\right)\right|^2\,dx - \dfrac{1}{2} \int_{\Omega} \nabla\left|\nabla\left(\theta+\epsilon\psi\right)\right|^2 \cdot x\,dx\\
	& ~~~~~~~~~~ + \int_{\partial \Omega} \left(\nabla\left(\theta+\epsilon\psi\right) \cdot \vec{n}\right) \left(\nabla\left(\theta+\epsilon\psi\right)\cdot x\right)\,dS\\
	& ~~~~ = \dfrac{1}{2} \int_{\Omega} \left|\nabla \left(\theta+\epsilon \psi\right)\right|^2 \,dx + \dfrac{1}{2} \int_{\partial \Omega} \left|\nabla\left(\theta+\epsilon \psi \right)\right|^2 x \cdot \vec{n}\,dS,
	\end{align*}
	where we have used the fact $ \nabla \left(\theta+\epsilon\psi\right) \parallel \vec{n} $ on $ \partial \Omega $ in the last equality (see \eqref{prop:parallelrelation}). Therefore, together with the assumption on the geometry of $ \Omega $, which yields the right is the sum of two nonnegative integrals, this identity shows $ \nabla \left(\theta + \epsilon \psi\right) \equiv 0 $ in $ {\Omega} $.
\end{pf}
\begin{rmk}
Indeed, $ \nabla \left(\theta+\epsilon\psi\right) \parallel \vec{n} $ on $ \partial \Omega $ implies that $ \theta + \epsilon \psi $ is constant on the boundary and hence from \eqref{harmonicid01} and the maximum principle for harmonic equations, $ \theta + \epsilon \psi $ is constant in $ \Omega $ \footnote{This is pointed out by Prof. YanYan Li in a private conversation.}. In particular, $ \nabla \theta + \epsilon \nabla \psi = 0 $ in $ \Omega $.
\end{rmk}
In the next lemma, we will show that the property \eqref{idtt01} implies that \eqref{rdmd01} can be reduced to the classical system for gaseous stars \eqref{rdprb01}, for some constants $ \tilde K > 0, \gamma > 1 $. 
\begin{lm}\label{equvlc01}
	For any regular solution $\left( \rho, \theta, \Omega\right) $ to \eqref{rdmd01} with $ 0 < \epsilon K < 1 $, $ \exists \tilde{K} > 0, ~ \gamma > 1 $ such that $ \left( \rho,\Omega\right) $ satisfies the system \eqref{rdprb01}. Moreover, there is a positive constant $ S > 0 $ such that
	\begin{equation}\label{pmtr}
	\rho^{ - \epsilon K } \theta^{1-\epsilon K} = S  ~~~~~~ x \in \Omega.
	\end{equation}
\end{lm}
\begin{pf}
	Let $\left( \rho, \theta, \Omega\right) $ be a regular solution to \eqref{rdmd01}. Then from \eqref{idtt01} in Lemma \ref{improp}, 
	\begin{displaymath}
	\nabla \psi = - \dfrac{1}{\epsilon}\nabla \theta.
	\end{displaymath}
	Together with $ \eqref{rdmd01}_1 $, it holds, 
	\begin{equation}
	K \theta \nabla \rho + K \rho \nabla \theta = \dfrac{\rho}{\epsilon} \nabla \theta,
	\end{equation}
	which yields, for $ \epsilon K \neq 1 $,
	\begin{displaymath}
	\nabla \left(\rho^{- \epsilon K } \theta^{1 - \epsilon K }\right) = 0 ~~~~ x \in \Omega.
	\end{displaymath}
	Since $ \rho, \theta $ are both positive and finite in $ \Omega $, this means $ \exists S > 0 $, such that \eqref{pmtr} holds.
	Define $
	\tilde{K} = S^{\frac{1}{1-\epsilon K}} K > 0, ~ \gamma = 1 + \dfrac{\epsilon K }{1-\epsilon K} > 1 $.
	Then
	\begin{equation}\label{rdpm}
	\theta = \left(S \rho^{\epsilon K}\right)^{\frac{1}{1-\epsilon K}} ~~ \text{and} ~~ P = K \rho \theta = K \rho  \left(S\rho^{\epsilon K}\right)^{\frac{1}{1-\epsilon K}} =  \tilde{K} \rho^{\gamma}.
	\end{equation}
	This finishes the proof.  
\end{pf}


In the following lemma, we prove $ 1/6 < \epsilon K < 1 $ is not only a necessary, but also a sufficient condition to guarantee the existence of regular solutions to \eqref{rdmd01}.

\begin{lm}\label{exists}
	For $ \dfrac{1}{6} <  \epsilon K < 1 $ and any constant $ M > 0 $, there exist infinitely many regular solutions to \eqref{rdmd01} with \eqref{totalmas} and \eqref{ztmb}. Moreover, all regular solutions are spherically symmetric.
\end{lm}

\begin{pf}
	In the following, we write  \begin{displaymath}
		\gamma_\epsilon = 1 + \frac{\epsilon K}{1 - \epsilon K} = \frac{1}{1-\epsilon K}.
	\end{displaymath}
	
	{\bf Case 1.}
	$ \frac{1}{6} < \epsilon K < \frac{1}{4} $ or $ \frac{1}{4} < \epsilon K < 1 $.  
	
	In this case,  $ \frac{6}{5} < \gamma_\epsilon < \frac{4}{3} $ or $ \gamma_\epsilon > \frac{4}{3} $. For any fixed $ S > 0 $, let $ \tilde{K}_{S} = S^{\gamma_\epsilon} K $. From Proposition \ref{exiisentro}, there is a unique regular ball type solution $ (\rho_S, \Omega_S) $ to \eqref{rdprb01} and \eqref{totalmas} with $ \tilde{K} = \tilde{K}_{S} $. Inspired by \eqref{rdpm}, define 
	\begin{displaymath} \theta_S = \frac{\tilde{K}_{S}}{K} \rho_S^{\gamma_{\epsilon}-1} = S^{\gamma_\epsilon} \rho_S^{\gamma_\epsilon-1}. \end{displaymath} 
	Then it is easy to verify that $ (\rho_S, \theta_S, \Omega_S) $ is a solution to \eqref{rdmd01} with \eqref{totalmas} and \eqref{ztmb} for any $ S>0 $. Therefore, by choosing different values of $ S > 0 $, there are infinitely many regular solutions to \eqref{rdmd01}. Therefore, $ S $ given above is a parametrisation of these regular solutions. Notice, $ \rho_S, \theta_S, S $ satisfy the relation \eqref{pmtr}.
	
	{\bf Case 2.}
	$ \epsilon K = \frac{1}{4} $, $ \gamma_\epsilon = \frac{4}{3} $. 
	
	From Proposition \ref{4/3}, there is a $ \tilde{K}_{4} = \bigl(\frac{M}{M_1}\bigr)^{2/3} $, such that 
	\begin{equation}\label{par4/3001}
		M_c(\tilde{K}_4) = \tilde{K}_4^{3/2} M_1 =  M.
	\end{equation} 
	By taking $ \tilde{K} = \tilde{K}_4 $ in \eqref{4/300}, Proposition \ref{exiisentro} yields there are infinitely many regular ball type solutions $ \left( \rho,\Omega \right) $ to \eqref{4/300} with prescribed total mass $ M = M_c(\tilde K_4) $. Similarly, by defining 
	\begin{displaymath} \theta = \dfrac{\tilde{K}_{4}}{K} \rho^{1/3}, \end{displaymath}
	 then $ \left( \rho, \theta, \Omega \right) $ is a regular solution to \eqref{rdmd01}. The multiplicity of regular ball type solutions to \eqref{4/300} yields the multiplicity of regular solutions to \eqref{rdmd01}. Moreover, $ S $ given by \eqref{pmtr} satisfies	 \begin{equation}\label{par4/3}
	 	S = \rho^{-1/4} \theta^{1-1/4} = \biggl(\dfrac{\tilde{K}_{4}}{K}\biggr)^{3/4} = \biggl( \dfrac{M}{K^{3/2} M_1} \biggr)^{1/2},
	 \end{equation}
	 which depends only on the total mass $ M $ and the coefficient $ K $. In this case, $ S $ does not parametrize the solutions.
\end{pf}
\begin{rmk}
	In the case $ 1/6 \leq \epsilon K < 1/4 $ , similar arguments actually also show the existence of singular solutions or ground-state type solutions to \eqref{rdmd01}.
\end{rmk}

We claim that the regular solutions obtained in Lemma \ref{exists} are all the regular solutions to \eqref{rdmd01}. This is a direct consequence of Lemma \ref{equvlc01}.

\begin{lm}\label{lm:alsol} For $ \dfrac{1}{6} <  \epsilon K < 1 $ and $ M > 0 $, the regular solutions obtained in Lemma \ref{exists} are all the regular solutions to \eqref{rdmd01}.\end{lm}

Notice, $ S $ given in \eqref{pmtr} does not always parametrize the regular solutions to \eqref{rdmd01}. In the next lemma, we show that the solutions can be parametrized by a scaling variable.

\begin{lm}[Homology Invariance]\label{lm:hit}
Given a regular solution $ (\rho, \theta, \Omega) $ to \eqref{rdmd01} with \eqref{totalmas} and \eqref{ztmb}, for $ s > 0 $, define
\begin{equation}\label{hlinvrt}
\begin{gathered}
	\rho_s(x) = s^3 \rho(s x), ~ \theta_s(x) = s \theta(s x),  \\
	\Omega_s = \left\lbrace x \in \mathbb{R}^3 | \rho_s(x) > 0\right\rbrace = \left\lbrace x = s^{-1}y | \rho(y)>0 \right\rbrace.
	\end{gathered}
\end{equation}
Then $(\rho_s, \theta_s, \Omega_s)$ is also a regular solution to \eqref{rdmd01} with the same total mass $ M > 0 $.
\end{lm}

\begin{pf}
The proof follows from a direct calculation and is omitted here.
\end{pf}

%

\begin{lm}[Parametrization of Regular Solutions]\label{lm:chrs}
	For any fixed $ \epsilon $ in \eqref{rdmd01} with $ \frac{1}{6} < \epsilon K < 1 $ and $ M>0 $, all regular solutions to \eqref{rdmd01} with \eqref{totalmas} and \eqref{ztmb} are parametrised by $ s > 0 $ as in \eqref{hlinvrt}, where $ (\rho, \theta, \Omega) $ is one of the regular solutions to \eqref{rdmd01} with \eqref{totalmas} and \eqref{ztmb} given in Lemma \ref{exists}.
\end{lm}

\begin{pf}
The proof follows from the construction of regular solutions Lemma \ref{exists}, Lemma \ref{lm:alsol}, Lemma \ref{lm:hit} and the parametrization of the solutions to \eqref{4/300} in Proposition \ref{ch4/3}. 
\end{pf}

A direct consequence of the homology invariance \eqref{hlinvrt} is the following. 

\begin{cor}
	For fixed $ M > 0 $ and $ \frac{1}{6} < \epsilon K < 1 $, let $ \left(\rho, \theta, \Omega = B(R) \right) $ be the regular solution to \eqref{rdmd01} with \eqref{totalmas} and \eqref{ztmb}. The radius $ R $ of $ \Omega $ is monotonically decreasing with respect to the central density $ \rho(0) $ and the central temperature $ \theta(0) $. And the central density $ \rho(0) $ is monotonically increasing with respect to the central temperature $ \theta(0) $. In particular, there exist constants $ C_1 $ and $ C_2 $, depending on the total mass $M$ and $\epsilon K $,  such that $ R \cdot \theta(0) = C_1 $ and $ \rho(0) = C_2 \theta^{3}(0) $, and consequently $ R^3 \rho(0) = C_1^3 C_2 $. 
\end{cor}

In the following lemma, we point out the vacuum property of $ \rho, \theta $ near the gas-vacuum interface $ \partial \Omega $ for the regular solution $(\rho, \theta,\Omega) $ to \eqref{rdmd01}. 

\begin{lm}\label{stationary:boundary}
For any regular solution $ (\rho, \theta, \Omega) $ to \eqref{rdmd01} with $ 1/6 < \epsilon K < 1 $, $ \theta $ is Lipschitz continuous across $ \partial \Omega $, while $ \rho $ is $ \frac{1-\epsilon K}{\epsilon K} $-H\"{o}lder continuous across $ \partial \Omega $. Moreover, it obeys the following asymptotic behaviour near the origin and the boundary respectively,
	\begin{enumerate}
		\item $ \rho, \theta $ is analytic at the origin, and in particular, $ \theta = a + b r^2 + O(r^4),~ r \sim 0 $ with $ \dr^{2k+1} \theta(0) = 0 $ for any nonnegative integer $ k \geq 0 $ ; 
		\item $ \rho ,  \theta \bigr|_{\partial\Omega} = 0, ~-\infty < \dr \rho^{\frac{\epsilon K }{1-\epsilon K}} , \dr \theta \bigr|_{\partial\Omega} \leq - C < 0 $;
		\item $ \abs{\dr^k \rho}{} \leq O(\rho^{\frac{1-(1+k)\epsilon K}{1-\epsilon K}}) 
		$, $ \abs{\dr^k \theta}{} \leq O(1) + O(\theta^{\frac{1-\epsilon K}{\epsilon K} - k + 2}) $ for any $ k > 0 $ near the boundary,
	\end{enumerate}
	where $ \dr $ represents the derivative along the radial direction of $ \Omega $.
\end{lm}
\begin{pf}
	Similar as before, denote the constants
	\begin{equation*}
			S = \rho^{-\epsilon K} \theta^{1-\epsilon K}, \gamma_\epsilon = 1 + \dfrac{\epsilon K}{1-\epsilon K}, \alpha_\epsilon = \dfrac{1-\epsilon K}{\epsilon K}, \bar K = \dfrac{1}{\epsilon K} S^{\frac{1}{1-\epsilon K}} K = \dfrac{S^{1/(1-\epsilon K)}}{\epsilon}.
	\end{equation*}
	Then, $ \alpha_\epsilon \in (0,5), \gamma_\epsilon \in (6/5, \infty) $. 
	We study the following ODE, 
	\begin{equation}\label{ODE:Lane-Emden}
		\bar K \left( \partial_{rr}+ \dfrac{2}{r} \partial_r \right) u + u^{\alpha_\epsilon} = 0, ~ u(0) = \rho(0)^{\epsilon K /(1-\epsilon K)} , \dr u (0) = 0.
	\end{equation}
	Also, let $ 0 < R < \infty $ be the first zero of $ u $. It can be verified, similar as before, $ \rho = u^{\alpha_\epsilon}, \theta = S^{1/(1-\epsilon K)} u $ for $ r \in \left\lbrack 0, R \right\rbrack $. Without loss of generality, we assume $ \bar K = 1 $ and $ \rho(0) = 1 $. We study $ u $ first. Similar to \cite[Lemma 3.3]{Jang2014}, we claim
	\begin{enumerate}
		\item[a)] $ u $ is analytic at the origin. $ u(r) = 1 - b r^2 + O(r^4), ~ r \sim 0 $,
		for some positive constant $ b > 0 $. Also, $ \dr^{2k+1} u(0) = 0 $ for any nonnegative integer $ k \geq 0 $;
		\item[b)]  $ \abs{\dr^k u}{} \leq O(1) + O(u^{\alpha_\epsilon - k + 2})$ for $ k\geq 1,  r \sim R $, and $ \abs{\dr u}{} = O(1) $ for $ r \sim R $. 
	\end{enumerate}
	Then the lemma follows after a substitution. Indeed, the property for $ \theta $ is a direct consequence of a) and b). We briefly demonstrate how to get the asymptotic behaviour of $ \dr^k \rho $ near the boundary. Notice, $ u = \rho^{\gamma_\epsilon-1} $. It is easy to verify $ \abs{\dr\rho}{}= O(\rho^{2-\gamma_\epsilon} )= O	( \rho^{\frac{1-2\epsilon K}{1-\epsilon K}} ) $. Moreover, it holds,
	\begin{equation}\label{sub01}
		\abs{\rho^{\gamma_\epsilon -2}\dr^k \rho}{} 
		- L_k \leq \abs{\dr^k u}{} 
		\leq O(1) + O(\rho^{1-(k-2)(\gamma_\epsilon-1)}),
	\end{equation}
	where $ L_k $ contains the terms look like the absolute value of
	\begin{equation*}
		\rho^{\gamma_\epsilon-1-n} (\dr^{\alpha_1^k}\rho)(\dr^{\alpha_2^k}\rho) \cdots (\dr^{\alpha_n^k}\rho),
	\end{equation*}
	with $ \alpha_1^k + \alpha_2^k + \cdots + \alpha_n^k = k, 2 \leq n \leq k, 1 \leq \alpha_i^k < k $ for $ i = 1,2,\cdots,n $. Then, if it is assumed $ \abs{\dr^l \rho}{} \leq O(\rho^{\frac{1-(1+l)\epsilon K}{1-\epsilon K}}) $ holds for $ 0 < l < k $, $ L_k $ can be bounded by
	\begin{equation*}
		\rho^{\gamma_\epsilon-1-k \frac{\epsilon K}{1-\epsilon K}}.
	\end{equation*}
	Therefore, it shall hold from \eqref{sub01}
	\begin{equation*}
		\abs{\dr^k \rho}{} \leq O(\rho^{2-\gamma_\epsilon}) + O( \rho^{3-\gamma_\epsilon - (k-2)(\gamma_\epsilon - 1)} ) + O(\rho^{1-k\frac{\epsilon K}{1-\epsilon K}}) \leq O(\rho^{1-k\frac{\epsilon K}{1-\epsilon K}}).
	\end{equation*}
	This will finish the proof via the induction principle. What is left is to show the claim. In fact, it is standard and following the same program as in \cite{Jang2014}. We attach the proof here for the sake of convenience. 
	
	To show a), rewrite \eqref{ODE:Lane-Emden} as
	\begin{equation}\label{ODE:Lane-Emden001}
		r u_{rr} + 2 u_{r} + c r u^{\alpha_\epsilon} = 0,
	\end{equation}
	for some constant $ c > 0 $.
	After dividing this expression by $ r $ and taking the limit $ r \rightarrow 0^+ $, it holds $ 3 u_{rr}(0) + c = 0 $. Next, differentiating \eqref{ODE:Lane-Emden001} once yields the identity,
	\begin{equation}\label{ODE:Lane-Emden002}
		r u_{rrr} + 3 u_{rr} + c(\alpha_\epsilon r u^{\alpha_\epsilon-1} u_r + u^{\alpha_{\epsilon}}) = 0.
	\end{equation}
	Again, dividing \eqref{ODE:Lane-Emden002} with $ r $ will give the following,
	\begin{equation*}
		u_{rrr} + 3\dfrac{u_{rr} - u_{rr}(0)}{r} + c \left( \alpha_\epsilon u^{\alpha_\epsilon-1} u_r + \dfrac{u^{\alpha_\epsilon}-u^{\alpha_\epsilon}(0)}{r} \right) = 0,
	\end{equation*}
	where it has been applied the relation $ 3u_{rr}(0) + c = 0 $ and $ u(0) = 1 $. Passing the limit $ r \rightarrow 0^+ $ then yields $ 4u_{rrr}(0) = 0 $. To calculate $ u_{rrrr}(0) $, from \eqref{ODE:Lane-Emden002},
	\begin{equation*}
		r u_{rrrr} + 4 u_{rrr} + c(\alpha_\epsilon(\alpha_\epsilon-1) r u^{\alpha_\epsilon-2} u_r^2 + \alpha_{\epsilon}u^{\alpha_\epsilon-1} ( \underbrace{ r u_{rr} + 2u_r}_{-cru^{\alpha_\epsilon}} )) = 0.
	\end{equation*}
	Thus $ 5u_{rrrr}(0) - \alpha_{\epsilon} c^2 = 0 $. In general, for any integer $ k \geq 1 $, apply $ \dr^{2k-2} $ to \eqref{ODE:Lane-Emden001} and use Leibniz's rule to get
	\begin{equation}\label{ODE:Lane-Emden003}
		r \dr^{2k} u + 2k \dr^{2k-1} u + c \dr^{2k-2} (r u^{\alpha_\epsilon} ) = 0.
	\end{equation}
	Now we claim for any positive integer $ k \geq 1 $
	\begin{equation}\label{ODE:Lane-Emden004} 
		\dr^{2k-2}(ru^{\alpha_\epsilon})(0) = 0 ~~ \text{and} ~~ \dr^{2k-1} u(0) = 0. 
	\end{equation} 
	This has been shown for $ k = 1,2 $. By the principle of mathematical induction, we only need to show the relation \eqref{ODE:Lane-Emden004} with $ k = n+1 $ under the assumption that \eqref{ODE:Lane-Emden004} holds for all $ k \leq n $. Consider \eqref{ODE:Lane-Emden003} with $ k = n $, which can be written as.
	\begin{equation*}
		\dr^{2n} u + 2n\dfrac{\dr^{2n-1}u}{r} + c \dfrac{\dr^{2n-2}(r u^{\alpha_\epsilon})}{r} = 0.
	\end{equation*}
	Passing the limit $ r \rightarrow 0^+ $ and noticing \eqref{ODE:Lane-Emden004}, it holds 
	\begin{equation}\label{ODE:Lane-Emden005}
		(2n+1) \dr^{2n} u (0) + c l_n = 0,
	\end{equation} with $ l_n = \lim_{r\rightarrow 0^+} \frac{\dr^{2n-2}(ru^{\alpha_{\epsilon}})}{r} = \dr^{2n-1}(r u^{\alpha_\epsilon})(0) $. Again, differentiate \eqref{ODE:Lane-Emden003} with $ k = n $ with respect to $ r $, and it holds
	\begin{equation*}
		r \dr^{2n+1} u + (2n + 1) \dr^{2n} u + c \dr^{2n-1}(r u^{\alpha_\epsilon}) = 0,
	\end{equation*}
	or equivalently, by making use of \eqref{ODE:Lane-Emden005}
	\begin{equation*}
		\dr^{2n+1} u + (2n+1) \dfrac{\dr^{2n}u - \dr^{2n}u(0)}{r} + c \dfrac{\dr^{2n-1}(r u^{\alpha_\epsilon}) - \dr^{2n-1}(r u^{\alpha_\epsilon})(0) }{r}=0.
	\end{equation*}
	Passing the limit $ r \rightarrow 0^+ $ then yields $(2n+2) \dr^{2n+1} u(0) + c \dr^{2n} (r u^{\alpha_\epsilon})(0) = 0 $. Meanwhile, direct calculation shows
	\begin{equation}
		\dr^{2n}(r u^{\alpha_\epsilon})(0) = (r \dr^{2n}u^{\alpha_\epsilon}) (0) + 2n (\dr^{2n-1} u^{\alpha_\epsilon})(0) = 0 + 0 = 0,
	\end{equation}
	where we have used the fact that the expansion of $\dr^{2n-1} u^{\alpha_\epsilon}$ contains at least one multiplier with an odd-order derivative of $ u $. Thus $\dr^{2n+1} u(0) = - \frac{c}{2n+2}\dr^{2n}(ru^{\alpha_\epsilon})(0) =  0 $. This finishes the proof of \eqref{ODE:Lane-Emden004}. 
	
	In order to show b), it suffices to study $ u $ near $ r = R $. We already know $ \abs{u_{r}}{} = O(1) $ for $ r \simeq R $. Hence, $ \abs{u_{rr}}{} \leq O(1) + O(u^{\alpha_\epsilon}) $ for $ r \simeq R $ as a consequence of \eqref{ODE:Lane-Emden001}. Inductively, $ \abs{u_{rrr}}{} \leq O(1) + O(u^{\alpha_\epsilon-1}) $, $ \abs{u_{rrrr}}{} \leq  O(1) + O(u^{\alpha_\epsilon - 2}) $, ... $ \abs{\dr^k u}{} \leq O(1) + O(u^{\alpha_\epsilon - k + 2}) $ for $ r\simeq R $.
\end{pf}

As a corollary of the homology invariance \eqref{hlinvrt}, we will construct a class of self-similar solutions describing the expanding or collapsing configuration for the radiational gaseous stars in the following.
\begin{lm}\label{lm:expanding}
For $ \iota = 0, ~ c_\nu = 3K $ and $ 1/6 < \epsilon K < 1 $, there is a globally expanding solution to \eqref{eq:dynm} with $ \Omega(t) = B_{a+bt}(0) $ for any positive constants $ a, ~ b > 0 $, where $ B_s(0) $ represents a ball centred at the origin with radius $ s $. Also, for any $ a > 0, b < 0 $, there is a collapsing solution with $ \Omega(t) = B_{a+bt}(0) $ and the collapsing time is given by $ \abs{a/b}{} $.
\end{lm}

\begin{pf}
	Given a regular solution $ (\bar\rho, \bar\theta, \bar\Omega) $ to \eqref{rdmd01}, it is spherically symmetric. Without loss of generality one can assume it is in the form $$ (\bar\rho, \bar\theta, \bar\Omega) = (\bar\rho(r), \bar\theta(r), \bar B_1), $$ with $ r = |x| $ and $ \bar B_1 $ being a ball centred at the origin with radius $ \bar R =1 $. 
	
	Define the ansatz $ (\rho(t, x), u(t,x), \theta(t,x), \Omega(t) )$ as follows.
	\begin{equation}\label{dy:ansatz}
		\begin{cases}
			\rho(t,x) = \alpha^3(t) \bar\rho ( \alpha (t) \cdot r ),\\
			\theta(t, x) = \alpha(t) \bar\theta ( \alpha(t) \cdot r ),\\
			u(t,x) = - \dfrac{\alpha'(t)}{\alpha(t)} x ,\\
			x \in \Omega(t) = \alpha^{-1}(t)\bar B_1 = \left\lbrace x \in \mathbb{R}^3 | 0 \leq |x| \leq \alpha^{-1}\right\rbrace,
		\end{cases}
	\end{equation}
	with $ r = |x| $ and $ \alpha(t) > 0 $ to be determined. Plug \eqref{dy:ansatz} into \eqref{eq:dynm}. Notice, the homology invariance \eqref{hlinvrt} implies the following identity
	\begin{equation}
		\begin{cases}
			 \nabla (K\rho \theta) = - \rho \nabla \psi & x\in \Omega(t), \\
			 - \Delta \theta = \epsilon \rho  & x \in \Omega(t), \\
			 \Delta \psi = \rho & x \in \mathbb{R}^3. \\
		\end{cases}
	\end{equation}
	Therefore, \eqref{eq:dynm} is reduced to
	\begin{equation}\label{dy:rdastz}
		\begin{cases}
			\rho_t + \dv (\rho u ) = 0 & x \in \Omega(t), \\
		(\rho u )_t + \dv (\rho u \otimes u ) = 0 & x \in \Omega(t), \\
		c_\nu (\rho \theta )_t + c_\nu \dv(\rho \theta u) + K \rho \theta \dv u = 0 & x \in \Omega(t), \\
		\end{cases}
	\end{equation}
	or equivalently
	\begin{equation}\label{dy:rdastz01}
		\begin{cases}
			\left( - \dfrac{\alpha'}{\alpha}\right)' + \left(\dfrac{\alpha'}{\alpha}\right)^2 = 0,\\
			c_\nu  \alpha'  - 3 K\alpha' = 0.
		\end{cases}
	\end{equation}
	When $ c_\nu = 3 K $, \eqref{dy:rdastz01} has an non-trivial solution satisfying
	\begin{displaymath}
		\alpha'' \alpha = 2 \alpha'^2,
	\end{displaymath}
	and therefore the non-trivial solution is
	\begin{equation}
		\alpha(t) = \dfrac{1}{a+b t}~~~~~~~~\text{for some constants $ a, ~ b $}.
	\end{equation}
	Moreover, it follows from the ansatz \eqref{dy:ansatz},
	\begin{equation}
		u (t,x ) = \dfrac{b}{a+bt} x ~~~ \text{and} ~~~ \Omega(t) = (a+bt)\bar B_1.
	\end{equation}
	This will finish the proof after choosing $ a,b > 0 $ or $ a>0, b<0 $.
\end{pf}

\section{A Prior Estimates}\label{sec:aprior}

In this section, we will establish the a prior estimate for the problem \eqref{eq:LagrangianCoordinates} with the initial density and the initial temperature satisfying \eqref{PV:lagrangian}. In the following, we denote $ \alpha = \frac{1-\epsilon K}{\epsilon K} $ and $ \sigma = R_0 - x $ for $ 0 \leq x\leq R_0 $ being the distance to the vacuum boundary. Then we have $ \rho_0 \simeq \sigma^\alpha $ and $ \Theta_0 \simeq \sigma $. Therefore, $ \rho_0, \Theta_0 $ satisfy
\begin{equation}\label{PV:lagrangianboundary}
	\abs{\dx^k\rho_0}{} \lesssim \sigma^{\alpha-k}, ~ \abs{\dx^k\Theta_0}{} \lesssim 1 + \sigma^{\alpha-k+2}.
\end{equation}
By denoting $$ \Phi = \dfrac{1}{x^3} \int_0^x s^2 \rho_0(s)\,ds, ~ \mathfrak{B} = (2\mu+\lambda)\dfrac{v_x}{r_x} + 2\lambda\dfrac{v}{r}, $$ we rewrite \eqref{eq:LagrangianCoordinates} as
\begin{equation}\label{eq:Lg000}
	\begin{cases}
		x^2 \rho_0 \dt v +  \bigl(K\dfrac{x^2\rho_0}{r_x} \Theta \bigr)_x - 2K x \rho_0 \Theta \dfrac{x}{r} = -  \dfrac{x^5\rho_0}{r^2} \Phi & \\ ~~~~~~~~~~ + r^2 \mathfrak{B}_x  
		+ 4 \mu r^2  \bigl(\dfrac{v}{r}\bigr)_x  & x \in (0, R_0), \\
		c_\nu x^2 \rho_0 \dt \Theta + K x^2 \rho_0 \Theta \dfrac{\dt(r^2 r_x)}{r^2 r_x} - \bigl(\dfrac{r^2}{r_x} \Theta_x\bigr)_x = \epsilon x^2 \rho_0 & \\ ~~~~~~~~~~ + 2\mu r^2 r_x \bigl( \bigl(\dfrac{v_x}{r_x}\bigr)^2 + 2 \bigl( \dfrac{v}{r}\bigr)^2 \bigr) + \lambda r^2 r_x \bigl( \dfrac{v_x}{r_x} + 2 \dfrac{v}{r} \bigr)^2 & x \in (0, R_0).
	\end{cases}
\end{equation}
Notice
\begin{equation}
	\norm{\Phi}{\supnorm} < \infty, ~ \left. \mathfrak B\right|_{x=R_0} = 0,
\end{equation}
and it is easy to verify that for a strong (or classical) solution, $ \Theta > 0 $ in $ (0,R_0) $. 

We will organize this section as follows. In Section \ref{sec:pointwise}, we obtain some point-wise estimates, which will manipulate the nonlinearity. In Section \ref{sec:elliptic}, we obtain the regularity estimates for the strong solutions in the spatial variable. In Section \ref{sec:energy}, we close the energy estimates for the strong solutions. In Section \ref{sec:classical}, we will perform the corresponding higher order estimates for the classical solutions. 

\subsection{Point-wise Estimates}\label{sec:pointwise}


In this section, the goal is to derive the bound the quantity 
\begin{equation}\label{APrioriAsum}
	\Lambda_0 = \Lambda_0(r,v,\Theta) : = \sup_{x\in(0,R_0)} \bigl\lbrace \abs{\dfrac{1}{r_x}}{},\abs{\dfrac{x}{r}}{},\sum_{i=0}^{1} \abs{\dt^i r_x}{},  \sum_{i=0}^{1} \abs{\dfrac{\dt^i r}{x}}{}, \abs{\dfrac{\Theta}{\sigma}}{}
		\bigr\rbrace.
\end{equation}
Also, we denote
\begin{equation} \label{APriorAsum2}
		M_0 = M_0 (r) := \sup_{x\in(0,R_0)} \bigl\lbrace \abs{r_x}{} ,\abs{\dfrac{r}{x}}{} , \abs{\dfrac{1}{r_x}}{}, \abs{\dfrac{x}{r}}{} \bigr\rbrace.
\end{equation}
In the following, it is assumed $ M_0, \Lambda_0 < \infty $.

\begin{lm}\label{lm:pointwise}
For a strong solution to \eqref{eq:LagrangianCoordinates} in the space $ \mathfrak X $ defined in \eqref{functionalspace:strong}, we have
	\begin{equation}\label{uniformese02oct}
		\norm{v_x}{\supnorm}, \norm{\dfrac{v}{x}}{\supnorm}, \norm{\dfrac{\Theta}{\sigma}}{\supnorm} \lesssim P(M_0) \bigl(\mathcal{E}_{1} + 1\bigr).
	\end{equation}
	Also,
	\begin{equation}\label{uniformese02oct02}
		\Lambda_0, \norm{x^{1/2} \Theta_x}{\supnorm} \lesssim P(M_0) \bigl(\mathcal{E}_1 + 1\bigr).
	\end{equation} 
	In particular, 
	\begin{equation}\label{uniformese01June}
	x^{1/2} \Theta_x \bigr|_{x=0} = 0.	
	\end{equation}

\end{lm}

\begin{pf}
For a strong solution in the space $ \mathfrak X $ defined in \eqref{functionalspace:strong}, the embedding theory yields,
$$ \norm{x{\Theta_x}}{L_t^\infty L_x^\infty(0,R_0/4)} \lesssim \norm{\Theta_x}{L^\infty_t L^2_x} + 
		\norm{ x \Theta_{xx}}{L^\infty_t L^2_x} < \infty. $$
Therefore, $ \frac{r^2}{r_x}\Theta_x\bigl|_{x=0} = 0 $. Then integrate \subeqref{eq:Lg000}{2} over $ x \in (0,x) $,
	\begin{align}\label{uniformese0000}
		& \dfrac{r^2}{r_x} \Theta_x = \int_{0}^{x} \bigl( \dfrac{r^2}{r_x} \Theta_x\bigr)_x \,dx = c_v \int_0^x x^2 \rho_0 \Theta_t \,dx + K \int_0^x \dfrac{x^2\rho_0}{r^2r_x} \Theta \bigl( r^2 v\bigr)_x \,dx {\nonumber} \\
		& ~~ - \epsilon \int_0^x x^2 \rho_0 \,dx - \int_0^x \biggl\lbrace 2\mu r^2 r_x \bigl( \bigl(\dfrac{v_x}{r_x}\bigr)^2 + 2 \bigl( \dfrac{v}{r}\bigr)^2 \bigr) + \lambda r^2 r_x \bigl( \dfrac{v_x}{r_x} + \dfrac{v}{r} \bigr)^2 \biggr\rbrace \,dx  {\nonumber} \\
		& : = \sum_{i=1}^{4} J_i^1.
	\end{align}
%
	By applying H\"older's inequality, we have the following estimates on the right,
	\begin{align*}
		& J_1^1 \leq c_\nu \biggl( \int_0^x x^2 \,dx \biggr)^{1/2} \biggl(\int_0^x x^2 \rho_0 \Theta_t^2 \biggr)^{1/2} \lesssim x^{3/2} \mathcal{E}_1^{1/2}, \\
		& J_2^1 = K \int_0^x \dfrac{x^2\rho_0}{r^2r_x} \Theta \bigl( r^2 v_x + 2 r r_x v\bigr)\,dx \\
		& ~~~~ \lesssim x^2 \norm{\dfrac{1}{r_x}}{\supnorm}  \biggl( \int  \Theta^2\,dx \biggr)^{1/2}\biggl( \int v_x^2\,dx \biggr)^{1/2} \\
		& ~~~~~~  + x^2 \norm{\dfrac{x}{r}}{\supnorm} \biggl( \int \Theta^2\,dx \biggr)^{1/2}\biggl(  \int \abs{\dfrac{v}{x}}{2}\,dx \biggr)^{1/2},\\
		& J_3^1 \lesssim x^3, \\
		& J_4^1 \lesssim x^2 \norm{\dfrac{r^2}{x^2r_x}}{\supnorm} \int v_x^2\,dx + x^2 \norm{r_x}{\supnorm} \int \abs{\dfrac{v}{x}}{2}\,dx.
	\end{align*}
	Then we have \begin{equation*}
		\dfrac{r^2}{r_x} \Theta_x = \sum_{i=1}^{4} J_i^1 \lesssim x^{3/2} P(M_0) \biggl\lbrace \biggl(\int_0^x x^2 \rho_0 \Theta_t^2 \,dx\biggr)^{1/2} + \int_0^x \bigl( \Theta^2 + v_x^2 + \abs{\frac{v}{x}}{2} \bigr) \,dx \biggr\rbrace.
	\end{equation*}
	By noticing
	\begin{align*}
		& \int v_x^2 \,dx, \int \abs{\dfrac{v}{x}}{2} \,dx \lesssim \int \chi v_x^2 \,dx + \int \chi \abs{\dfrac{v}{x}}{2}\,dx + \int x^2 v_x^2 \,dx + \int v^2 \,dx \lesssim \mathcal{E}_1, \\
		& \int \Theta^2\,dx \lesssim \int x^2 \Theta_x^2 \,dx \lesssim  \mathcal E_{1},
	\end{align*}
	the following estimate holds,
	\begin{equation}\label{uniformese0001}
	\begin{aligned}
		& x^{1/2} \Theta_x \lesssim P(M_0) \biggl\lbrace \biggl(\int_0^x x^2 \rho_0 \Theta_t^2 \,dx\biggr)^{1/2} + \int_0^x \bigl( \Theta^2 + v_x^2 + \abs{\frac{v}{x}}{2} \bigr) \,dx \biggr\rbrace \\
		& ~~~~ \lesssim P(M_0) ( \mathcal{E}_1 + 1 ).
	\end{aligned}
	\end{equation}
	In particular, by taking $ x \rightarrow 0^+ $, one can get \eqref{uniformese01June}.
	
	From \eqref{uniformese0001} we shall obtain the estimate of $ \norm{\frac{\Theta}{\sigma}}{\supnorm} $. On one hand, by taking $ x > R_0/4 $, \eqref{uniformese0001} yields the bound of $ \Theta_x $ in the corresponding area. Together with the boundary condition $ \Theta(R_0,t) = 0$, it yields the bound of $ \Theta/\sigma $ away from the centre. On the other hand, for $ 0 \leq x \leq R_0/2 $, 
	\begin{align*}
		& \abs{\Theta}{} \leq \int_x^{R_0} \abs{\Theta_x}{}\,dx \lesssim \norm{x^{1/2}\Theta_x}{\supnorm} \int_x^{R_0} x^{-1/2} \,dx \lesssim P(M_0) ( \mathcal{E}_1 + 1).
	\end{align*}
	Thus we have shown
	\begin{equation}\label{uniformese0002}
		\norm{\dfrac{\Theta}{\sigma}}{\supnorm} \leq P(M_0) ( \mathcal{E}_1 + 1 ).
	\end{equation}
	
On the other hand, integrate $ \frac{1}{r^2} \eqref{eq:Lg000}_{1} $ over the interval $ (x,R_0) $ . It holds,
	\begin{align*}
		(2\mu+\lambda)\dfrac{(r^2v)_x}{r^2r_x} = 4\mu \biggl. \dfrac{v}{r} \bigr|_{x=R_0} - \int_x^{R_0} \bigl(\dfrac{x}{r}\bigr)^2 \rho_0 v_t\,dx + K \dfrac{x^2\rho_0}{r^2r_x} \Theta - \int_x^{R_0} \dfrac{x^5\rho_0}{r^4} \Phi.
	\end{align*}
	Therefore, by using H\"older inequality,
	\begin{equation}\label{uniformese0003}
		\begin{aligned}
			& \abs{\dfrac{v_x}{r_x} + 2 \dfrac{v}{r}}{} \lesssim P(M_0) \biggl\lbrace \abs{v(R_0,t)}{} + \biggl(\int \rho_0 v_t^2 \,dx\biggr)^{1/2} + \abs{\Theta}{} + 1 \biggr\rbrace\\
			& ~~~~ \lesssim P(M_0) ( \mathcal{E}_1 + 1 ),
		\end{aligned}
	\end{equation}
	where it has been made used of \eqref{uniformese0002} and the facts
	\begin{align*}
		& \abs{v(R_0,t)}{2} \lesssim \int v^2 \,dx + \int x^2 v_x^2 \,dx \lesssim \mathcal E_1 ,\\
		& \int\rho_0v_t^2 \lesssim \int \chi \rho_0 v_t^2\,dx + C \int x^2 \rho_0 v_t^2 \,dx \lesssim \mathcal{E}_{1}.
	\end{align*}
	Next step is to multiply \subeqref{eq:Lg000}{1} with $ r $ and integrate the resulting over $ x \in (0,x) $. Then it follows,
	\begin{align*}
		& (2\mu+\lambda)\biggl\lbrace\int_0^x \biggl( \dfrac{(r^2v)_x}{r^2r_x} \biggr)_x r^3 \,dx \biggr\rbrace = \int_0^x x^2 r \rho_0 v_t\,dx + \int_0^x \bigl( K\dfrac{x^2 \rho_0}{r^2r_x}\Theta\bigr)_x r^3 \,dx \\
		& ~~~~+ \int_0^x \dfrac{x^5\rho_0\Phi}{r}\,dx = \int_0^x x^2 r\rho_0 v_t\,dx + K\dfrac{x^2\rho_0}{r^2r_x}\Theta r^3 - \int_0^x K \dfrac{x^2\rho_0}{r^2r_x}\Theta(r^3)_x\,dx \\
		& ~~~~ +  \int_0^x \dfrac{x^5\rho_0\Phi}{r}\,dx.
	\end{align*}
	Similar as before, it holds,
	\begin{align*}
		& \abs{\int_0^x \biggl( \dfrac{(r^2v)_x}{r^2r_x} \biggr)_x r^3 \,dx}{} \lesssim P(M_0) \biggl\lbrace r^3 \biggl(\int \rho_0 v_t^2 \,dx\biggr)^{1/2} + r^3 \norm{\Theta}{\supnorm} + r^3 \biggr\rbrace \\
		& ~~~~ \lesssim P(M_0) r^3(\mathcal{E}_{1} + 1),
	\end{align*}
	where the following is used
	\begin{align*}
		& \abs{\int_0^x K \dfrac{x^2\rho_0}{r^2r_x}\Theta(r^3)_x\,dx}{} \lesssim P(M_0) \norm{\Theta}{\supnorm} \abs{\int_0^x (r^3)_x \,dx}{} \lesssim P(M_0) r^3 \norm{\Theta}{\supnorm}.
	\end{align*}
	Moreover, direct calculation with integration by parts yields,
	\begin{align*}
		& \int_0^x \biggl( \dfrac{(r^2v)_x}{r^2r_x} \biggr)_x r^3 \,dx = \dfrac{(r^2v)_x}{r^2r_x} r^3 - \int_0^x \dfrac{(r^2v)_x}{r^2r_x} (r^3)_x\,dx \\
		& = r^3 \biggl( \dfrac{v_x}{r_x} + 2 \dfrac{v}{r} \biggr) - 3 \int_0^x (r^2v)_x \,dx =  r^3 \biggl( \dfrac{v_x}{r_x} + 2 \dfrac{v}{r} \biggr) - 3 r^2 v  = r^3 \biggl( \dfrac{v_x}{r_x} - \dfrac{v}{r} \biggr).
	\end{align*}
	Combining these calculations then gets
	\begin{equation}\label{uniformese0004}
		\abs{\dfrac{v_x}{r_x} - \dfrac{v}{r}}{}\lesssim P(M_0) (\mathcal{E}_{1} + 1).
	\end{equation}
	Together with \eqref{uniformese0003} then yields,
	\begin{equation*}
		\norm{v_x}{\supnorm}, \norm{\dfrac{v}{x}}{\supnorm} \lesssim P(M_0) (\mathcal{E}_{1}  + 1).
	\end{equation*}
\end{pf}

The following corollary is a direct consequence of the fact
\begin{equation*}
	v_x = \dt r_x, \frac{v}{x} = \dr \dfrac{r}{x},
\end{equation*}
and Lemma \ref{lm:pointwise}. 
\begin{cor}\label{cor:pointwise}
	For a strong solution to \eqref{eq:LagrangianCoordinates} in the space  $ \mathfrak X $, there is a time 
	\begin{equation} \label{cor:lowerboundoftime01}
		T_1 = T_1(\mathcal E_1) \geq \dfrac{C_1}{\mathcal E_1 + 1},
	\end{equation} 
	such that if $ T \leq T_1 $, 
	\begin{equation} \label{cor:upperbound01}
		M_0 \leq 2 ~~ \text{and} ~~ \Lambda_0 \leq C_1(\mathcal E_1+1), 
	\end{equation} 
	for some positive constant $ C_1 < \infty $.
\end{cor}
\begin{pf}
	From \eqref{uniformese02oct} and the definition of $ M_0 $ in \eqref{APriorAsum2}, we have
	\begin{equation*}
		\abs{r_x-1}{}, \abs{r/x-1}{} \leq T P_1(M_0)(\mathcal E_1 + 1).
	\end{equation*}
	some polynomial $ P_1 $. Let 
	\begin{equation}
		T_1 : = \dfrac{1}{2P_1(M_0)(\mathcal E_1 + 1)}.
	\end{equation}
	Then for $ T \leq T_1 $, we have $ \abs{r_x-1}{}, \abs{r/x-1}{} \leq 1/2 $ and $ M_0 \leq 2 $. Therefore $$ T_1 \geq \dfrac{1}{2P_1(2)(\mathcal{E}_1+1)}. $$
\end{pf}

\subsection{Elliptic Estimates}\label{sec:elliptic}

In this section, we aim to establish the regularity 
 \eqref{Regularity:Strong}. In other words, we will show that any solution to \eqref{eq:LagrangianCoordinates} with $ \mathcal E_1 < \infty $ will belong to the space $ \mathfrak X $.
 
We start with presenting the estimates for $ r, v $. This can be done using the technique in \cite{LuoXinZeng2016}. In fact, 
\subeqref{eq:Lg000}{1} can be rewritten as
\begin{equation}\label{eq:altformofmmteq}
		 (2\mu+\lambda) \mathcal G_{xt} + K \dfrac{x^2\rho_0 \Theta }{r^2r_x} \mathcal G_x = K \dfrac{x^2}{r^2r_x} \bigl( \rho_0 \Theta\bigr)_x  + \bigl( \dfrac{x}{r}\bigr)^2 \rho_0 v_t + \dfrac{x^5 \rho_0}{r^4} \Phi,
\end{equation}
where
\begin{equation}\label{def:relentp}
	\mathcal G = \ln \bigl( \dfrac{r^2r_x}{x^2}\bigr).
\end{equation}

\begin{lm} For a solution to \eqref{eq:altformofmmteq}, it holds
	\begin{align}
			& \dfrac{d}{dt} \int \mathcal G_{x}^2 \,dx + \int \dfrac{x^2\rho_0\Theta}{r^2r_x}\mathcal G_{x}^2 \,dx \lesssim \int \mathcal G_x^2 \,dx + P(\Lambda_0)\int \Theta_x^2\,dx {\nonumber}\\
			& ~~~~~~~~  + P(\Lambda_0) \bigl( \mathcal E_1 + 1 \bigr), {\label{ell:010}}\\
			& \int \mathcal G_{xt}^2 \,dx \lesssim P(\Lambda_0) \int \mathcal G_x^2 \,dx + P(\Lambda_0) (\mathcal E_1 + 1 ) + P(\Lambda_0) \int \Theta_x^2\,dx {\label{ell:020}} .
	\end{align}
\end{lm}
\begin{pf}
	Multiply \eqref{eq:altformofmmteq} with $ \mathcal G_{x} $ and integrate the resulting equation over $ x\in (0, R_0) $ in the following,
	\begin{equation}\label{ell:001}
		\begin{aligned}
			& \dfrac{d}{dt}\biggl\lbrace \dfrac{2\mu+\lambda}{2}\int \mathcal G_x^2 \,dx \biggr\rbrace + K \int \dfrac{x^2\rho_0\Theta}{r^2r_x} \mathcal G_x^2 \,dx = K \int \dfrac{x^2}{r^2r_x}(\rho_0\Theta)_x \mathcal G_x \,dx \\
			& ~~~~~~ + \int \bigl( \dfrac{x}{r}\bigr)^2 \rho_0v_t \mathcal G_x \,dx + \int \dfrac{x^5\rho_0}{r^4}\Phi \mathcal G_{x}\,dx \lesssim  \int \mathcal G_x^2\,dx \\
			& ~~~~ + \int \biggl( \dfrac{x^2}{r^2r_x} (\rho_0\Theta)_x \biggr)^2 \,dx + \norm{\dfrac{x}{r}}{\supnorm}^4 \int \rho_0 v_t^2 \,dx + \norm{\dfrac{x^4}{r^4}}{\supnorm}^2 \int x^2 \rho_0\,dx.
		\end{aligned}
	\end{equation}
	In the meantime, by noticing
	\begin{align*}
		\abs{(\rho_0\Theta)_x}{} \lesssim \abs{(\rho_0)_x \Theta}{} + \abs{\rho_0 \Theta_x}{}\lesssim \sigma^{\alpha-1} \abs{\Theta}{} + \abs{\Theta_x}{},
	\end{align*}
	as a consequence of \eqref{PV:lagrangianboundary}, 
	\begin{align*}
	& \int \biggl( \dfrac{x^2}{r^2r_x} (\rho_0\Theta)_x \biggr)^2 \,dx \lesssim \norm{\dfrac{x^2}{r^2r_x}}{\supnorm}^2 \int \sigma^{2\alpha-2} \Theta^2 \,dx + \norm{\dfrac{x^2}{r^2r_x}}{\supnorm}^2 \int \Theta_x^2\,dx.
	\end{align*}
	For $ \alpha \neq 1/2 $, by applying the inequalities a), b) in Lemma \ref{lm:hardy},
	\begin{equation*}
		\int \sigma^{2\alpha-2} \Theta^2 \,dx \lesssim \int x^2 \sigma^{2\alpha} \Theta_x^2 \,dx \lesssim \int x^2 \Theta_x^2 \,dx.
	\end{equation*}
	For $ \alpha = 1/2 $, 
	\begin{equation*}
		\int \sigma^{2\alpha-2} \Theta^2 \,dx = \int \sigma^{-1}\Theta^2\,dx \lesssim \norm{\dfrac{\Theta}{\sigma}}{\supnorm}^2 + \int \Theta^2 \,dx \lesssim \Lambda_0^2 + \int x^2 \Theta_x^2\,dx.
	\end{equation*}
	Therefore, \eqref{ell:010} is a direct consequence of \eqref{ell:001} and the fact
	\begin{equation*}
		\int \rho_0 v_t^2 \, dx \lesssim \int \chi \rho_0 v_t^2 \,dx + \int x^2 \rho_0 v_t^2 \,dx \lesssim \mathcal E_1 .
	\end{equation*}
	Meanwhile, from \eqref{eq:altformofmmteq}
	\begin{align*}
		& \int \mathcal G_{xt}^2 \,dx \lesssim 
		P(\Lambda_0) \int \mathcal G_{x}^2 \,dx + P(\Lambda_0) + P(\Lambda_0) \int x^2 \Theta_x^2 \,dx \\
		& ~~~~~~~~~~+ P(\Lambda_0) \int \rho_0 v_t^2 \, dx + P(\Lambda_0) \int \Theta_x^2\,dx,
	\end{align*}
	and thus \eqref{ell:020} holds. 	
\end{pf}

The next lemma is considering the regularity estimates for $ \Theta $.
\begin{lm}\label{lm:L2-theta-xx}
For a solution to \subeqref{eq:Lg000}{2}, it holds 
	\begin{align}
		& \int \chi \dfrac{r^2}{r_x} \biggl( \bigl(\dfrac{\Theta_x}{r_x}\bigr)_x\biggr)^2 \,dx + \int \chi r_x \bigl( \dfrac{\Theta_x}{r_x}\bigr)^2 \,dx \lesssim P(\Lambda_0) ( \mathcal E_{1}+1 ), {\label{ell:030}} \\
		& \int \biggl( \dfrac{1}{r\sqrt{r_x}} \bigl( \dfrac{r^2}{r_x} \Theta_x\bigr)_x  \biggr)^2 \,dx  \lesssim P(\Lambda_0) ( \mathcal E_{1} + 1 ). {\label{ell:040}}
	\end{align}
\end{lm}

\begin{pf}
	By multiplying \subeqref{eq:Lg000}{2} with $ (r\sqrt{r_x})^{-1} $, it holds the following 
	\begin{align*}
		& \dfrac{1}{r\sqrt{r_x}} \bigl( \dfrac{r^2}{r_x} \Theta_x\bigr)_x = c_\nu \dfrac{x^2 \rho_0}{r\sqrt{r_x}} \Theta_t + K \dfrac{x^2 \rho_0}{r\sqrt{r_x}} \Theta \dfrac{(r^2v)_x}{r^2 r_x} \\
		& ~~~~~~ - r \sqrt{r_x} \bigl\lbrack 2\mu \bigl( \bigl(\dfrac{v_x}{r_x}\bigr)^2 + 2 \bigl( \dfrac{v}{r}\bigr)^2 \bigr) + \lambda \bigl( \dfrac{v_x}{r_x} + \dfrac{v}{r} \bigr)^2 \bigr\rbrack - \epsilon \dfrac{x^2 \rho_0}{r\sqrt{r_x}}.
	\end{align*}
	Thus, direct calculation yields
	\begin{align*}
		& \int \chi \biggl( \dfrac{1}{r\sqrt{r_x}} \bigl( \dfrac{r^2}{r_x} \Theta_x\bigr)_x  \biggr)^2 \,dx \lesssim \int \biggl( \dfrac{1}{r\sqrt{r_x}} \bigl( \dfrac{r^2}{r_x} \Theta_x\bigr)_x  \biggr)^2 \,dx \\
		& ~~~~ \lesssim \norm{\dfrac{x^2}{r^2r_x}}{\supnorm} \int x^2 \rho_0 \Theta_t^2 \,dx + \bigl( \norm{\dfrac{x^2 v^2}{r^4 r_x}}{\supnorm} + \norm{\dfrac{x^2 v_x^2}{r^2 r_x^3}}{\supnorm} \bigr) \int x^2 \rho_0 \Theta^2 \,dx \\
		& ~~~~+ \norm{\dfrac{x^2}{r^2r_x}}{\supnorm} \int x^2 \rho_0\,dx + \norm{\dfrac{r^2 v_x^2}{x^2 r_x^3}}{\supnorm} \int x^2 v_x^2 \,dx\\
		& ~~~~ + \norm{\dfrac{v^2 r_x}{r^2}}{\supnorm} \int x^2 \bigl(\dfrac{v}{x}\bigr)^2 \, dx \lesssim P(\Lambda_0) (\mathcal E_{1} + 1 ).
	\end{align*}
	Meanwhile, by applying integration by parts and noticing $ \chi' \leq 0 $,
	\begin{align*}
		& \int \chi \biggl( \dfrac{1}{r\sqrt{r_x}} \bigl( \dfrac{r^2}{r_x} \Theta_x\bigr)_x  \biggr)^2 \,dx = \int \chi \dfrac{r^2}{r_x} \biggl( \bigl(\dfrac{\Theta_x}{r_x}\bigr)_x\biggr)^2 \,dx + 4 \int \chi r_x \bigl( \dfrac{\Theta_x}{r_x}\bigr)^2 \,dx \\
		& ~~~~~~ + 4 \int \chi r \bigl(\dfrac{\Theta_x}{r_x}\bigr) \bigl( \dfrac{\Theta_x}{r_x}\bigr)_x\,dx = \int \chi \dfrac{r^2}{r_x} \biggl( \bigl(\dfrac{\Theta_x}{r_x}\bigr)_x\biggr)^2 \,dx + 2 \int \chi r_x \bigl( \dfrac{\Theta_x}{r_x}\bigr)^2 \,dx \\
		& ~~~~~~ - 2 \int \chi' r \bigl(\dfrac{\Theta_x}{r_x}\bigr)^2 \,dx - 2 \underbrace{ \chi r \bigl(\dfrac{\Theta_x}{r_x}\bigr)^2\biggr|_{x=0} }_{=0\footnotemark} \\
		& ~~~~ \geq \int \chi \dfrac{r^2}{r_x} \biggl( \bigl(\dfrac{\Theta_x}{r_x}\bigr)_x\biggr)^2 \,dx + 2 \int \chi r_x \bigl( \dfrac{\Theta_x}{r_x}\bigr)^2 \,dx.
	\end{align*}
	\footnotetext{See \eqref{uniformese01June}}
	Combining the above inequalities to finish the proof.
\end{pf}

Now we will combine the two lemmas above to derive the desired regularity estimates.
\begin{lm}[Elliptic Estimates]\label{lm:EllipticRegularity} For a solution to \eqref{eq:Lg000} and $ 0 < T < 1 $, it holds,
	\begin{equation}\label{ell:100}
		\begin{gathered}
			\norm{r_{xx}}{\stnorm{\infty}{2}}, \norm{\bigl(\dfrac{r}{x}\bigr)_x}{\stnorm{\infty}{2}}, \norm{v_{xx}}{\stnorm{\infty}{2}}, \norm{\bigl(\dfrac{v}{x}\bigr)}{\stnorm{\infty}{2}}, \\
			\norm{\Theta_x}{\stnorm{\infty}{2}}, \norm{x\Theta_{xx}}{\stnorm{\infty}{2}} \lesssim P(\Lambda_0)(\mathcal E_1 + 1).
		\end{gathered}
	\end{equation}
	This will show that the solution with $ \mathcal E_1 < \infty $ is in the space $ \mathfrak X $. 
\end{lm}

\begin{pf}
	From \eqref{ell:030}, it holds
	\begin{equation}\label{ell:070}
		\begin{aligned}
			& \int \Theta_x^2 \,dx \lesssim P(\Lambda_0) \int \chi r_x \bigl( \dfrac{\Theta_x}{r_x}\bigr)^2 \,dx + \int x^2 \Theta_x^2\,dx \lesssim P(\Lambda_0) (\mathcal E_1 + 1 ) .
		\end{aligned}
	\end{equation}
	Therefore \eqref{ell:010} implies
	\begin{align*}
		& \dfrac{d}{dt} \int \mathcal G_{x}^2 \,dx + \int \dfrac{x^2 \rho_0 \Theta}{r^2 r_x} \mathcal G_x^2 \,dx \lesssim \int \mathcal G_x^2 \,dx + P(\Lambda_0) ( \mathcal{E}_1 + 1 ).
	\end{align*}
	Notice, from the definition of $ \mathcal G $, it holds $ \mathcal G_x(x, 0) = 0 $. Gr\"onwall's inequality then yields, for $ 0 < T < 1 $, 
	\begin{equation}\label{ell:050}
		\int\mathcal G_{x}^2 \,dx  \lesssim T \exp(c T) P(\Lambda_0) ( \mathcal E_1 + 1) \lesssim P(\Lambda_0) ( \mathcal E_1 + 1) .
	\end{equation}
	Also, from \eqref{ell:020}, \eqref{ell:070}, \eqref{ell:050},
	\begin{equation}\label{ell:060}
		\int \mathcal G_{xt}^2\,dx \lesssim P(\Lambda_0) ( \mathcal E_1 + 1).
	\end{equation}
	From \eqref{ell:050}, \eqref{ell:060}, we shall derive the bound of $ r_{xx}, v_{xx}, \bigl(\frac{r}{x}\bigr)_x, \bigl(\frac{v}{x}\bigr)_x $ in $ L^2 $ space. 
	Notice
	\begin{align*}
		& \mathcal G_x = \dfrac{x^2}{r^2r_x} \biggl( \dfrac{r^2r_x}{x^2} \biggr)_x= \dfrac{x}{rr_x} \biggl( 2 r_x \bigl( \dfrac{r}{x} \bigr)_x + \dfrac{r}{x} r_{xx}  \biggr), \\
		& \mathcal G_{xt} = \dfrac{x}{rr_x} \biggl( 2 r_x \bigl(\dfrac{v}{x}\bigr)_x + v_{xx}\bigl(\dfrac{r}{x}\bigr) + 2 v_x\bigl(\dfrac{r}{x}\bigr)_x + r_{xx}\bigl(\dfrac{v}{x}\bigr) \biggr)\\
		& ~~~~~~ - \dfrac{x}{rr_x} \biggl( \dfrac{v_x}{r_x} + \dfrac v r \biggr)\biggl( 2 r_x \bigl(\dfrac r x \bigr)_x + r_{xx}\bigl(\dfrac r x \bigr) \biggr).
	\end{align*}
	Thus \eqref{ell:050} implies 
	\begin{align*}
		& P(\Lambda_0)(\mathcal E_1 + 1) \gtrsim  \int \biggl( 2 r_x \bigl( \dfrac{r}{x} \bigr)_x + \dfrac{r}{x} r_{xx} \biggr)^2 \,dx = \int r_{xx}^2 \bigl(\dfrac{r}{x}\bigr)^2 \,dx \\
		& ~~~~ + \int 4 r_x r_{xx} \bigl(\dfrac{r}{x}\bigr)\bigl(\dfrac{r}{x}\bigr)_x + 4 r_x^2 \biggl(\bigl(\dfrac{r}{x}\bigr)_x\biggr)^2 \,dx : = A+B.
	\end{align*}
	After rewriting 
	\begin{align*}
		r_x = x\bigl(\dfrac{r}{x}\bigr)_x + \dfrac{r}{x}, ~ r_{xx} = x \bigl(\dfrac{r}{x}\bigr)_{xx} + 2 \bigl(\dfrac{r}{x}\bigr)_x,
	\end{align*}
	applying integration by parts to $ B $ then yields,
	\begin{align*}
		& B = \biggl. \biggl( 2 x \bigl(\dfrac{r}{x}\bigr)^2 \biggl( \bigl(\dfrac{r}{x}\bigr)_x \biggr)^2 + \dfrac{4}{3} x^2 \bigl(\dfrac{r}{x}\bigr) \biggl(\bigl(\dfrac{r}{x}\bigr)_x\biggr)^3 \biggr) \biggr|_{x=R_0} + 10 \int \bigl(\dfrac{r}{x}\bigr)^2 \biggl(\bigl(\dfrac{r}{x}\bigr)_x\biggr)^2\,dx \\
		& ~~~ + \dfrac{28}{3}\int x \bigl( \dfrac{r}{x} \bigr) \biggl(\bigl(\dfrac{r}{x}\bigr)_x\biggr)^3\,dx + \dfrac{8}{3} \int x^2 \biggl(\bigl(\dfrac{r}{x}\bigr)_x\biggr)^4\,dx \\
		& = \biggl. \bigl( 2 x \biggl(\dfrac{r}{x}\bigr)^2 \biggl( \bigl(\dfrac{r}{x}\bigr)_x \biggr)^2 + \dfrac{4}{3} x^2 \bigl(\dfrac{r}{x}\bigr) \biggl(\bigl(\dfrac{r}{x}\bigr)_x\biggr)^3 \biggr) \biggr|_{x=R_0} \\
		& ~~~ + \dfrac{2}{3} \int\biggl(2x \biggl(\bigl(\dfrac{r}{x}\bigr)_x\biggr)^2 +\dfrac{7}{2}\dfrac{r}{x}\bigl(\dfrac{r}{x}\bigr)_x\biggr)^2 \,dx + \dfrac{11}{6} \int \bigl(\dfrac{r}{x}\bigr)^2\biggl(\bigl(\dfrac{r}{x}\bigr)_x\biggr)^2\,dx.
	\end{align*}
	Meanwhile, since
	\begin{equation*}
		\bigl(\dfrac{r}{x}\bigr)_x = \dfrac{r_x}{x} - \dfrac{r}{x^2},
	\end{equation*}
	the boundary term in $ B $ can be bounded by $ P(\Lambda_0) $. Combining these calculation then yields,
	\begin{equation}\label{ell:090}
		\int r_{xx}^2 \,dx + \int \biggl(\bigl(\dfrac{r}{x}\bigr)_x\biggr)^2\,dx \lesssim P(\Lambda_0) (\mathcal{E}_{1} + 1). 
	\end{equation}
	Similarly, \eqref{ell:060} implies
	\begin{align*}
		& P(\Lambda_0) (\mathcal{E}_{1} + 1) \gtrsim \int \biggl( 2 r_x \bigl(\dfrac{v}{x}\bigr)_x + v_{xx}\bigl(\dfrac{r}{x}\bigr)\biggr)^2 \,dx \\ 
		& ~~~ = \int 4 r_x^2 \biggl( \bigl(\dfrac v x \bigr)_x\biggr)^2 + \bigl(\dfrac r x \bigr)^2 v_{xx}^2 \,dx  + 4 \int r_x \bigl(\dfrac{r}{x}\bigr)\bigl(\dfrac v x \bigr)_x v_{xx}\,dx : = C + D.
	\end{align*}
	Again,
	\begin{align*}
		& D = 4 \int r_x \bigl( \dfrac{r}{x} \bigr) \bigl(\dfrac{v}{x}\bigr)_x \biggl( x \bigl(\dfrac{v}{x}\bigr)_{xx} + 2 \bigl(\dfrac v x\bigr)_x\biggr)\,dx\\
		& ~~~ = \biggl. \biggl( 2 x r_x \bigl( \dfrac r x \bigr) \biggl( \bigl( \dfrac v x \bigr)_x \biggr)^2 \biggr) \biggr|_{x=R_0} + 6 \int r_x \bigl( \dfrac r x \bigr) \biggl( \bigl( \dfrac v x \bigr)_x\biggr)^2\,dx\\
		& ~~~~~ - 2 \int x r_{xx} \bigl(\dfrac r x \bigr) \biggl(\bigl(\dfrac v x \bigr)_x\biggr)^2 \,dx - 2 \int x r_x \bigl( \dfrac r x \bigr)_x \biggl( \bigl( \dfrac v x \bigr)_x \biggr)^2\,dx,
	\end{align*}
	where the boundary term can be bounded by $ P(\Lambda_0) $ due to the fact
	\begin{equation*}
		\bigl(\dfrac{v}{x}\bigr)_x = \dfrac{v_x}{x} - \dfrac{v}{x^2}. 
	\end{equation*}
	Therefore, we have the following
	\begin{equation}\label{ell:080}
		\begin{aligned}
			& \int \bigl(\dfrac r x \bigr)^2 v_{xx}^2 \,dx + \int\biggl(4 r_x^2 + 6 r_x \bigl(\dfrac r x \bigr) - 2 x r_x \bigl( \dfrac r x \bigr)_x \biggr. \\& ~~~~~~ \biggl.  - 2 x r_{xx}\bigl(\dfrac r x \bigr) \biggr) \biggl( \bigl( \dfrac v x \bigr)_x\biggr)^2\,dx \lesssim P(\Lambda_0) (\mathcal{E}_{1} + 1 ).
		\end{aligned}
	\end{equation}
	With a direct calculation, we have the following inequality concerning the coercivity of coefficients
	\begin{equation}\label{coecivity}
		\begin{aligned}
			& 4 r_x^2 + 6 r_x \bigl(\dfrac r x \bigr) - 2 x r_x \bigl( \dfrac r x \bigr)_x - 2 x r_{xx}\bigl( \dfrac r x \bigr) = 10 \bigl( \dfrac r x \bigr)^2 + 2  x^2 \biggl(\bigl(\dfrac r x \bigr)_x\biggr)^2 \\
			& ~~~~ +12 x \bigl( \dfrac r x \bigr) \bigl( \dfrac{r}{x}\bigr)_x - 2 x r_{xx} \bigl( \dfrac r x \bigr) \geq 9 \bigl( \dfrac r x \bigr)^2 - C x^2 \biggl(\bigl(\dfrac r x \bigr)_x\biggr)^2 - C x^2 r_{xx}^2.
		\end{aligned}
	\end{equation}
	Hence from \eqref{ell:080}, 
	\begin{equation}
		\begin{aligned}
			& \int \bigl(\dfrac r x \bigr)^2 v_{xx}^2 \,dx + 9 \int \bigl(\dfrac r x \bigr)^2 \biggl( \bigl(\dfrac v x \bigr)_x\biggr)^2 \,dx \lesssim P(\Lambda_0) (\mathcal{E}_{1} + 1)\\
			& ~~~~~~~~ + C \int \biggl(r_{xx}^2 + \bigl(\dfrac{r}{x}\bigr)_x^2 \biggr) x^2 \biggl(\bigl(\dfrac v x \bigr)_x\biggr)^2\,dx \\
			& ~~~~ \lesssim  P(\Lambda_0) (\mathcal{E}_{1} + 1) + P(\Lambda_0) \int r_{xx}^2 + \biggl( \bigl( \dfrac r x \bigr)_x\biggr)^2\,dx \lesssim P(\Lambda_0) (\mathcal{E}_{1} + 1),
		\end{aligned}
	\end{equation}
	where it has been applied the fact
	\begin{equation*}
		x^2 \biggl( \bigl( \dfrac v x \bigr)_x\biggr)^2 = \biggl( v_x - \dfrac{v}{x} \biggr)^2 \lesssim P(\Lambda_0).
	\end{equation*}
	Combining these estimate then yields,
	\begin{equation}
			\int v_{xx}^2 \,dx + \int \biggl(\bigl(\dfrac{v}{x}\bigr)_x\biggr)^2\,dx \lesssim P(\Lambda_0) ( \mathcal E_1 + 1 ). 
	\end{equation}
	Moreover, from \eqref{ell:040}, \eqref{ell:070}, \eqref{ell:090},
	\begin{align*}
		& \int x^2 \Theta_{xx}^2 \,dx \lesssim P(\Lambda_0) \biggl( \int \Theta_x^2 \,dx + \int x^2 \Theta_x^2 r_{xx}^2 \,dx + \int \biggl( \dfrac{1}{r\sqrt{r_x}} \bigl(\dfrac{r^2}{r_x}\Theta_x\bigr)_x\biggr)^2\,dx \biggr) \\
		 & ~~~ \lesssim P(\Lambda_0) (\mathcal E_1 + 1) + P(\Lambda_0) \norm{x\Theta_x}{\supnorm}^2 \int r_{xx}^2\,dx \lesssim P(\Lambda_0) ( \mathcal E_1 + 1 ).
	\end{align*}
	Thus we finish the proof of the lemma.
\end{pf}

We have finished the elliptic estimates for the strong solutions to \eqref{eq:Lg000}. The energy estimates will be continued in the following and will close the a prior estimates.


\subsection{Energy Estimates}\label{sec:energy}

In this section, we will work on the energy estimates on \eqref{eq:Lg000} (or equivalently \eqref{eq:LagrangianCoordinates} ). In particular, we write down the temporal derivative version of \eqref{eq:Lg000} in the following. 
\begin{equation} \label{eq:Lg010}
	\begin{cases}
		x^2 \rho_0 \dt v_t - r^2 \mathfrak B_{xt} - 4\mu r^2 \bigl( \dfrac{v_t}{r} \bigr)_x = - \bigl( K x^2 \rho_0 I_1^1 \bigr)_x + 2K x\rho_0 I_2^1\\
		~~~~~~~~  + x^2 I_3^1 + x^2 I_4^1 + x^2 \rho_0 I_5^1, \\
		c_v x^2 \rho_0 \dt\Theta_t -  \bigl( \dfrac{r^2}{r_x} \Theta_{xt} \bigr)_x = - K x^2 \rho_0 I_6^1 + x^2 I_7^1 -  \bigl( x^2 I_8^1 \Theta_x \bigr)_x,
	\end{cases}
\end{equation}
where
\begin{align}
		& I_1^1 = \bigl( \dfrac{\Theta}{r_x} \bigr)_t, ~ I_2^1 =  \bigl( \Theta \dfrac{x}{r}\bigr)_t ,~ I_3^1 =  \bigl(\dfrac{r^2}{x^2} \mathfrak B_x \bigr)_t - \dfrac{r^2}{x^2} \mathfrak B_{xt}, {\nonumber} \\
		& I_4^1 = \biggl( 4\mu \dfrac{r^2}{x^2}  \bigl( \dfrac{v}{r}\bigr)_x \biggr)_t - 4\mu\dfrac{r^2}{x^2} \bigl( \dfrac{v_t}{r}\bigr)_x, ~ I_5^1 = - \bigl(\dfrac{x^3}{r^2} \Phi \bigr)_t, {\nonumber}\\
		& I_6^1 =  \biggl( \Theta \dfrac{(r^2r_x)_t}{r^2r_x} \biggr)_t, ~ \\
		& I_7^1 = \biggl( 2\mu \dfrac{r^2 r_x}{x^2} \biggl( \bigl(\dfrac{v_x}{r_x}\bigr)^2 + 2 \bigl( \dfrac{v}{r}\bigr)^2 \biggr) + \lambda \dfrac{r^2 r_x}{x^2} \biggl( \dfrac{v_x}{r_x} + 2 \dfrac{v}{r} \biggr)^2 \biggr)_t, {\nonumber}\\
		& I_8^1 = - \dfrac{2rv}{x^2r_x}+\dfrac{r^2 v_x}{x^2r_x^2}. {\nonumber}
\end{align}
The following boundary condition is a consequence of differentiating \eqref{boundarycndtnLG} in the temporal variable,
\begin{equation}\label{boundarycndtnLG1}
	v_t(0,t) = 0, ~ \Theta_t(R_0,t) = 0, ~ \mathfrak B_t (R_0,t) = 0, ~ t \geq 0.
\end{equation}
Also, it is assumed $ \mathcal E_1^0 < \infty $ defined in \eqref{initialenergy-strong} through out this section. 

Now we will derive the estimates away from the center $ x = 0 $ in the next lemma.

\begin{lm}\label{lm:energyest-strong-001}
For a solution to \eqref{eq:Lg000} (or equivalently \eqref{eq:LagrangianCoordinates}) and \eqref{eq:Lg010}, it holds
\begin{align}
		& \norm{x\sqrt{\rho_0} v_t}{\stnorm{\infty}{2}}^2 + \norm{x \sqrt{\rho_0}\Theta_t}{\stnorm{\infty}{2}}^2 + \norm{x v_{xt}}{\stnorm{2}{2}}^2 + \norm{v_t}{\stnorm{2}{2}}^2 {\nonumber} \\
		& ~~~~~~ + \norm{x \Theta_{xt}}{\stnorm{2}{2}}^2 \lesssim T P(\Lambda_0) (\mathcal E_1 + 1) + P(M_0) \mathcal E_1^0, {\label{ene:1st-010}} \\
		&  \norm{xv_x}{\stnorm{\infty}{2}}^2  + \norm{v}{\stnorm{\infty}{2}}^2 + \norm{x \Theta_x }{\stnorm{\infty}{2}}^2 + \norm{x\sqrt{\rho_0}v_t}{\stnorm{2}{2}}^2 {\nonumber} \\
		& ~~~~~~ + \norm{x\sqrt{\rho_0}\Theta_t}{\stnorm{2}{2}}^2 \lesssim TP(\Lambda_0) (\mathcal E_1 + 1 ) + P(M_0) \mathcal E_1^0.  {\label{ene:1st-020}}
\end{align}
\end{lm}
\begin{pf}
	Multiply \subeqref{eq:Lg010}{1} with $ v_t $ and \subeqref{eq:Lg010}{2} with $ \Theta_t $ respectively and integrate the resulting equations over $ x \in (0, R_0) $. After integration by parts, it holds the following,
	\begin{align}
		& \dfrac{d}{dt} \biggl\lbrace \dfrac{1}{2} \int x^2 \rho_0 v_t^2 \,dx \biggr\rbrace + \int \mathfrak B_t \bigl( r^2 v_t\bigr)_x\,dx - \int 4\mu r^2 \bigl(\dfrac{v_t}{r}\bigr)_x v_t\,dx {\nonumber}\\
		& ~~~~~~ = K \int x^2 \rho_0 I_1^1 v_{xt} \,dx + 2 K \int x\rho_0 I_2^1 v_t \,dx + \int x^2 I_3^1 v_t\,dx {\nonumber} \\
		& ~~~~~~~~~~ + \int x^2 I_4^1 v_t\,dx + \int x^2 \rho_0 I_5^1 v_t\,dx : = \sum_{i=1}^{5} J_i^2,  {\label{ene:1st-001}} \\
		& \dfrac{d}{dt} \biggl\lbrace \dfrac{c_\nu}{2} \int x^2 \rho_0 \Theta_t^2 \,dx  \biggr\rbrace + \int \dfrac{r^2}{r_x} \Theta_{xt}^2\,dx = - K \int x^2 \rho_0 I_6^1 \Theta_t\,dx {\nonumber} \\
		& ~~~~~~~~~~ + \int x^2 I_7^1 \Theta_t\,dx + \int x^2 I_8^1 \Theta_x \Theta_{xt} : = \sum_{i=6}^8 J_i^2. {\label{ene:1st-002}}
	\end{align}
	To evaluate the left of \eqref{ene:1st-001}, notice
	\begin{equation*}
		\begin{aligned}
			& \mathfrak B_t (r^2 v_t)_x - 4\mu r^2 \bigl(\dfrac{v_t}{r} \bigr)_x v_t = r^2 r_x \biggl\lbrace 2\mu \bigl( \dfrac{v_{xt}^2}{r_x^2}  + 2 \dfrac{v_t^2}{r^2} \bigr) + \lambda \bigl( \dfrac{v_{xt}}{r_x}+2\dfrac{v_t}{r} \bigr)^2 \biggr\rbrace \\
			& ~~~~~~ - \bigl( (2\mu+\lambda)\dfrac{v_x^2}{r_x^2} + 2\lambda\dfrac{v^2}{r^2} \bigr) \bigl( v_{xt} r^2 + 2 v_t r r_x \bigr).
		\end{aligned}
	\end{equation*}
	Thus \eqref{ene:1st-001}, \eqref{ene:1st-002} can be written as
	\begin{align*}
		& \dfrac{d}{dt} \biggl\lbrace\dfrac{1}{2} \int x^2 \rho_0 v_t^2 \,dx \biggr\rbrace +  \int r^2 r_x \biggl\lbrace 2\mu \bigl( \dfrac{v_{xt}^2}{r_x^2}  + 2 \dfrac{v_t^2}{r^2} \bigr) + \lambda \bigl( \dfrac{v_{xt}}{r_x}+2\dfrac{v_t}{r} \bigr)^2 \biggr\rbrace \,dx \\
		& ~~~~~  = \sum_{i=1}^5 J_i^2 + J_9^2, \\
		& \dfrac{d}{dt} \biggl\lbrace\dfrac{c_\nu}{2}\int x^2 \rho_0 \Theta_t^2 \,dx \biggr\rbrace + \int \dfrac{r^2}{r_x} \Theta_{xt}^2 \,dx = \sum_{i=6}^{8} J_i^2,
	\end{align*}
	with
	\begin{equation*}
		J_9^2 =  \int \bigl( (2\mu+\lambda)\dfrac{v_x^2}{r_x^2} + 2\lambda\dfrac{v^2}{r^2} \bigr) \bigl( v_{xt} r^2 + 2 v_t r r_x \bigr) \,dx.
	\end{equation*}
	To evaluate the $ J_i^2 $'s, by applying the Cauchy's inequality, it holds
	\begin{align*}
		& J_1^2 + J_2^2 + J_4^2 + J_5^2 + J_9^2 \lesssim \delta \int \bigl( x^2 v_{xt}^2 + v_t^2\bigr) \,dx  + C_\delta P(\Lambda_0) \bigl( \mathcal E_{1} + 1 ),  \\
		& J_6^2 + J_7^2 + J_8^2 \lesssim \delta \int x^2 \Theta_{xt}^2 + \int\bigl( x^2 v_{xt}^2 + v_t^2 \bigr) \,dx + P(\Lambda_0) \int x^2 \Theta_t^2\,dx \\
		& ~~~~~~ + C_\delta P(\Lambda_0) ( \mathcal E_{1} + 1).
	\end{align*}
	For $ 0 < \omega < R_0/2 $, by applying $ c), d) $ in Lemma \ref{lm:hardy}, 
	\begin{equation}\label{ene:1st-003}
	\begin{aligned}
		& \int x^2 \Theta_t^2 \,dx = \int_0^{R_0-\omega} x^2 \Theta_t^2 \,dx + \int_{R_0-\omega}^{R_0} x^2 \Theta_t^2 \,dx \lesssim \omega^{-\alpha} \int_0^{R_0-\omega} x^2 \rho_0 \Theta_t^2\,dx \\
		& ~~~~~~ + \omega^2 \int_{R_0-\omega}^{R_0} x^2 \Theta_{xt}^2 \,dx.
	\end{aligned}\end{equation}
	Therefore, after choosing $ \omega $ so that $ \omega^2 P(\Lambda_0) \simeq \delta $, we have
	\begin{equation*}
		P(\Lambda_0) \int x^2 \Theta_t^2 \,dx \lesssim \delta \int x^2 \Theta_{xt}^2 \,dx + P(\Lambda_0) \mathcal E_1.
	\end{equation*}
	On the other hand, by applying integration by parts and Cauchy's inequality, 
	\begin{align*}
		& J_3^2 
		= - 2 \int \mathfrak B \left( r v v_t \right)_x \, dx \lesssim \delta \int \bigl( x^2 v_{xt}^2 + v_t^2 \bigr) \,dx + C_\delta P(\Lambda_0) \mathcal E_{1}. 
	\end{align*}	
	Therefore, after choosing $ \delta $ small enough, integration in the temporal variable of \eqref{ene:1st-001} and \eqref{ene:1st-002} yields.
	\begin{align*}
		& \norm{x\sqrt{\rho_0} v_t}{\stnorm{\infty}{2}}^2 + \norm{x v_{xt}}{\stnorm{2}{2}}^2 + \norm{v_t}{\stnorm{2}{2}}^2 \lesssim T P(\Lambda_0) (\mathcal E_1 + 1) + P(M_0) \mathcal E_1^0,\\
		& \norm{x \sqrt{\rho_0}\Theta_t}{\stnorm{\infty}{2}}^2 + \norm{x \Theta_{xt}}{\stnorm{2}{2}}^2 \lesssim \norm{xv_{xt}}{\stnorm{2}{2}}^2 + \norm{v_t}{\stnorm{2}{2}}^2 + TP(\Lambda_0)(\mathcal E_1 + 1)\\
		& ~~~~~~ + P(M_0) \mathcal E_1^0. 
	\end{align*}
	Thus we have shown \eqref{ene:1st-010}. \eqref{ene:1st-020} is the consequence of \eqref{ene:1st-010} and the fundamental theory of calculus.
\end{pf}

The next step is to derive the estimates around the symmetric center $ x = 0 $. In fact, we have the following lemma. 

\begin{lm}\label{lm:energyest-strong-002}
	For a solution to \eqref{eq:Lg000}(or equivalently \eqref{eq:LagrangianCoordinates}) and \eqref{eq:Lg010}, it holds
	\begin{align}
	& \norm{\sqrt{\chi\rho_0}v_t}{\stnorm{\infty}{2}}^2 + \norm{\sqrt{\chi}v_{xt}}{\stnorm{2}{2}}^2 + \norm{\sqrt{\chi}\dfrac{v_t}{x}}{\stnorm{2}{2}}^2 {\nonumber}\\
	& ~~~~ \lesssim TP(\Lambda_0) (\mathcal E_1 + 1 ) + P(M_0) \mathcal E_1^0,	{\label{ene:1th-101}}\\
	&  \norm{\sqrt{\chi}v_x}{\stnorm{\infty}{2}}^2 + \norm{\sqrt{\chi}\dfrac{v}{x}}{\stnorm{\infty}{2}}^2 + \norm{\sqrt{\chi\rho_0}v_t}{\stnorm{2}{2}}^2 {\nonumber}\\
	& ~~~~ \lesssim TP(\Lambda_0) (\mathcal E_1 + 1 ) + P(M_0) \mathcal E_1^0. {\label{ene:0th-201}}
	\end{align}
\end{lm}
\begin{pf}
	Multiply \subeqref{eq:Lg010}{1} with $ \chi \dfrac{v}{r^2} $ and integrate the resulting equation over $ x \in ( 0, R_0 ) $. After integration by parts, it holds
	\begin{equation}
		\begin{aligned}
			& \int \chi \dfrac{x^2}{r^2} \rho_0 \dt v_t v_t \,dx + \int \mathfrak B_t (\chi v_{t})_x - 4\mu \chi \bigl(\dfrac{v_t}{r}\bigr)_x v_t \,dx \\
			& ~~~~ = K \int x^2 \rho_0  I_1^1 \bigl(\chi \dfrac{v_t}{r^2} \bigr)_x \,dx + 2 K \int\chi x\rho_0 I_2^1 \dfrac{v_t}{r^2} \,dx + \int \chi \dfrac{x^2}{r^2} v_t I_3^1 \,dx\\
			& ~~~~~~~ + \int \chi \dfrac{x^2}{r^2} v_t  I_4^1 \,dx + \int \chi  \dfrac{x^2}{r^2} \rho_0 v_t I_5^1 \,dx : = \sum_{i=1}^{5} J_i^3 . 
		\end{aligned}
	\end{equation}
	To evaluate the left hand side, notice
	\begin{equation*}
		\begin{aligned}
			& \int \chi \dfrac{x^2}{r^2} \rho \dt v_t v_t \,dx = \dfrac{d}{dt}\biggl\lbrace \dfrac{1}{2} \int \chi \dfrac{x^2}{r^2} \rho_0 v_t^2 \,dx \biggr\rbrace - J_{6}^3, \\
			& \int \mathfrak B_t (\chi v_t)_x - 4\mu \chi \bigl(\dfrac{v_t}{r}\bigr)_x v_t \,dx = (2\mu+\lambda) \int\chi r_x \biggl( \bigl(\dfrac{v_{xt}}{r_x}\bigr)^2 + \bigl( \dfrac{v_t}{r} \bigr)^2 \biggr)   \,dx - J_{7}^3,
		\end{aligned}
	\end{equation*}
	with
	\begin{align*}
		& J_{6}^3 =  - \int \chi \dfrac{x^2v}{r^3} \rho_0 v_t^2\,dx, \\
		& J_{7}^3 =  - (2\mu+\lambda) \int \chi' \biggl( \dfrac{v_t^2}{r} + \dfrac{v_tv_{xt}}{r_x} \biggr) \,dx \\
		& ~~~~~~~ + \int \biggl( (2\mu+\lambda)\bigl(\dfrac{v_x}{r_x}\bigr)^2 +  2 \lambda \bigl(\dfrac{v}{r}\bigr)^2 \biggr) \bigl( \chi v_{xt} + \chi' v_t \bigr) \,dx .
	\end{align*}
	Now we estimate $ J_i^3 $'s. To begin with, 
	\begin{align*}
		& J_1^3 \lesssim P(M_0)\int x^2 \rho_0 \bigl(\abs{\Theta_t}{} + \abs{\Theta v_x}{}\bigr) \bigl(\abs{\chi' \dfrac{v_t}{x^2}}{} + \abs{\chi \dfrac{v_{xt}}{x^2}}{} + \abs{\chi \dfrac{v_t}{x^3}}{} \bigr) \,dx \lesssim P(M_0) \mathcal E_1 \\
		& ~~~~ + \delta \int \chi \bigl(v_{xt}^2 + \bigl( \dfrac{v_t}{x}\bigr)^2 \bigr) \,dx  + C_\delta \biggl( P(M_0) \int \rho_0 \Theta_t^2 \,dx + P(\Lambda_0) \int \Theta^2\,dx  \biggr)\\
		& \lesssim \delta \int \chi \bigl(v_{xt}^2 + \bigl( \dfrac{v_t}{x}\bigr)^2 \bigr) \,dx + (1+C_\delta) P(\Lambda_0) \mathcal E_1 + C_\delta P(M_0) \int x^2 \Theta_{xt}^2 \,dx,
	\end{align*}
	where we have used the following in the last inequality,
	\begin{equation*}
	 \int \Theta^2 \,dx \lesssim \int x^2 \Theta_x^2 \,dx \lesssim \mathcal E_1, ~ \int \rho_0 \Theta_t^2 \,dx \lesssim \int x^2 \Theta_{xt}^2 \,dx,
	\end{equation*}
	 by applying $ a), c) $. For $ J_7^3 $,
	 \begin{align*}
	 	& J_7^3 \lesssim P(M_0)\int \bigl( x^2v_{xt}^2 + v_t^2 \bigr) \,dx + \delta \int \chi \bigl( v_{xt}^2 + \bigl(\dfrac{v_t}{x}\bigr)^2 \bigr)\,dx \\
	 	& ~~~~ + (1 + C_\delta) P(\Lambda_0) \mathcal E_1.
	 \end{align*}
	Similarly, by applying Cauchy's inequality
	\begin{align*}
	& J_1^3 + J_2^3 + J_4^3 + J_5^3 + J_6^3 + J_7^3 \lesssim \delta \int \chi \bigl(v_{xt}^2 + \bigl( \dfrac{v_t}{x}\bigr)^2 \bigr) \,dx + (1+C_\delta) P(\Lambda_0) \mathcal E_1\\
	& ~~~~ + (1+ C_\delta) P(M_0) \biggl( \int x^2 \Theta_{xt}^2 \,dx + \int \bigl(x^2 v_{xt}^2 + v_t^2\bigr) \,dx \biggr).
	\end{align*}
	For $ J_3^3 $, after integration by parts, we employ the Cauchy's inequality
	\begin{align*}
		& J_3^3 = \int \chi \dfrac{x^2}{r^2}v_t \biggl( 2 \dfrac{rv}{x^2} \mathfrak B_x \biggr) \,dx = - 2 \int \mathfrak B \biggl( \chi \dfrac{vv_t}{r} \biggr)_x \,dx 
		 \lesssim \delta \int \chi \bigl( v_{xt}^2 + \bigl(\dfrac{v_t}{x}\bigr)^2\bigr) \,dx \\ 
		 & ~~~~ + C_\delta P(\Lambda_0) \mathcal E_1.
	\end{align*}
	Summing up these estimates with a small enough $ \delta $ yields,
	\begin{align*}
		& \dfrac{d}{dt} \biggl\lbrace \dfrac{1}{2}\int \chi \dfrac{x^2}{r^2}\rho_0 v_t^2 \,dx \biggr\rbrace + \dfrac{(2\mu+\lambda)}{2} \int\chi r_x \biggl( \bigl(\dfrac{v_{xt}}{r_x}\bigr)^2 + \bigl( \dfrac{v_t}{r} \bigr)^2 \biggr)   \,dx \\
		 & ~~~~ \lesssim P(\Lambda_0)\mathcal E_1 +  P(M_0) \biggl( \int x^2 \Theta_{xt}^2 \,dx + \int \bigl(x^2 v_{xt}^2 + v_t^2\bigr) \,dx \biggr).
	\end{align*}
	Therefore, integrating in the temporal variable implies,
	\begin{align*}
		& \norm{\sqrt{\chi\rho_0}v_t}{\stnorm{\infty}{2}}^2 + \norm{\sqrt{\chi}v_{xt}}{\stnorm{2}{2}}^2 + \norm{\sqrt{\chi}\dfrac{v_t}{x}}{\stnorm{2}{2}}^2 \lesssim TP(\Lambda_0) \mathcal E_1 \\
		& ~~~~ + P(M_0) \bigl(\norm{x\Theta_{xt}}{\stnorm{2}{2}}^2 + \norm{xv_{xt}}{\stnorm{2}{2}}^2 + \norm{v_t}{\stnorm{2}{2}}^2 \bigr) + P(M_0) \mathcal E_1^0.
	\end{align*}
	Then together with \eqref{ene:1st-010}, this will establish \eqref{ene:1th-101}. \eqref{ene:0th-201} is the consequence of \eqref{ene:1th-101} and the fundamental theory of calculus.
\end{pf}

Now we can show the main a prior estimate for the solution to \eqref{eq:LagrangianCoordinates} with $ \mathcal E_1^0 < \infty $. 
\begin{lm}\label{lm:energyest-strong-003}
For a solution to \eqref{eq:LagrangianCoordinates}, there is a time 
\begin{equation}\label{aprioriest:strong-time}
	T_* = T_*(\mathcal E_1^0) \geq \dfrac{1}{P_*(\mathcal E_1^0 + 1)},
\end{equation}
such that if $ T \leq \min\bigl\lbrace T_*, 1 \bigr\rbrace $, 
\begin{equation}\label{aprioriest:strong-energy}
	\mathcal E_1 \leq 1 + C_* \mathcal E_1^0,
\end{equation}
and consequently, 
\begin{equation}\label{aprioriest:strong-pointwise}
	M_0 \leq 2 ~~ \text{and} ~~ \Lambda_0 \leq C_*( 1 + \mathcal E_1^0 ),
\end{equation}
for some positive polynomial $ P_* = P_*(\cdot) $ and some positive constant $ 1 < C_* < \infty $. 
\end{lm}
\begin{pf}
	As a consequence of Lemma \ref{lm:energyest-strong-001}, Lemma \ref{lm:energyest-strong-002} and the fundamental theory of calculus, we have
	\begin{equation*}
		\mathcal E_1 \lesssim TP(\Lambda_0) (\mathcal E_1 + 1 ) + P(M_0) \mathcal E_1^0.
	\end{equation*}
	For $ T \leq T_1 $ (defined in Corollary \ref{cor:pointwise}), this will imply, after employing \eqref{cor:upperbound01},
	\begin{equation*}
		\mathcal E_1 \leq T P_2(\mathcal E_1 + 1) + C_2 \mathcal E_1^0 ,
	\end{equation*}
	for some polynomial $ P_2 = P_2(\cdot) $ and some positive constant $ 1 < C_2 < \infty $. Let
	\begin{equation*}
		T_2 : = \dfrac{1}{P_2(\mathcal E_1 + 1)}.
	\end{equation*}  
	Then for $ T \leq \min \bigl\lbrace T_1, T_2, 1  \bigr\rbrace $, we have
	\begin{equation*}
		\mathcal E_1 \leq 1 + C_2 \mathcal E_1^0, ~ \text{and} ~ T_2, T_1 \geq \dfrac{1}{P_3(\mathcal E_1^0 + 1)},
	\end{equation*}
	with some polynomial $ P_3 = P_3(\cdot) $. This finishes the proof.
\end{pf}

\subsection{Higher Regularity}\label{sec:classical}
In this section, we will perform higher order estimates of the solutions to \eqref{eq:LagrangianCoordinates}. In the following, it is assumed $ \mathcal E_1 < \infty, \Lambda_0 < \infty, M_0 \leq 2 $ and $ \mathcal E_2^0 < \infty $. We will show that for a short time, $ \mathcal E_2 < \infty $. In particular, we study the regularity of classical solutions in the space $ \mathfrak Y $. In order to do so, we first establish some elliptic estimates for the classical solutions, and then head to the energy estimates.

\subsubsection*{Elliptic Estimates for Classical Solutions}

In the following, we shall perform higher order elliptic estimates for the classical solutions. In particular, we will consider the following system, which consists of the temporal derivative version of \eqref{eq:altformofmmteq} and the rearrangement of \subeqref{eq:Lg010}{2},
\begin{equation}\label{ell:011}
	\begin{cases}
		& (2\mu+\lambda) \mathcal G_{xtt} = - K \dfrac{x^2\rho_0 \Theta}{r^2r_x} \mathcal G_{xt} - K \bigl( \dfrac{x^2\rho_0\Theta}{r^2r_x}\bigr)_t \mathcal G_{x}  +  K \dfrac{x^2}{r^2r_x} (\rho_0\Theta_t)_x   \\
		& ~~~~~~- K \dfrac{x^2 \dt(r^2r_x)}{(r^2r_x)^2} (\rho_0\Theta)_x + \dfrac{x^2}{r^2} \rho_0 v_{tt} - 2 \dfrac{x^2v}{r^3} \rho_0 v_t - 4 \dfrac{x^5v\rho_0}{r^5} \Phi, \\
		& \dfrac{1}{r\sqrt{r_x}} \bigl( \dfrac{r^2}{r_x}\Theta_{xt} \bigr)_x = c_\nu \dfrac{x^2\rho_0}{r\sqrt{r_x}} \Theta_{tt} + K \dfrac{x^2 \rho_0}{r\sqrt{r_x}} I_6^1 - \dfrac{x^2}{r\sqrt{r_x}} I_7^1 + \dfrac{( x^2 I_8^1 \Theta_x )_x}{r\sqrt{r_x}}. 
	\end{cases}
\end{equation}
In the meantime, the spatial derivatives of \eqref{eq:altformofmmteq} and \subeqref{eq:Lg000}{2} can be written as the following,
\begin{equation}\label{ell:012}
	\begin{cases}
		(2\mu+\lambda) \mathcal G_{xxt} = \bigl( \dfrac{x^2}{r^2} \rho_0 v_t\bigr)_x + \bigl( K \dfrac{x^2\rho_0}{r^2r_x}\Theta\bigr)_{xx}  + \bigl( \dfrac{x^5 \rho_0}{r^4} \Phi\bigr)_x,\\
		\bigl( r^2 \dfrac{\Theta_x}{r_x}\bigr)_{xx} = c_\nu ( x^2 \rho_0 \Theta_t )_x + K \bigl( x^2 \rho_0 \Theta \dfrac{\dt(r^2r_x)}{r^2r_x} \bigr)_x - \epsilon ( x^2 \rho_0 )_x - \mathbf S_x,
	\end{cases}
\end{equation}
where 
\begin{equation}
	\mathbf S = 2\mu r^2 r_x \biggl( \bigl(\dfrac{v_x}{r_x}\bigr)^2 + 2 \bigl( \dfrac{v}{r}\bigr)^2 \biggr) + \lambda r^2 r_x \bigl( \dfrac{v_x}{r_x} + 2 \dfrac{v}{r} \bigr)^2.
\end{equation}

The estimates on \eqref{ell:011} are presented in the next lemma.
\begin{lm}\label{lm:highregularity001} For a solution to \eqref{ell:011} and $ 0 < T < 1 $, it holds,
\begin{equation}\label{ell:200}
		\begin{gathered} 
			\norm{x\Theta_{xxt}}{\stnorm{\infty}{2}}^2 + \norm{\Theta_{xt}}{\stnorm{\infty}{2}}^2 + \norm{v_{xxt}}{\stnorm{\infty}{2}}^2 + \norm{\bigl(\dfrac{v_t}{x}\bigr)_x}{\stnorm{\infty}{2}}^2 \\
			~~~~~~~~ \lesssim P(\Lambda_0,\mathcal E_1)(1+\mathcal E_2).
		\end{gathered}
	\end{equation}
	In particular,
	\begin{equation}\label{ell:300}
		\begin{gathered}
			\norm{\Theta_t}{\supnorm}^2 +  \norm{x \Theta_{xt}}{\supnorm}^2 + \norm{v_t}{\supnorm}^2 + \norm{v_{xt}}{\supnorm}^2 \\
			~~~~~~ + \norm{\dfrac{v_t}{x}}{\supnorm}^2 \lesssim P(\Lambda_0,\mathcal E_1)(1+\mathcal E_2).
		\end{gathered}
	\end{equation}
	As a consequence, for $ R_0 / 2 < x \leq R_0 $, it shall have
	\begin{equation}\label{recoverofPVT}
		-\infty < \Theta_x \leq - c < 0, ~ {\dfrac{\Theta}{\sigma}}{} \geq c > 0,
	\end{equation}
	for some constant $ c > 0 $ and $ 0 < T < T_3 $ where
	\begin{equation}\label{aprioriest:classical-time-01} T_3 := \dfrac{1}{P_4(\Lambda_0,\mathcal E_1, \mathcal E_2)}. \end{equation}
	for some positive polynomial $ P_4 = P_4(\cdot) $.  
\end{lm}

\begin{pf}
	By taking square of \subeqref{ell:011}{2}, it holds
	\begin{align*}
		& \int \abs{\dfrac{1}{r\sqrt{r_x}} \bigl( \dfrac{r^2}{r_x}\Theta_{xt}\bigr)_x }{2}\,dx \lesssim P(\Lambda_0) \int x^2 \rho_0 \Theta_{tt}^2\,dx + P(\Lambda_0) \int x^2 \rho_0^2 \abs{I_6^1}{2}\,dx \\
		& ~~~~~~ + P(\Lambda_0) \int x^2 \abs{I_7^1}{2}\,dx + \int \abs{\dfrac{( x^2 I_8^1 \Theta_x )_x}{r\sqrt{r_x}}}{2}\,dx \\
		& ~~~~\lesssim P(\Lambda_0) \mathcal E_{2} + P(\Lambda_0) \int \abs{\dfrac{( x^2 I_8^1 \Theta_x )_x}{x}}{2}\,dx.
	\end{align*}
	In the meantime, notice
	\begin{align*}
		& \abs{ (x^2 I_8^1 \Theta_x)_x }{} \lesssim P(\Lambda_0) \abs{x \Theta_x}{} + P(\Lambda_0) \abs{x^2 \bigl( \abs{\dfrac{v_x}{x}}{} + \abs{\dfrac{1}{x}}{} + \abs{v_{xx}}{} + \abs{r_{xx}}{} \bigr) \Theta_x}{} \\
		& ~~~~ + P(\Lambda_0) \abs{x^2 \Theta_{xx} }{}.
	\end{align*}
	We then have
	\begin{align*}
		& \int \abs{\dfrac{( x^2 I_8^1 \Theta_x )_x}{x}}{2}\,dx \lesssim P(\Lambda_0) \int \bigl(\abs{\Theta_x}{2} + x^2 \abs{\Theta_{xx}}{2}\bigr) \,dx \\
		&~~~~ + P(\Lambda_0) \norm{x\Theta_x}{\supnorm}^2 \int \abs{v_{xx}}{2} + \abs{r_{xx}}{2} \,dx \lesssim P(\Lambda_0,\mathcal E_1),
	\end{align*}
	by noticing \eqref{uniformese02oct02}, \eqref{ell:100} in the last inequality. Thus we have shown
	\begin{equation}\label{ell:013}
		\int \abs{\dfrac{1}{r\sqrt{r_x}} \bigl( \dfrac{r^2}{r_x}\Theta_{xt}\bigr)_x }{2}\,dx \lesssim P(\Lambda_0) \mathcal{E}_2 + P(\Lambda_0,\mathcal E_1).
	\end{equation}
	Similarly, the square of \subeqref{ell:011}{1} then yields,
	\begin{align*}
			& \int \abs{\mathcal G_{xtt}}{2} \,dx \lesssim P(\Lambda_0) \bigl( \int \abs{\mathcal G_{xt}}{2} + \abs{\mathcal G_x}{2} \bigr) \,dx + P(\Lambda_0) \norm{\Theta_t}{\supnorm}^2 \int \abs{\mathcal G_x}{2} \,dx \\
			& ~~~~ + P(\Lambda_0) \bigl( \int \abs{\Theta_{xt}}{2}\,dx + \int \sigma^{2\alpha-2} \abs{\Theta_t}{2}\,dx + \int \Theta_x^2 \,dx\\
			& ~~~~ + \int \sigma^{2\alpha-2} \Theta^2\,dx + \int \rho_0 (v_{tt}^2 + v_t^2 + v^2 )\,dx \bigr) \\
			&   \lesssim P(\Lambda_0,\mathcal E_1) + P(\Lambda_0) \mathcal E_2 + P(\Lambda_0,\mathcal E_1) \bigl( \norm{\Theta_t}{\supnorm}^2 \\
			& ~~~~ + \int \abs{\Theta_{xt}}{2}\,dx + \int \sigma^{2\alpha-2} \abs{\Theta_t}{2}\,dx + \int \sigma^{2\alpha-2} \abs{\Theta}{2}\,dx  \bigr),
	\end{align*}
	where in the last inequality it has been applied \eqref{ell:050}, \eqref{ell:060}, \eqref{ell:100} and the fact $$ \int \rho_0 (v_{tt}^2 + v_t^2 + v^2 )\,dx \lesssim \int \chi \rho_0 (v_{tt}^2 + v_t^2 + v^2 )\,dx + \int x^2 \rho_0 (v_{tt}^2 + v_t^2 + v^2 )\,dx \lesssim \mathcal E_2.  $$
	From the embedding theory, we have $$ \norm{\Theta_t}{\supnorm}^2 \lesssim \int \Theta_t^2 \,dx + \int\Theta_{xt}^2 \,dx \lesssim \int \Theta_{xt}^2 \,dx. $$
	To evaluate the last two integrals on the right, for $ \alpha \neq 1/2 $, the inequalities $ a), b), c) $ in Lemma \ref{lm:hardy} then imply
	\begin{equation*}
		\int \sigma^{2\alpha-2}\abs{\Theta_t}{2}\,dx \lesssim \int x^2 \sigma^{2\alpha}\Theta_{xt}^2 \,dx \lesssim \int \Theta_{xt}^2\,dx;
	\end{equation*}
	for $ \alpha = 1/2 $, by applying Cauchy's inequality, Poincar\'e inequality and the inequality $ b) $ in Lemma \ref{lm:hardy},
	\begin{align*}
		& \int \sigma^{2\alpha-2}\abs{\Theta_t}{2}\,dx = \int \sigma^{-1} \abs{\Theta_t}{2} \,dx \lesssim \int \abs{\Theta_t}{2}\,dx + \int \sigma^{-2} \abs{\Theta_t}{2} \,dx \\
		& ~~~~ \lesssim \int\Theta_{xt}^2 \,dx.
	\end{align*}
	Similarly for $\int \sigma^{2\alpha-2} \Theta^2\,dx $. Thus it holds
	\begin{equation}\label{ell:014}
		\int \abs{\mathcal G_{xtt}}{2}\,dx \lesssim P(\Lambda_0, \mathcal E_1) ( 1 + \int \Theta_{xt}^2 \,dx + \mathcal E_2 ).
	\end{equation}
	
	From \eqref{ell:013} and \eqref{ell:014} we shall derive the desired estimates.
	With similar calculations as before, the left of \eqref{ell:013} satisfies the following,
	\begin{align*}
		& \int \abs{\dfrac{1}{r\sqrt{r_x}} \bigl( \dfrac{r^2}{r_x}\Theta_{xt}\bigr)_x }{2}\,dx \gtrsim \int \chi \abs{\dfrac{1}{r\sqrt{r_x}} \bigl( \dfrac{r^2}{r_x}\Theta_{xt}\bigr)_x }{2}\,dx \\
		& = \int \chi \dfrac{r^2}{r_x} \biggl( \bigl(\dfrac{\Theta_{xt}}{r_x}\bigr)_x \biggr)^2\,dx + 4 \int \chi r_x \bigl( \dfrac{\Theta_{xt}}{r_x}\bigr)^2 \,dx + 4 \int \chi r \bigl( \dfrac{\Theta_{xt}}{r_x}\bigr) \bigl( \dfrac{\Theta_{xt}}{r_x}\bigr)_x \,dx \\
		& = \int \chi \dfrac{r^2}{r_x} \biggl( \bigl(\dfrac{\Theta_{xt}}{r_x}\bigr)_x\biggr)^2 \,dx + 2 \int \chi r_x \bigl(\dfrac{\Theta_{xt}}{r_x}\bigr)^2 \,dx - 2 \int \chi' r \bigl(\dfrac{\Theta_{xt}}{r_x}\bigr)^2 \,dx \\
		& \geq \int \chi \dfrac{r^2}{r_x} \biggl( \bigl(\dfrac{\Theta_{xt}}{r_x}\bigr)_x\biggr)^2 \,dx + 2 \int \chi r_x \bigl(\dfrac{\Theta_{xt}}{r_x}\bigr)^2 \,dx.
	\end{align*}
	Therefore, \eqref{ell:013} then yields
	\begin{equation}\label{ell:015}
		\int \chi x^2 \biggl( \bigl(\dfrac{\Theta_{xt}}{r_x}\bigr)_x\biggr)^2 \,dx + \int \chi \Theta_{xt}^2 \,dx \lesssim P(\Lambda_0) \mathcal{E}_2 + P(\Lambda_0,\mathcal E_1),
	\end{equation}
	and
	\begin{equation}\label{ell:016}
		\int \Theta_{xt}^2 \,dx \lesssim \int \chi \Theta_{xt}^2 \,dx + \int x^2 \Theta_{xt}^2\,dx \lesssim P(\Lambda_0) \mathcal{E}_2 + P(\Lambda_0,\mathcal E_1).
	\end{equation}
	From \eqref{ell:013} and \eqref{ell:016},
	\begin{align}
		& \int x^2 \Theta_{xxt}^2 \,dx \lesssim P(\Lambda_0) \biggl( \int\Theta_{xt}^2\,dx + \norm{x\Theta_{xt}}{\supnorm}^2 \int r_{xx}^2 \,dx {\nonumber} \\
		& ~~~~ + \int \abs{\dfrac{1}{r\sqrt{r_x}} \bigl( \dfrac{r^2}{r_x}\Theta_{xt}\bigr)_x }{2}\,dx \biggr) \lesssim P(\Lambda_0) \mathcal{E}_2 + ( 1 + \norm{x\Theta_{xt}}{\supnorm}^2 ) P(\Lambda_0,\mathcal E_1). \label{ell:0193}
	\end{align}
	By applying the following inequality,
	\begin{equation}\label{ell:019}
		f^2(x) \lesssim \int f^2\,dx + \int 2 \abs{f f_x}{} \,dx \lesssim \delta \int f_x^2\,dx + (1+C_\delta)\int f^2 \,dx,
	\end{equation}
	we shall have on the right of the above inequality,
	\begin{equation}\label{ell:0191}
		\norm{x\Theta_{xt}}{\supnorm}^2  \lesssim \delta \int ( x^2 \Theta_{xxt}^2 + \Theta_{xt}^2 ) \,dx + ( C_\delta + 1) \int x^2 \Theta_{xt}^2 \,dx.
	\end{equation}
	Therefore by choosing $ \delta $ small enough so that $ \delta P(\Lambda_0,\mathcal E_1) << 1 $ , it holds from \eqref{ell:016}, \eqref{ell:0193} and \eqref{ell:0191},
	\begin{equation}\label{ell:0192}
		\int x^2 \Theta_{xxt}^2 \,dx \lesssim P(\Lambda_0) \mathcal{E}_2 + P(\Lambda_0,\mathcal E_1).
	\end{equation}
	
	In the meantime,
	\begin{displaymath}
		\mathcal{G}_{xtt} = \dfrac{x}{rr_x} \biggl( 2 r_x \bigl(\dfrac{v_t}{x}\bigr)_x + v_{xxt}\bigl(\dfrac{r}{x}\bigr)\biggr) + l_1 + l_2 + l_3 + l_4 + l_5,
	\end{displaymath}
	with
	\begin{align*}
			& l_1 = \dfrac{x}{rr_x} \biggl\lbrack 2\biggl(2v_x\bigl(\dfrac{v}{x}\bigr)_x + v_{xx}\bigl(\dfrac{v}{x}\bigr)\biggr) + 2 v_{xt}\bigl(\dfrac{r}{x}\bigr)_x + r_{xx}\bigl(\dfrac{v_t}{x}\bigr) \biggr\rbrack,\\
			& l_2 = \dfrac{x}{rr_x} \bigl(\dfrac{v_x}{r_x} + \dfrac{v}{r}\bigr)\biggl\lbrack 2r_x \bigl(\dfrac{v}{x}\bigr)_x + v_{xx}\bigl(\dfrac{r}{x}\bigr) + 2 v_x \bigl(\dfrac{r}{x}\bigr)_x + r_{xx}\bigl(\dfrac{v}{x}\bigr) \biggr\rbrack,\\
			& l_3 = \dfrac{x}{rr_x} \bigl( \dfrac{v_x}{r_x} + \dfrac{v}{r}\bigr)^2 \biggl\lbrack 2 r_x \bigl(\dfrac{r}{x}\bigr)_x + r_{xx}\bigl(\dfrac{r}{x}\bigr)\biggr\rbrack,\\
			& l_4 = -\dfrac{x}{rr_x} \biggl\lbrack \dfrac{v_{xt}}{r_x} + \dfrac{v_t}{r} - \biggl( \bigl(\dfrac{v_x}{r_x}\bigr)^2 + \bigl(\dfrac{v}{r}\bigr)^2 \biggr) \biggr\rbrack \biggl\lbrack2r_x\bigl(\dfrac{r}{x}\bigr)_x +r_{xx}\bigl(\dfrac{r}{x}\bigr)\biggr\rbrack, \\
			& l_5 = - \dfrac{x}{rr_x}\bigl(\dfrac{v_x}{r_x}+\dfrac{v}{r}\bigr)\biggl\lbrack 2 v_x \bigl(\dfrac{r}{x}\bigr)_x + v_{xx}\bigl(\dfrac{r}{x}\bigr) +2 r_x \bigl(\dfrac{v}{x}\bigr)_x +r_{xx}\bigl(\dfrac{v}{x}\bigr)\biggr\rbrack,
	\end{align*}
	satisfying 
	\begin{align*}
		& |l_1|, |l_2|, |l_3|, |l_4|, |l_5| \lesssim P(\Lambda_0) \biggl( \bigl(\dfrac{v}{x}\bigr)_x+ v_{xx} + \bigl(\dfrac{r}{x}\bigr)_x + r_{xx} \\
		& ~~~~~~~~ + \bigl( v_{xt} + \dfrac{v_t}{x}\bigr)\biggl(\bigl(\dfrac{r}{x}\bigr)_x + r_{xx} \biggr)  \biggr).
	\end{align*}
	From \eqref{ell:014}, \eqref{ell:100}, \eqref{ell:016},
	\begin{equation}\label{ell:017}
		\begin{aligned}
			& \int \biggl( 2 r_x \bigl(\dfrac{v_t}{x}\bigr)_x + v_{xxt}\bigl(\dfrac{r}{x}\bigr)\biggr)^2 \,dx \lesssim P(\Lambda_0,\mathcal E_1) \mathcal{E}_2 \\
			& ~~~~~~ + \bigl( 1 + \norm{v_{xt}}{\supnorm}^2 + \norm{\dfrac{v_t}{x}}{\supnorm}^2 \bigr)  P(\Lambda_0,\mathcal E_1).
		\end{aligned}
	\end{equation}
	Again, the left of \eqref{ell:017} can be evaluated as,
	\begin{align*}
		& \int \biggl(  2 r_x \bigl(\dfrac{v_t}{x}\bigr)_x + v_{xxt}\bigl(\dfrac{r}{x}\bigr) \biggr)^2\,dx = \int \bigl(\dfrac{r}{x}\bigr)^2 v_{xxt}^2\,dx + 4 \int r_x^2 \bigl(\dfrac{v_t}{x}\bigr)_x^2\,dx \,dx\\
		& ~~~~ + 4 \int r_x \bigl(\dfrac{r}{x}\bigr) \bigl(\dfrac{v_t}{x}\bigr)_x v_{xxt}\,dx,
	\end{align*}
	with
	\begin{align*}
		& \int r_x \bigl(\dfrac{r}{x}\bigr) \bigl(\dfrac{v_t}{x}\bigr)_x v_{xxt}\,dx = \dfrac{1}{2}\int x r_x \bigl(\dfrac{r}{x}\bigr) \biggl\lbrack\bigl(\dfrac{v_t}{x}\bigr)_{x}^2\biggr\rbrack_{x} \,dx + 2 \int r_x \bigl(\dfrac{r}{x}\bigr) \bigl(\dfrac{v_t}{x}\bigr)_x^2\,dx\\
		& =  \dfrac{1}{2}x r_x \bigl(\dfrac{r}{x}\bigr) \bigl(\dfrac{v_t}{x}\bigr)_x^2 \bigr|_{x=R_0} + \dfrac{3}{2}\int r_x \bigl(\dfrac{r}{x}\bigr)\bigl(\dfrac{v_t}{x}\bigr)_x^2\,dx - \dfrac{1}{2} \int x r_{xx}\bigl(\dfrac{r}{x}\bigr) \bigl(\dfrac{v_t}{x}\bigr)_x^2\,dx\\
		& ~~~~ - \dfrac{1}{2} \int x r_x \bigl(\dfrac{r}{x}\bigr)_x\bigl(\dfrac{v_t}{x}\bigr)_x^2\,dx.
	\end{align*}
	The boundary term in the above can be rewritten as
	\begin{equation*}
		   \dfrac{1}{2}x r_x \bigl(\dfrac{r}{x}\bigr) \bigl(\dfrac{v_t}{x}\bigr)_x^2 \bigr|_{x=R_0} = \biggl. \dfrac{1}{2}x r_x \bigl( \dfrac{r}{x} \bigr) \bigl( \dfrac{v_{xt}}{x} - \dfrac{v_t}{x^2}\bigr)^2 \bigr|_{x=R_0}, 
	\end{equation*}
	which can be bounded by $ P(\Lambda_0) ( \norm{v_{xt}}{\supnorm}^2 + \norm{v_t/x}{\supnorm}^2 ) $.
	Summing up from \eqref{ell:017} then gives us, together with the inequality \eqref{coecivity}, and \eqref{ell:100},
	\begin{align}\label{ell:018}
			& \int \bigl(\dfrac{r}{x}\bigr)^2 v_{xxt}^2\,dx + 9 \int \bigl( \dfrac{r}{x}\bigr)^2 \biggl\lbrack\bigl(\dfrac{v_t}{x}\bigr)_x\biggr\rbrack^2 \,dx \lesssim P(\Lambda_0,\mathcal E_1) \mathcal{E}_2 {\nonumber} \\
			& ~~~~~~ + \bigl( 1 + \norm{v_{xt}}{\supnorm}^2 + \norm{\dfrac{v_t}{x}}{\supnorm}^2 \bigr)  P(\Lambda_0,\mathcal E_1) {\nonumber} \\
			& ~~~~~~ + C \int \biggl( r_{xx}^2 + \bigl(\dfrac{r}{x}\bigr)_x^2 \biggr) x^2 \biggl\lbrack \bigl(\dfrac{v_t}{x}\bigr)_x\biggr\rbrack^2\,dx {\nonumber}\\
			& ~~~~ \lesssim P(\Lambda_0,\mathcal E_1) \mathcal{E}_2 + \bigl( 1 + \norm{v_{xt}}{\supnorm}^2 + \norm{\dfrac{v_t}{x}}{\supnorm}^2 \bigr)  P(\Lambda_0,\mathcal E_1),
	\end{align}
	where the last inequality follows from
	\begin{equation*}
		x^2 \biggl\lbrack \bigl(\dfrac{v_t}{x}\bigr)_x\biggr\rbrack^2 = \bigl\lbrack v_{xt} - \dfrac{v_t}{x} \bigr\rbrack^2 \lesssim \norm{v_{xt}}{\supnorm}^2 + \norm{\dfrac{v_t}{x}}{\supnorm}^2.
	\end{equation*}
	After applying the inequality \eqref{ell:019} with $ \delta $ sufficiently small so that $ \delta P(\Lambda_0,\mathcal E_1 ) <<1 $, \eqref{ell:018} then yields,	
	\begin{equation}
		\int v_{xxt}^2 \, dx + \int \biggl\lbrack \bigl( \dfrac{v_t}{x}\bigr)_x\biggr\rbrack^2 \,dx \lesssim P(\Lambda_0,\mathcal E_1)(1+\mathcal E_2),
	\end{equation}
	and thus
	\begin{equation}
		 \norm{v_{xt}}{\supnorm}^2 + \norm{\dfrac{v_t}{x}}{\supnorm}^2 \lesssim P(\Lambda_0,\mathcal E_1)(1+\mathcal E_2).
	\end{equation}
	This finishes the proof of \eqref{ell:200}. 
	The rests in \eqref{ell:300} follow easily by employing the standard embedding theory. 
	To show \eqref{recoverofPVT}, from \eqref{ell:300}, we have $$ \norm{\Theta_{xt}}{\supnorm(R_0/2,R_0)}, \norm{\dfrac{\Theta_t}{\sigma}}{\supnorm} < C, $$ for some $ C > 0 $. Therefore, \eqref{recoverofPVT} follows from the following ODE,
	\begin{align*}
		& \dt \Theta_x = \Theta_{xt}, ~ \dt \dfrac{\Theta}{\sigma} = \dfrac{\Theta_t}{\sigma},
	\end{align*}
	and the initial constraints on $ \Theta $ (see \eqref{PV:lagrangian}).
	This finishes the proof. 
\end{pf}

As a corollary, we shall have
\begin{cor}\label{cor:e1est}
	For a solution to \eqref{eq:LagrangianCoordinates} with $ \Lambda_0, \mathcal E_1, \mathcal E_2 < \infty $, 
	\begin{equation}\label{dxThetaCen}
		\lim_{x\rightarrow 0^+}  \Theta_x(x,t) = 0.
	\end{equation}
\end{cor}
\begin{pf}
The calculation as in \eqref{uniformese0000} shows
	\begin{equation*}
		x^2 \Theta_x \lesssim P(\Lambda_0,\mathcal E_1,\mathcal E_2)(1+ \norm{\Theta_t}{\supnorm}) \int_0^x x^2 \,dx,
	\end{equation*}
	of which the right is bounded by $ x^3 $. Thus
	$$ \Theta_x \lesssim x P(\Lambda_0,\mathcal E_1,\mathcal E_2)(1+\norm{\Theta_t}{\stnorm{\infty}{2}}) $$ Therefore \eqref{ell:300} and taking $ x\rightarrow 0^+ $ in the above inequality will yield  \eqref{dxThetaCen}.
\end{pf}

Now we present the estimates on \eqref{ell:012}. 

\begin{lm}\label{lm:highregularity002}
	Under the assumption that $ 2\alpha - 2 > -1 $ (equivalently $ \epsilon K < 2/3 $) so that
	\begin{equation}
		\int \abs{(\rho_0)_x}{2} \,dx, \int \sigma^2 \abs{(\rho_0)_{xx}}{2} \,dx < \int \sigma^{2\alpha-2} \,dx < +\infty,
	\end{equation}
	we have, for $ 0 < T < 1 $ 
	\begin{equation}\label{ell:400}
		\begin{gathered}
			\norm{r_{xxx}}{\stnorm{\infty}{2}}^2 + \norm{\bigl(\dfrac{r}{x}\bigr)_{xx}}{\stnorm{\infty}{2}}^2 + \norm{v_{xxx}}{\stnorm{\infty}{2}}^2 + \norm{\bigl(\dfrac{v}{x}\bigr)_{xx}}{\stnorm{\infty}{2}}^2 \\
			+ \norm{\Theta_{xx}}{\stnorm{\infty}{2}}^2 + \norm{x\Theta_{xxx}}{\stnorm{\infty}{2}}^2 \lesssim P(\Lambda_0,\mathcal E_1)(1 + \mathcal E_2).
		\end{gathered}
	\end{equation}
	Consequently, 
	\begin{equation}\label{ell:500}
	\begin{gathered}
		\norm{r_{xx}}{\supnorm},  \norm{\bigl(\dfrac r x \bigr)_x}{\supnorm}, \norm{v_{xx}}{\supnorm}, \norm{\bigl(\dfrac v x \bigr)_x}{\supnorm}, \\
		~~~~~~ \norm{\Theta_x}{\supnorm}, \norm{x \Theta_{xx}}{\supnorm} \lesssim \tilde P(\Lambda_0,\mathcal E_1)(1 + \mathcal E_2).
	\end{gathered}
	\end{equation}
\end{lm}

\begin{pf}
	By taking square of \subeqref{ell:012}{1}, clumsy but direct calculation yields the following,
	\begin{align}\label{ell:0190}
			& \int \mathcal G_{xxt}^2\,dx \lesssim P(\Lambda_0) \biggl( \norm{v_t}{\supnorm}^2 \int \bigl\lbrack \bigl(\dfrac{r}{x}\bigr)_x\bigr\rbrack^2 \,dx + \int \abs{(\rho_0)_x}{2} \,dx + \int v_{xt}^2 \,dx\biggr) {\nonumber}\\
			& ~~ + P(\Lambda_0) \biggl( \int \sigma^2 (\rho_0)_{xx}^2 \,dx + \int (\rho_0)_x^2 \Theta_x^2 \,dx + \int \Theta_{xx}^2 \,dx {\nonumber}\\
			& ~~~~ + ( \norm{\bigl(\dfrac{r}{x}\bigr)_x}{\supnorm}^2 + \norm{r_{xx}}{\supnorm}^2 ) \int \bigl( \abs{(\rho_0)_x}{2} + \Theta_x^2 \bigr) \,dx \biggr) \\
			& ~~ + P(\Lambda_0) \int \bigl( \bigl( \dfrac{r}{x}\bigr)_{xx}^2 +  r_{xxx}^2\bigr) \,dx {\nonumber}\\
			& ~~  \lesssim P(\Lambda_0,\mathcal E_1)(1+\mathcal E_2) +  P(\Lambda_0,\mathcal E_1)\biggl\lbrace \int \bigl( \bigl( \dfrac{r}{x}\bigr)_{xx}^2 + r_{xxx}^2\bigr) \,dx + \int \Theta_{xx}^2 \,dx \biggr\rbrace, {\nonumber}
	\end{align}
	where in the last inequality, it has been employed \eqref{PV:lagrangianboundary}, \eqref{ell:100}, \eqref{ell:300} and the embedding inequality,
	\begin{align*}
		& \norm{\bigl(\dfrac{r}{x}\bigr)_x}{\supnorm}^2 \lesssim \int \bigl(\dfrac r x \bigr)_x^2 \,dx + \int \bigl(\dfrac r x \bigr)_{xx}^2 \,dx, \\
		& \norm{r_{xx}}{\supnorm}^2 \lesssim \int r_{xx}^2 \,dx + \int r_{xxx}^2 \,dx,
	\end{align*}
	and the following inequality by applying $ a) $ in Lemma \ref{lm:hardy} 
	\begin{equation*}
		\int \abs{(\rho_0)_x\Theta_x}{2}  \,dx \lesssim\int \sigma^{2\alpha-2} \Theta_x^2 \,dx \lesssim \int \Theta_x^2 \,dx + \int \Theta_{xx}^2 \,dx.	\end{equation*}
		
	Similarly, the square of $ \dfrac{\eqref{ell:012}_{2}}{r\sqrt{r_x}} $ implies
	\begin{equation}\label{ell:021}
		\int \abs{ \dfrac{1}{r\sqrt{r_x}} \bigl( r^2 \dfrac{\Theta_x}{r_x}\bigr)_{xx} }{2}\,dx \lesssim P(\Lambda_0,\mathcal E_1)(1+\mathcal E_2),
	\end{equation} where \eqref{ell:200} is applied.
	To evaluate the left of \eqref{ell:021}, notice first
	\begin{align*}
		& \bigl( r^2 \dfrac{\Theta_x}{r_x}\bigr)_{xx} = r^2 \bigl( \dfrac{\Theta_x}{r_x} \bigr)_{xx} + 4 r r_x \bigl( \dfrac{\Theta_x}{r_x} \bigr)_x + 2 r_x^2 \dfrac{\Theta_x}{r_x} + 2 r r_{xx} \dfrac{\Theta_x}{r_x} \\
		& ~~ = r^2 \bigl( \dfrac{\Theta_x}{r_x} \bigr)_{xx} + 3 r r_x \bigl( \dfrac{\Theta_x}{r_x} \bigr)_x + 2 r r_{xx} \dfrac{\Theta_x}{r_x}  + \dfrac{r_x}{r} \bigl( r^2 \dfrac{\Theta_x}{r_x} \bigr)_x,
	\end{align*}
	and from \subeqref{eq:Lg000}{2} and \eqref{ell:300},
	\begin{align*}
		\abs{\dfrac{1}{r^2} \bigl( r^2 \dfrac{\Theta_x}{r_x} \bigr)_x}{} \lesssim P(\Lambda_0,\mathcal E_1)(1+\mathcal E_2).
	\end{align*}
	Thus from \eqref{ell:021} and \eqref{ell:100}, it holds,
	\begin{equation}\label{ell:0211}
		\begin{gathered}
			\int \abs{ \dfrac{r}{\sqrt{r_x}}  \bigl( \dfrac{\Theta_x}{r_x} \bigr)_{xx} + 3 \sqrt{r_x} \bigl( \dfrac{\Theta_x}{r_x} \bigr)_x }{2}\,dx \\
			~~~~~~~~~~~~~ \lesssim P(\Lambda_0,\mathcal E_1)(1 +  \norm{r_{xx}}{\supnorm}^2 + \mathcal E_2).
		\end{gathered}
	\end{equation}
	In the meantime, we have
	\begin{align*}
		& \int \chi \abs{ \dfrac{r}{\sqrt{r_x}}  \bigl( \dfrac{\Theta_x}{r_x} \bigr)_{xx} + 3 \sqrt{r_x} \bigl( \dfrac{\Theta_x}{r_x} \bigr)_x }{2}\,dx = \int \chi \dfrac{r^2}{r_x} \bigl\lbrack \bigl(\dfrac{\Theta_x}{r_x}\bigr)_{xx}\bigr\rbrack^2 \,dx \\
		& ~~~ + 9\int \chi r_x \bigl\lbrack \bigl( \dfrac{\Theta_x}{r_x}\bigr)_x\bigr\rbrack^2\,dx + 3 \int \chi r \bigl\lbrack \bigl(\bigl( \dfrac{\Theta_x}{r_x}\bigr)_x \bigr)^2\bigr\rbrack_x \,dx \\
		&  = \int \chi \dfrac{r^2}{r_x} \bigl\lbrack \bigl(\dfrac{\Theta_x}{r_x}\bigr)_{xx}\bigr\rbrack^2 \,dx + 6 \int \chi r_x \bigl\lbrack \bigl( \dfrac{\Theta_x}{r_x}\bigr)_x\bigr\rbrack^2\,dx - 3\int \chi' r \bigl\lbrack \bigl( \dfrac{\Theta_x}{r_x}\bigr)_x\bigr\rbrack^2 \,dx \\
		& \gtrsim \int \chi \dfrac{r^2}{r_x} \bigl\lbrack \bigl(\dfrac{\Theta_x}{r_x}\bigr)_{xx}\bigr\rbrack^2 \,dx  + 6 \int \chi r_x \bigl\lbrack \bigl( \dfrac{\Theta_x}{r_x}\bigr)_x\bigr\rbrack^2\,dx.
	\end{align*}
	Thus we shall have from \eqref{ell:100} and \eqref{ell:0211},
	\begin{align*}
		& \int \chi \Theta_{xx}^2 \,dx \lesssim P(\Lambda_0) \biggl( \int \chi r_x \bigl\lbrack \bigl( \dfrac{\Theta_x}{r_x}\bigr)_x\bigr\rbrack^2 \, dx +  \norm{r_{xx}}{\supnorm}^2 \int \Theta_x^2 \,dx \biggr) \\
		& ~~~~~~ \lesssim P(\Lambda_0,\mathcal E_1)(1 +  \int {r_{xxx}}^2 \,dx + \mathcal E_2),\\
		& \int \chi x^2 \Theta_{xxx}^2 \,dx \lesssim P(\Lambda_0) \biggl( \int \chi \dfrac{r^2}{r_x} \bigl\lbrack \bigl(\dfrac{\Theta_x}{r_x}\bigr)_{xx}\bigr\rbrack^2 \,dx +  \norm{r_{xx}}{\supnorm}^2 \biggl( \int x^2 \Theta_{xx}^2 \,dx \\
		& ~~~~ + \int r_{xx}^2 \,dx \biggr) +  \norm{x^2\Theta_{x}}{\supnorm}^2 \int r_{xxx}^2\,dx \biggr) \lesssim P(\Lambda_0,\mathcal E_1)(1 \\
		& ~~~~~~~~~ +  \int {r_{xxx}}^2 \,dx + \mathcal E_2), \\
		& \int_{R_0/2}^{R_0} x^2 \Theta_{xxx}^2\,dx \lesssim P(\Lambda_0) \int_{R_0/2}^{R_0}  \abs{ \dfrac{r}{\sqrt{r_x}}  \bigl( \dfrac{\Theta_x}{r_x} \bigr)_{xx} + 3 \sqrt{r_x} \bigl( \dfrac{\Theta_x}{r_x} \bigr)_x }{2}\,dx \\
		& ~~~~ + P(\Lambda_0) (1+ \norm{r_{xx}}{\supnorm}^2) \int ( x^2\Theta_{xx}^2 + \Theta_x^2 + r_{xx}^2)\,dx + P(\Lambda_0,\mathcal E_1) \int r_{xxx}^2 \,dx \\
		& ~~~~~~~ \lesssim P(\Lambda_0,\mathcal E_1)(1 +  \int {r_{xxx}}^2 \,dx + \mathcal E_2),
	\end{align*}
	where we have applied the embedding inequalities
	\begin{equation}\label{ell:0212}
	\begin{gathered} 
		\norm{x\Theta_x}{\supnorm}^2 \lesssim \norm{\Theta_x}{\stnorm{\infty}{2}}^2 + \norm{x\Theta_{xx}}{\stnorm{\infty}{2}}^2,\\
		\norm{r_{xx}}{\supnorm}^2 \lesssim \norm{r_{xx}}{\stnorm{\infty}{2}}^2 + \norm{r_{xxx}}{\stnorm{\infty}{2}}^2. \end{gathered}\end{equation}
	Therefore
	\begin{align}
		& \int \Theta_{xx}^2 \,dx \lesssim \int \chi \Theta_{xx}^2 \,dx + \int x^2 \Theta_{xx}^2 \,dx {\nonumber}\\
		& ~~~~~~~ \lesssim P(\Lambda_0,\mathcal E_1)(1 +  \int {r_{xxx}}^2 \,dx + \mathcal E_2), {\label{ell:022}}\\
		& \int x^2 \Theta_{xxx}^2 \,dx \lesssim \int \chi x^2 \Theta_{xxx}^2 \,dx + \int_{R_0/2}^{R_0} x^2 \Theta_{xxx}^2\,dx {\nonumber}\\
		& ~~~~~~~ \lesssim P(\Lambda_0,\mathcal E_1)(1 +  \int {r_{xxx}}^2 \,dx + \mathcal E_2). {\label{ell:023}}
	\end{align}
	
	On the other hand, from \eqref{ell:0190} and \eqref{ell:022}, 
	\begin{equation}\label{ell:024}
		\int \mathcal G_{xxt}^2\,dx  \lesssim P(\Lambda_0,\mathcal E_1)(1+\mathcal E_2) + P(\Lambda_0,\mathcal E_1) \int  \bigl( \bigl(  \dfrac{r}{x}\bigr)_{xx}^2 + r_{xxx}^2\bigr) \,dx.
	\end{equation}
	In the following, notice
	\begin{equation}
		\mathcal{G}_{xxt} = \dfrac{x}{rr_x}\biggl( 2r_x \bigl(\dfrac{v}{x}\bigr)_{xx} + v_{xxx}\bigl(\dfrac{r}{x}\bigr) \biggr) + l_1 + l_2 + l_3 + l_4 ,
	\end{equation}
	where
	\begin{align}
		& l_1 = \dfrac{x}{rr_x}\biggl( 2v_x \bigl(\dfrac{r}{x}\bigr)_{xx} + r_{xxx}\bigl(\dfrac{v}{x}\bigr) \biggr) - \dfrac{x}{rr_x}\biggl(\dfrac{v}{r} + \dfrac{v_x}{r_x} \biggr)\biggl( 2r_x \bigl(\dfrac{r}{x}\bigr)_{xx} + r_{xxx}\bigl(\dfrac{r}{x}\bigr) \biggr) ,{\nonumber} \\
		& l_2 = \dfrac{3 x}{rr_x} \biggl( v_{xx}\bigl(\dfrac{r}{x}\bigr)_x + r_{xx}\bigl(\dfrac{v}{x}\bigr)_x\biggr) - \dfrac{3x}{rr_x}\biggl(\dfrac{v}{r}+\dfrac{v_x}{r_x}\biggr) r_{xx}\bigl(\dfrac{r}{x}\bigr)_x, {\nonumber}\\
		& l_3 = \bigl(\dfrac{x}{rr_x}\bigr)_x \biggl( 2 v_x \bigl(\dfrac{r}{x}\bigr)_x + 2 r_x \bigl(\dfrac{v}{x}\bigr)_x + v_{xx}\bigl(\dfrac{r}{x}\bigr) + r_{xx} \bigl(\dfrac{v}{x}\bigr)\biggr),\\
		& l_4 = - \biggl( \dfrac{x^2}{r^2r_x^2} \bigl(r_x\dfrac{v}{x} + \dfrac{r}{x} v_x\bigr) \biggr)_x\biggl( 2r_x\bigl(\dfrac{r}{x}\bigr)_x + r_{xx}\bigl(\dfrac{r}{x}\bigr)\biggr). {\nonumber}
	\end{align}
	It shall holds from \eqref{ell:024} and \eqref{ell:0212},
	\begin{align}\label{ell:025}
			& \int \biggl(2 r_x \bigl(\dfrac{v}{x}\bigr)_{xx} + v_{xxx} \bigl(\dfrac r x \bigr) \biggr)^2 \,dx \\
			& ~~~~~~ \lesssim P(\Lambda_0,\mathcal E_1)\bigl(1+\mathcal E_2 + \int  \bigl( \bigl(  \dfrac{r}{x}\bigr)_{xx}^2 + r_{xxx}^2\bigr) \,dx \bigr).
	\end{align}
	Integration by parts as before then yields,
	\begin{align*}
		& \int \bigl( \dfrac{r}{x}\bigr)^2 v_{xxx}^2 \,dx + \int \biggl( 4 r_x^2 + 10r_x \dfrac{r}{x} - 2 x r_x \bigl( \dfrac{r}{x}\bigr)_x - 2 x r_{xx}\bigl(\dfrac{r}{x}\bigr) \biggr) \bigl(\dfrac{v}{x}\bigr)_{xx}^2\,dx \\
		& ~~~~~ = \int \biggl(2 r_x \bigl(\dfrac{v}{x}\bigr)_{xx} + v_{xxx} \bigl(\dfrac r x \bigr) \biggr)^2 \,dx- \biggl( 2 x r_x \bigl(\dfrac{r}{x}\bigr)  \bigl(\dfrac{v}{x}\bigr)_{xx}^2 \biggr)\biggr|_{x=R_0} \\
		& \lesssim P(\Lambda_0) \norm{v_{xx}}{\supnorm}^2 + P(\Lambda_0,\mathcal E_1)\bigl(1+\mathcal E_2 + \int  \bigl( \bigl(  \dfrac{r}{x}\bigr)_{xx}^2 + r_{xxx}^2\bigr) \,dx \bigr).
	\end{align*}
	In the meantime,
	\begin{align*}
		& 4 r_x^2 + 10r_x \dfrac{r}{x} - 2 x r_x \bigl( \dfrac{r}{x}\bigr)_x - 2 x r_{xx}\bigl(\dfrac{r}{x}\bigr) = 14 \bigl(\dfrac{r}{x}\bigr)^2 +16 x \bigl(\dfrac{r}{x}\bigr) \bigl(\dfrac{r}{x}\bigr)_x \\
		& ~~ + 2 x^2\bigl(\dfrac r x \bigr)_x^2 - 2 x r_{xx}\bigl(\dfrac r x \bigr)  \geq 13 \bigl(\dfrac{r}{x}\bigr)^2 - P(\Lambda_0) \biggl( x^2  \bigl( \dfrac r x \bigr)_x^2 + x^2 r_{xx}^2 \biggr) .
	\end{align*}
	Thus from \eqref{ell:025} and the calculation above of the left, we shall have,
	\begin{align}\label{ell:026}
		& \int v_{xxx}^2 \,dx + \int \bigl(\dfrac{v}{x}\bigr)_{xx}^2 \,dx \lesssim P(\Lambda_0) \int x^2 \biggl( r_{xx}^2 + \bigl(\dfrac r x  \bigr)_x^2 \biggr) \bigl(\dfrac{v}{x}\bigr)_{xx}^2 \,dx {\nonumber}\\
		& ~~~ + P(\Lambda_0) \norm{v_{xx}}{\supnorm}^2 + P(\Lambda_0,\mathcal E_1)(1+\mathcal E_2 + \int  \bigl( \bigl(  \dfrac{r}{x}\bigr)_{xx}^2 + r_{xxx}^2\bigr) \,dx \bigr) {\nonumber}\\
		& \lesssim \biggl( \norm{r_{xx}}{\supnorm}^2 + \norm{\bigl(\dfrac r x \bigr)_x}{\supnorm}^2 \biggr) \int \biggl( v_{xx}^2 + \bigl( \dfrac v x \bigr)_x^2\biggr) \,dx\\
		& ~~~ + P(\Lambda_0) \norm{v_{xx}}{\supnorm}^2 + P(\Lambda_0,\mathcal E_1)(1+\mathcal E_2 + \int  \bigl( \bigl(  \dfrac{r}{x}\bigr)_{xx}^2 + r_{xxx}^2\bigr) \,dx \bigr) {\nonumber},
	\end{align}
	where it has been used
	\begin{align*}
		& x \bigl( \dfrac v x \bigr)_{xx} = v_{xx} - 2 \bigl( \dfrac v x \bigr)_x.
	\end{align*}
	Then by applying \eqref{ell:019} to $ \norm{v_{xx}}{\supnorm}^2 $ with a small enough $ \delta $, together with \eqref{ell:100}, \eqref{ell:026} then yields
	\begin{equation}\label{ell:027}
		\int v_{xxx}^2 \,dx + \int \bigl(\dfrac{v}{x}\bigr)_{xx}^2 \,dx  \lesssim P(\Lambda_0,\mathcal E_1)(1+\mathcal E_2 + \int  \bigl( \bigl(  \dfrac{r}{x}\bigr)_{xx}^2 + r_{xxx}^2\bigr) \,dx \bigr).
	\end{equation}
	Or by noticing
	\begin{equation*}
		v_{xxx} = \dfrac{d}{dt} r_{xxx}, ~~ \bigl(\dfrac{v}{x}\bigr)_{xx} = \dfrac{d}{dt} \bigl(\dfrac{r}{x}\bigr)_{xx},
	\end{equation*}
	we have
	\begin{equation}\label{ell:028}
		\dfrac{d}{dt} \int \bigl( r_{xxx}^2 + \bigl( \dfrac r x \bigr)_{xx}^2 \bigr) \,dx \lesssim P(\Lambda_0,\mathcal E_1)(1+\mathcal E_2 +  \int \bigl( r_{xxx}^2 + \bigl( \dfrac r x \bigr)_{xx}^2 \bigr) \,dx ),
	\end{equation}
	which implies, after applying Gr\"onwall's inequality together with \eqref{ell:027},  
	\begin{equation}
		\int \bigl( r_{xxx}^2 + \bigl(\dfrac{r}{x}\bigr)_{xx}^2 \bigr) \,dx + \int \bigl( v_{xxx}^2 + \bigl(\dfrac{v}{x}\bigr)_{xx}^2 \bigr) \,dx \lesssim P(\Lambda_0,\mathcal E_1)(1 + \mathcal E_2),
	\end{equation}
	for $ 0 < T < 1 $.
	Plugging back to \eqref{ell:022} and \eqref{ell:023} then gives the following,
	\begin{equation}
		\int \Theta_{xx}^2 \,dx + \int x^2 \Theta_{xxx}^2 \,dx \lesssim P(\Lambda_0,\mathcal E_1)(1 + \mathcal E_2).
	\end{equation}
	This finishes the proof.
\end{pf}

\subsubsection*{Energy Estimates for Classical Solutions}

To show that $ \mathcal E_2 <\infty $ for the solution with $ \mathcal E_2^0 < \infty $, we study the following temporal derivative version of \eqref{eq:Lg010},
\begin{equation}\label{eq:Lg020}
	\begin{cases}
		x^2 \rho_0 \dt v_{tt} - r^2  \mathfrak B_{xtt} - 4 \mu r^2  \bigl( \dfrac{v_{tt}}{r} \bigr)_x = -  \bigl( K x^2 \rho_0 I_1^2 \bigr)_x + 2 K x \rho_0 I_2^2 \\
		~~~~~~~~ + x^2 I_3^2 + x^2 I_4^2 + x^2 \rho_0 I_5^2,\\
		c_\nu x^2 \rho_0 \dt \Theta_{tt} -  \bigl( \dfrac{r^2}{r_x}  \Theta_{xtt} \bigr)_x = - K x^2 \rho_0 I_6^2 + x^2 I_7^2 - \bigl( x^2 I_8^2 \bigr)_x ,
	\end{cases}
\end{equation}
where
\begin{align*}
		& I_1^2 = \dt^2\bigl( \dfrac{\Theta}{r_x} \bigr), ~ I_2^2 = \dt^2 \bigl( \Theta \dfrac{x}{r}\bigr) ,~ I_3^2 = \dt^2 \bigl(\dfrac{r^2}{x^2} \mathfrak B_x \bigr) - \dfrac{r^2}{x^2} \mathfrak B_{xtt}, {\nonumber}\\
		& I_4^2 = \dt^2 \bigl( 4\mu \dfrac{r^2}{x^2}  \bigl( \dfrac{v}{r}\bigr)_x \bigr) - 4\mu\dfrac{r^2}{x^2} \bigl( \dfrac{v_{tt}}{r}\bigr)_x, ~ I_5^2 = - \dt^2\bigl(\dfrac{x^3}{r^2} \Phi \bigr), {\nonumber}\\
		& I_6^2 = \dt^2 \bigl( \Theta \dfrac{(r^2r_x)_t}{r^2r_x} \bigr), ~ \\
		& I_7^2 = \dt^2 \biggl\lbrack 2\mu \dfrac{r^2 r_x}{x^2} \biggl( \bigl(\dfrac{v_x}{r_x}\bigr)^2 + 2 \bigl( \dfrac{v}{r}\bigr)^2 \biggr) + \lambda \dfrac{r^2 r_x}{x^2} \bigl( \dfrac{v_x}{r_x} + 2\dfrac{v}{r} \bigr)^2 \biggr\rbrack, {\nonumber}\\
		& I_8^2 = \dt \bigl( I_8^1  \Theta_x \bigr) - \dt \bigl( \dfrac{r^2}{x^2r_x} \bigr) \Theta_{xt}.{\nonumber}
\end{align*}
\eqref{eq:Lg020} should be complemented with the following boundary conditions
\begin{equation*}
	\dt^2 v(0,t) = 0 , ~ \dt^{2}\Theta(R_0,t) = 0, ~ \mathfrak B_{tt} (R_0,t) = 0, ~ t \geq 0.
\end{equation*}

Similar as before, we shall present energy estimates on \eqref{eq:Lg020}. 
Multiply \subeqref{eq:Lg020}{1} with $ v_{tt} $ and \subeqref{eq:Lg020}{2} with $ \Theta_{tt} $ respectively and integrate the resulting equations over $ x \in (0, R_0) $. After integration by parts, it holds the following,
		\begin{align}
		& \dfrac{d}{dt} \biggl\lbrace \dfrac{1}{2} \int x^2 \rho_0 v_{tt}^2 \,dx \biggr\rbrace + \int \mathfrak B_{tt} ( r^2 v_{tt} )_x \,dx - \int 4\mu r^2 \bigl( \dfrac{v_{tt}}{r}\bigr)_x v_{tt}\,dx {\nonumber}\\
		& ~~~~ = K \int x^2 \rho_0 I_1^2 v_{xtt} \,dx + 2 K \int x \rho_0 I_2^2 v_{tt}\,dx + \int x^2 I_3^2 v_{tt}\,dx {\nonumber}\\
		& ~~~~~~~ + \int x^2 I_4^2 v_{tt}\,dx + \int x^2 \rho_0 I_5^2 v_{tt}\,dx : = \sum_{i=1}^{5} J_i^4,  {\label{ene:2nd-001}} \\
		& \dfrac{d}{dt} \biggl\lbrace \dfrac{c_\nu}{2}\int x^2 \rho_0 \Theta_{tt}^2 \,dx \biggr\rbrace + \int \dfrac{r^2}{r_x} \Theta_{xtt}^2 \,dx = - K \int x^2 \rho_0 I_6^2 \Theta_{tt}\,dx {\nonumber}\\
		& ~~~~~~~ + \int x^2 I_7^2 \Theta_{tt}\,dx + \int x^2 I_8^2 \Theta_{xtt} : = \sum_{i=6}^8 J_i^4. \label{ene:2nd-002}
		\end{align}
	By noticing
	\begin{align*}
		& \mathfrak B_{tt} ( r^2 v_{tt})_x - 4\mu r^2 \bigl(\dfrac{v_{tt}}{r}\bigr)_x v_{tt} = r^2 r_x \biggl\lbrace 2\mu \bigl( \dfrac{v_{xtt}^2}{r_x^2} + 2 \dfrac{v_{tt}^2}{r^2} \bigr) + \lambda \bigl( \dfrac{v_{xtt}}{r_x} + 2 \dfrac{v_{tt}}{r} \bigr)^2 \biggr\rbrace \\
		& ~~~~ - (r^2 v_{tt})_x\biggl\lbrack 3 \biggl( (2\mu+\lambda)\dfrac{v_x}{r_x}\dfrac{v_{xt}}{r_x} + 2\lambda \dfrac{v}{r}\dfrac{v_t}{r} \biggr) - \biggl( (2\mu+\lambda) \dfrac{v_x^3}{r_x^3} + 2\lambda \dfrac{v^3}{r^3} \biggr) \biggr\rbrack,
	\end{align*}
	\eqref{ene:2nd-001} and \eqref{ene:2nd-002} can be written as
	\begin{align*}
		& \dfrac{d}{dt} \biggl\lbrace \dfrac{1}{2} \int x^2 \rho_0 v_{tt}^2 \,dx \biggr\rbrace + \int r^2r_x \biggl(2\mu \bigl(\dfrac{v_{xtt}}{r_x}\bigr)^2 + 2 \bigl( \dfrac{v_{tt}}{r}\bigr)^2 \bigr) + \lambda\bigl( \dfrac{v_{xtt}}{r_x} + 2 \dfrac{v_{tt}}{r}\bigr)^2 \biggr) \,dx\\
		& ~~~~~~~~~~~~~ = \sum_{i=1}^5 J_i^4 + J_9^4, \\
		& \dfrac{d}{dt} \biggl\lbrace \dfrac{c_\nu}{2} \int x^2 \rho_0\Theta_{tt}^2 \,dx  \biggr\rbrace + \int \dfrac{r^2}{r_x} \Theta_{xtt}^2 \,dx = \sum_{i=6}^8 J_i^4,
	\end{align*}
	with
	\begin{equation}
		J_9^4 = \int (r^2 v_{tt})_x\biggl\lbrack 3 \bigl( (2\mu+\lambda)\dfrac{v_x}{r_x}\dfrac{v_{xt}}{r_x} + 2\lambda \dfrac{v}{r}\dfrac{v_t}{r} \bigr) - \bigl( (2\mu+\lambda) \dfrac{v_x^3}{r_x^3} + 2\lambda \dfrac{v^3}{r^3} \bigr) \biggr\rbrack \,dx.
	\end{equation}
	Similar as before, $ J_i^4 $'s can be evaluated as follows,
	\begin{align*}
		&  J_1^4 + J_2^4 + J_4^4 + J_5^4 + J_9^4 \lesssim \delta \int \bigl( x^2 v_{xtt}^2 + v_{tt}^2 \bigr) \,dx + C_\delta P(\Lambda_0) \mathcal E_2,  \\
		& J_6^4 + J_7^4 + J_8^4 \lesssim \delta \biggl( \int \bigl(x^2 v_{xtt}^2 + v_{tt}^2 \bigr) \,dx + \int x^2 \Theta_{xtt}^2 \,dx \biggr)  \\
		& ~~~~~~ + (1+C_\delta)P(\Lambda_0)\bigl( 1 + \norm{v_{xt}}{\supnorm}^2 + \norm{\dfrac{v_t}{x}}{\supnorm}^2 \bigr) \mathcal E_2.
	\end{align*}
	To evaluate $ J_3^4 $, after noticing 
	\begin{equation*}
		I_3^2 = 4 \dfrac{rv}{x^2} \mathfrak B_{xt} + 2 \bigl( \dfrac{v^2}{x^2} + \dfrac{rv_t}{x^2} \bigr) \mathfrak B_x,
	\end{equation*}
	integration by parts yields,
	\begin{align*}
		& J_3^4 
		= - \int \mathfrak B_{t} ( 4 r v v_{tt} )_x + 2 \mathfrak B ( v^2 v_{tt} + r v_t v_{tt} )_x \,dx \\
		& ~ \lesssim \delta \int\bigl( x^2 v_{xtt}^2 + v_{tt}^2 \bigr) \,dx + C_\delta P(\Lambda_0) \mathcal E_2. \\
	\end{align*}
	Thus integration in the temporal variable in \eqref{ene:2nd-001} and \eqref{ene:2nd-002} then yields, after choosing $ \delta <<1 $ small enough, 
	\begin{align}
	& \norm{x\sqrt{\rho_0}v_{tt}}{\stnorm{\infty}{2}}^2 + \norm{xv_{xtt}}{\stnorm{2}{2}}^2 + \norm{v_{tt}}{\stnorm{2}{2}}^2 {\nonumber}\\
	& ~~~~ \lesssim TP(\Lambda_0)\mathcal E_2 + P(M_0) \mathcal E_2^0,
	{\label{ene:2nd-010}}   \\
	& \norm{x\sqrt{\rho_0}\Theta_{tt}}{\stnorm{\infty}{2}}^2 + \norm{x\Theta_{xtt}}{\stnorm{2}{2}}^2 \lesssim \norm{xv_{xtt}}{\stnorm{2}{2}}^2 + \norm{v_{tt}}{\stnorm{2}{2}}^2 {\nonumber}\\
	& ~~~~ + T P(\Lambda_0)\bigl( 1 + \norm{v_{xt}}{\supnorm}^2 + \norm{\dfrac{v_t}{x}}{\supnorm}^2 \bigr) \mathcal E_2 + P(M_0)\mathcal E_2^0.
	{\label{ene:2nd-020}}	
	\end{align}
	
	On the other hand, multiply \subeqref{eq:Lg020}{1} with $ \chi \dfrac{v_{tt}}{r^2} $ and integrate the resulting equation over $ x \in ( 0, R_0 ) $. After integration by parts, it holds,
	\begin{align*}
		& \int \chi \dfrac{x^2}{r^2} \rho_0 v_{tt} \dt v_{tt} \,dx + \int \bigl( \mathfrak B_{tt} + 4\mu \dfrac{v_{tt}}{r} \bigr) ( \chi v_{tt} )_x \,dx = K \int x^2 \rho_0 I_1^2 \bigl( \chi \dfrac{v_{tt}}{r^2} \bigr)_x \,dx \\
		&  ~~ + 2K \int \chi x\rho_0 I_2^2 \dfrac{v_{tt}}{r^2} \,dx + \int \chi \dfrac{x^2}{r^2} v_{tt} (I_3^2 + I_4^2)  \,dx + \int \chi \dfrac{x^2}{r^2} \rho_0 v_{tt} I_5^2 \,dx : = \sum_{i=10}^{13} J_i^4.
	\end{align*}
	The left side of the above identity can be rewritten as, after integration by parts,
	\begin{align*}
		& \dfrac{d}{dt} \biggl\lbrace \dfrac{1}{2} \int \chi \dfrac{x^2}{r^2}\rho_0 v_{tt}^2 \biggr\rbrace + (2\mu+\lambda) \int \chi r_x \biggl( \bigl(\dfrac{v_{xtt}}{r_x}\bigr)^2 + \bigl(\dfrac{v_{tt}}{r}\bigr)^2\biggr) \,dx - J_{14}^4, 
	\end{align*}
	where
	\begin{align*}
		& J_{14}^4 = - \int \chi \dfrac{x^2 v}{r^3} \rho_0 v_{tt}^2 \,dx - (2\mu+\lambda) \int \chi' v_{tt} \bigl( \dfrac{v_{xtt}}{r_x} + \dfrac{v_{tt}}{r} \bigr) \,dx \\
		&~~~~  + \int \biggl\lbrace 3 \biggl( (2\mu+\lambda)\dfrac{v_{xt}v_x}{r_x^2} + 2\lambda\dfrac{v_tv}{r^2}\biggr) - 2 \biggl( (2\mu+\lambda)\dfrac{v_x^3}{r_x^3} + 2\lambda \dfrac{v^3}{r^3}\biggr) \biggr\rbrace (\chi v_{tt})_x \,dx.
	\end{align*}
	Thus it ends up with
	\begin{equation*}
		\dfrac{d}{dt}\biggl\lbrace \dfrac{1}{2} \int \chi \dfrac{x^2}{r^2}\rho_0v_{tt}^2 \,dx \biggr\rbrace + (2\mu+\lambda) \int \chi r_x \biggl( \bigl(\dfrac{v_{xtt}}{r_x}\bigr)^2 + \bigl(\dfrac{v_{tt}}{r}\bigr)^2\biggr) \,dx = \sum_{i=10}^{14}J_i^4.
	\end{equation*}
	By applying Cauchy's inequality and the Hardy's inequality, the following estimate can be achieved,
	\begin{align*}
		& J_{10}^4 + J_{11}^4 + J_{13}^4 + J_{14}^4 \lesssim \delta \int \chi \bigl( v_{xtt}^2 + \bigl(\dfrac{v_{tt}}{x}\bigr)^2 \bigr) \,dx \\
		& ~~~~~~ + (1+C_\delta )P(\Lambda_0) \int \bigl( x^2 \Theta_{xtt}^2 + x^2 v_{xtt}^2 + v_{tt}^2\bigr) \,dx  + (1 + C_\delta) P(\Lambda_0)\mathcal E_2.
	\end{align*}
	In the meantime, integration by parts yields together with Cauchy's inequality, 
	\begin{align*}
		& J_{12}^4 = - \int \biggl( \bigl( 4\chi \dfrac{v}{r}v_{tt} \bigr)_x \mathfrak B_t + \bigl( 2 \chi \bigl( \dfrac{v^2}{r^2} + \dfrac{v_t}{r} \bigr) v_{tt}\bigr)_x \mathfrak B \biggr) \,dx \\
		& ~~~~ - \int \biggl( \bigl( 8\mu \chi \bigl( \dfrac{v^2}{r^2} + \dfrac{v_t}{r} \bigr) v_{tt} \bigr)_x \dfrac{v}{r} + \bigl( 16 \mu \chi \dfrac{v}{r} v_{tt} \bigr)_x \bigl( \dfrac{v_t}{r}-\dfrac{v^2}{r^2} \bigr) \\
		& ~~~~  + ( 4\mu \chi v_{tt} )_x \bigl( -3\dfrac{v_t v}{r^2} + 2 \dfrac{v^3}{r^3}\bigr) \biggr) \,dx \lesssim \delta \int \chi\bigl(v_{xtt}^2 + \bigl(\dfrac{v_{tt}}{x}\bigr)^2 \bigr) \,dx \\
		& ~~~~ + (1+C_\delta )P(\Lambda_0) \int \bigl(x^2 v_{xtt}^2 + v_{tt}^2\bigr) \,dx  + (1 + C_\delta) P(\Lambda_0)\mathcal E_2.
	\end{align*}
	Therefore, after choosing $ \delta <<1 $ small enough, integration in the temporal variable in the above yields,
	\begin{equation}\label{ene:2nd-030}
	\begin{aligned}
		& \norm{\sqrt{\chi \rho_0}v_{tt}}{\stnorm{\infty}{2}}^2 + \norm{\sqrt{\chi}v_{xtt}}{\stnorm{2}{2}}^2 + \norm{\sqrt{\chi} \dfrac{v_{tt}}{x}}{\stnorm{2}{2}}^2 \\
		& ~~~~ \lesssim P(\Lambda_0)(\norm{x\Theta_{xtt}}{\stnorm{2}{2}}^2 + \norm{xv_{xtt}}{\stnorm{2}{2}}^2 + \norm{v_{tt}}{\stnorm{2}{2}}^2 + T\mathcal E_2 ) + P(M_0) \mathcal E_2^0.
	\end{aligned}
	\end{equation} 
	
From \eqref{ene:2nd-010}, \eqref{ene:2nd-020} and \eqref{ene:2nd-030}, we will derive the a prior estimates for the classical solution.
\begin{lm} \label{lm:energyest-classical}
For a solution to \eqref{eq:LagrangianCoordinates}, there is a time 
\begin{equation}\label{aprioriest:classical-time}
	T_{**} = T_{**}(\mathcal E_1^0, \mathcal E_2^0) \geq \dfrac{1}{P_{**}(\mathcal E_1^0 + 1, \mathcal E_2^0 + 1)},
\end{equation}
such that if $ T \leq \min \bigl\lbrace T_{**},1 \bigr\rbrace $,
\begin{equation}\label{aprioriest:classical-energy}
	\mathcal E_2 \leq 1 + C_* \mathcal E_1^0 + C_{**} \mathcal E_2^0,
\end{equation}
for some positive polynomial $ P_{**} = P_{**}(\cdot) $ and some positive constant $ 1 < C_{**} < \infty $. 
\end{lm}

\begin{pf}From \eqref{ene:2nd-010}, \eqref{ene:2nd-020} and \eqref{ene:2nd-030}, after applying the fundamental theory of calculus in the temporal variable, we have
\begin{align*}
	& \mathcal E_2 \lesssim \mathcal E_1 + TP(\Lambda_0) (1+\norm{v_{xt}}{\supnorm}^2 + \norm{\dfrac{v_t}{x}}{\supnorm}^2 ) \mathcal E_2 + P(M_0) \mathcal E_2^0.
\end{align*}
Then, together with \eqref{ell:300}, it implies
\begin{equation*}
	\mathcal E_2 \leq \mathcal E_1 + T P_5(\Lambda_0) (\mathcal E_2 + 1)^2 + P_6(M_0)\mathcal E_2^0,
\end{equation*}
for some positive polynomial $ P_5 = P_5(\cdot), P_6 = P_6(\cdot) $. Let
\begin{equation}\label{aprioriest:classical-time-02}
	T_4 : = \dfrac{1}{P_5(\Lambda_0) (\mathcal E_2 + 1)^2}.
\end{equation}
Then for $ T \leq \min\bigl\lbrace T_*, T_3 ,T_4 , 1 \bigr\rbrace $, it follows,
\begin{equation*}
	\mathcal E_2 \leq 1 + C_* \mathcal E_1^0 + C_{**} \mathcal E_2^0,
\end{equation*}
for some constant $ 0 <  C_{**} < \infty $. Together with \eqref{aprioriest:strong-pointwise}, we have 
\begin{equation*}
	T_3, T_4 \geq \dfrac{1}{P_7(\mathcal E_1^0 + 1,\mathcal E_2^0 + 1)}.
\end{equation*}
This finishes the proof.
\end{pf}

\section{Local Well-posedness of \eqref{eq:LagrangianCoordinates}}\label{sec:well-posedness}
In this section, we will use the a priori estimates obtained in Section \ref{sec:aprior} to derive the local well-posedness of the problem \eqref{eq:LagrangianCoordinates}. In particular, we introduce the following Hilbert space
\begin{equation}
	\begin{aligned}
		\mathfrak M : = & \bigl\lbrace (r,v,\Theta) \bigr| r, v \in L_t^2H_x^2, v_t \in L_t^2H_x^1, \Theta \in L_t^2 H_x^1 , x\Theta \in L_t^2 H_x^2, \\
		& x\Theta_t \in  L_t^2 H_x^1      \bigr\rbrace.
	\end{aligned}
\end{equation}
$ \mathfrak M $ is endowed with the norm
\begin{equation}
\begin{aligned}
	\norm{(r, v, \Theta)}{\mathfrak M}^2: = & \norm{r}{L_t^2H_x^2}^2 + \norm{v}{L_t^2H_x^2}^2 + \norm{v_t}{L_t^2H_x^1}^2 + \norm{\Theta}{L_t^2H_x^1}^2 \\
	 & ~~~~ + \norm{x\Theta}{L_t^2H_x^2}^2 + \norm{x\Theta_t}{L_t^2H_x^1}^2.
 \end{aligned}
\end{equation}
A closed, bounded, convex subset of $ \mathfrak M $ is defined as
\begin{equation}
\begin{aligned}
	\mathfrak S = & \mathfrak S_{\bar M} : = \bigl\lbrace (r,v,\Theta)\in \mathfrak M \bigl| \norm{(r,v,\Theta)}{\mathfrak M} < \tilde P (\bar M), \mathcal E_1(v,\Theta) \leq \bar M, \\
	& \Lambda_0(r,v,\Theta) \leq  {\bar M}, M_0(r) \leq 2  \bigr\rbrace,
\end{aligned}
\end{equation}
for some positive polynomial $ \tilde P = \tilde P(\cdot) $ to be determined later. We will write down the linearized problem of \eqref{eq:LagrangianCoordinates} and construct a map
\begin{equation}
	\mathcal T : \mathfrak{S} \rightsquigarrow \mathfrak S.
\end{equation}
Then we can apply the Tychonoff fixed point theory (\cite{Coutand2011a}) to derive the well-posedness of \eqref{eq:LagrangianCoordinates}.

\subsection{The Linearized Problem}
We first introduce the linearized problem. Given $ (\bar r, \bar v,\bar\Theta) \in \mathfrak S $, consider the linearized problem,
\begin{equation}\label{eq:linearized}
	\begin{cases}
		\bigl(\dfrac{x}{\bar r}\bigr)^2 \rho_0 \dt v + \bigl(K\dfrac{x^2\rho_0}{\bar r^2 \bar r_x}  \Theta \bigr)_x = - \dfrac{x^2\rho_0}{\bar r^4} \int_0^x s^2 \rho_0(s)\,ds & \\ ~~~~~~~~~~ + (2\mu+\lambda)  \bigl( \dfrac{(\bar r^2 v)_x}{\bar r^2 \bar r_x}\bigr)_x  & x \in (0, R_0),  \\
		c_\nu x^2 \rho_0 \dt \Theta + K \dfrac{x^2 \rho_0}{\bar r^2\bar r_x} \bar \Theta (\bar r^2 v)_x - \bigl(\dfrac{\bar r^2}{\bar r_x} \Theta_x \bigr)_x = \epsilon x^2 \rho_0 & \\ ~~~~ + 2\mu \bar r^2  \bar r_x \bigl( \dfrac{\bar v_x v_x}{\bar r_x^2} + 2  \dfrac{\bar v v}{\bar r^2} \bigr) + \lambda \bar r^2 \bar r_x \bigl( \dfrac{\bar v_x}{\bar r_x} + 2 \dfrac{\bar v}{\bar r} \bigr) \cdot \bigl( \dfrac{v_x}{\bar r_x} + 2 \dfrac{v}{\bar r} \bigr)  & x \in (0, R_0),
	\end{cases}
\end{equation}
and $ r $ is defined by, for $ 0 < t < T $,
\begin{equation}\label{linearized:r}
	r(x,t) = x + \int_0^t v(x,s)\,ds.
\end{equation}
The boundary condition of \eqref{eq:linearized} is given as
\begin{equation}\label{linearized:boundarycondition}
	v(0,t)= 0, ~ \Theta(R_0, t) = 0, ~ \biggl\lbrack (2\mu+\lambda)\dfrac{v_x}{\bar r_x} + 2 \lambda \dfrac{v}{\bar r}\biggr\rbrack (R_0,t) = 0, ~ t \geq 0. 
\end{equation}
To solve the linearized problem, we consider the variables $$ (X,Y) : = (\rho_0 v, x^2 \rho_0 \Theta) $$ and rewrite \eqref{eq:linearized} in the new variables $ X,Y $. Then following standard Galerkin approximation (see \cite{Coutand2011a}),
denote the $ n \geq 1 $ level of approximating solutions as $ (X_n, Y_n) = (\sum_{i=0}^n \lambda_i^n(t)e_i,\sum_{i=0}^n \delta_i^n(t) e_i ) $, where $ (e_n)_{n\in \mathbb N} $ is a  Hilbert basis of $ H_0^1(0,R_0) $ and each $ e_n $ is of class $ H^k(0,R_0) $ for all $ k\geq 1 $. Then $( X_n, Y_n ) $ can be solved through the ODE system, for $ 1 \leq i \leq n $, 
\begin{equation*}
	\begin{cases}
		( \bigl(\dfrac{x}{\bar r}\bigr)^2 X_{nt} , e_i )_{L^2} - ( K\dfrac{Y_n}{\bar r^2 \bar r_x}  ,  e_{i,x} )_{L^2} = - ( \dfrac{x^2\rho_0}{\bar r^4} \int_0^x s^2 \rho_0(s)\,ds , e_i )_{L^2} \\ ~~~~~~~~ - ( (2\mu+\lambda)\dfrac{(X_n/\rho_0)_x}{\bar r_x} + \lambda \dfrac{X_n/\rho_0}{\bar r}, e_{i,x})_{L^2} + 4\mu ((\dfrac{X_n}{\rho_0\bar r})_x,e_i)_{L^2}, \\
		c_\nu ( Y_{nt},e_i)_{L^2} + (K \dfrac{x^2 \rho_0}{\bar r^2\bar r_x} \bar \Theta (\bar r^2 \dfrac{X_n}{\rho_0})_x,e_i )_{L^2} + (\dfrac{\bar r^2}{\bar r_x} ( \dfrac{Y_n}{x^2\rho_0})_x , e_{ix})_{L^2}\\ ~~~~ = \epsilon ( x^2 \rho_0, e_i)_{L^2}   + (2\mu \bar r^2  \bar r_x \bigl( \dfrac{\bar v_x ( X_n / \rho_0)_x}{\bar r_x^2} + 2  \dfrac{\bar v  X_n / \rho_0}{\bar r^2} \bigr),e_i)_{L^2} \\
		~~~~ + (\lambda \bar r^2 \bar r_x \bigl( \dfrac{\bar v_x}{\bar r_x} + 2 \dfrac{\bar v}{\bar r} \bigr) \cdot \bigl( \dfrac{( X_n / \rho_0)_x}{\bar r_x} + 2 \dfrac{ X_n / \rho_0}{\bar r} \bigr) , e_i)_{L^2} ,
	\end{cases}
\end{equation*}
with the initial data
$$ \lambda_i^n(0) = (\rho_0 v_0,e_i)_{L^2}, ~~\delta_i^n(0) = ( x^2 \rho_0 \Theta_0 ,e_i)_{L^2}, $$ where $ (\cdot,\cdot)_{L^2} $ denotes the $L^2$-inner product on $ (0, R_0) $. Here $ (v_0, \Theta_0) $ satisfies $ \mathcal E_1^0 <\infty $. 
By applying the Hardy's inequalities in Lemma \ref{lm:hardy}, the above Galerkin approximating problem is well-defined and can be solved via the ODE theory. In the following, define $ (v,\Theta) = (v_n,\Theta_n) : = ( X_n/\rho_0, Y_n/(x^2\rho_0) ) $. 
Also, define $ r = r_n $ by \eqref{linearized:r}. Then $ (r, v, \Theta ) \in \mathfrak M $.

In the rest of this section, we will sketch some $n$-independent a prior estimates to show $ ( r , v, \Theta) \in \mathfrak S $ for $ 0 < T < \tilde T_* $ where $ \tilde T_* = \tilde T_*(\bar M) > 0 $ for some $ \bar M > 0 $. Then by taking $ n \rightarrow \infty $, we solve the linearized problem \eqref{eq:linearized}. 

For the sake of convenience, denote $$ \mathcal{\bar E}_1 = \mathcal E_1(\bar v, \bar \Theta) \leq \bar M,  \bar \Lambda_0 = \Lambda_0(\bar r, \bar v, \bar \Theta) \leq \bar M, \bar M_0 =  M_0 (\bar r) \leq 2. $$ Following similar steps as in Section \ref{sec:energy}, the energy estimates yield
\begin{align*}
	 &\norm{x\sqrt{\rho_0} v_t}{\stnorm{\infty}{2}}^2 + \norm{x \sqrt{\rho_0}\Theta_t}{\stnorm{\infty}{2}}^2 + \norm{x v_{xt}}{\stnorm{2}{2}}^2 + \norm{v_t}{\stnorm{2}{2}}^2 {\nonumber} \\
	& ~~~~~~ + \norm{x \Theta_{xt}}{\stnorm{2}{2}}^2 \lesssim T(1+C_\delta)P(\bar\Lambda_0)(\mathcal E_1 + \mathcal{\bar E}_1 + 1) + \delta \mathcal{\bar E}_1 \Lambda_0^2 \\
	& ~~~~~~ + P(\bar M_0) \mathcal E_1^0, \\
	&  \norm{xv_x}{\stnorm{\infty}{2}}^2  + \norm{v}{\stnorm{\infty}{2}}^2 + \norm{x \Theta_x }{\stnorm{\infty}{2}}^2 + \norm{x\sqrt{\rho_0}v_t}{\stnorm{2}{2}}^2 {\nonumber} \\
	& ~~~~~~ + \norm{x\sqrt{\rho_0}\Theta_t}{\stnorm{2}{2}}^2 \lesssim T(1+C_\delta)P(\bar\Lambda_0)(\mathcal E_1 + \mathcal{\bar E}_1 + 1) \\
	& ~~~~~~ + \delta \mathcal{\bar E}_1 \Lambda_0^2  + P(\bar M_0) \mathcal E_1^0,\\
	& \norm{\sqrt{\chi\rho_0}v_t}{\stnorm{\infty}{2}}^2 + \norm{\sqrt{\chi}v_{xt}}{\stnorm{2}{2}}^2 + \norm{\sqrt{\chi}\dfrac{v_t}{x}}{\stnorm{2}{2}}^2 \\
	& ~~~~~~ \lesssim T(1+C_\delta)P(\bar\Lambda_0)(\mathcal E_1 + \mathcal{\bar E}_1 + 1) + \delta \mathcal{\bar E}_1 \Lambda_0^2  + P(\bar M_0) \mathcal E_1^0,	\\
	&  \norm{\sqrt{\chi}v_x}{\stnorm{\infty}{2}}^2 + \norm{\sqrt{\chi}\dfrac{v}{x}}{\stnorm{\infty}{2}}^2 + \norm{\sqrt{\chi\rho_0}v_t}{\stnorm{2}{2}}^2 \\
	& ~~~~~~ \lesssim T(1+C_\delta)P(\bar\Lambda_0)(\mathcal E_1 + \mathcal{\bar E}_1 + 1) + \delta \mathcal{\bar E}_1 \Lambda_0^2  + P(\bar M_0) \mathcal E_1^0. 
\end{align*}
We point out the estimates of the terms involving higher order derivative of $ \bar v,\bar \Theta $ in the above estimates. When calculating $ \int \dt \eqref{eq:linearized}_2 \Theta_t $, we engage with the following estimates on the right
\begin{align*}
	& \int \biggl\lbrace 2\mu \bar r^2  \bar r_x \bigl( \dfrac{\bar v_x v_x}{\bar r_x^2} + 2  \dfrac{\bar v v}{\bar r^2} \bigr) + \lambda \bar r^2 \bar r_x \bigl( \dfrac{\bar v_x}{\bar r_x} + 2 \dfrac{\bar v}{\bar r} \bigr) \cdot \bigl( \dfrac{v_x}{\bar r_x} + 2 \dfrac{v}{\bar r} \bigr) \biggr\rbrace_t \Theta_t \,dx \\
	& ~~ \lesssim \mathcal E_1 + \delta  \biggl( \Lambda_0^2 \int ( x^2 \bar v_{xt}^2 + \bar v_t^2) \,dx + \int (x^2 v_{xt}^2 + v_t^2) \,dx   \biggr)\\
	& ~~~~ + (1+C_\delta) P(\bar\Lambda_0) \int x^2 \Theta_t^2 \,dx,
\end{align*}
and on the left
\begin{align*}
	& \abs{\int k\bigl( \dfrac{x^2\rho_0}{\bar r^2 \bar r_x} \bar\Theta (\bar r^2 v)_x\bigr)_t \Theta_t\,dx}{} \lesssim  P(\bar M_0) \Lambda_0 \int x^2 \rho_0  \abs{\bar \Theta_t \Theta_t}{} \,dx \\
	& ~~~~ + (1+C_\delta)  P(\bar\Lambda_0) \int x^2 \Theta_t^2 \,dx + \delta \int (x^2 v_{xt}^2 + v_t^2 )  \,dx + \mathcal E_1 \\
	& ~~ \lesssim \delta \biggl( \Lambda_0^2 \int x^2 \bar \Theta_{xt}^2 \,dx + \int (x^2 v_{xt}^2 + v_t^2) \,dx  \biggr) \\
	& ~~~~ + \mathcal E_1 + (1+C_\delta)  P(\bar\Lambda_0) \int x^2 \Theta_t^2 \,dx.
\end{align*}
Then after applying \eqref{ene:1st-003} and integrating in the temporal variable, one can get the corresponding estimates. 

The next step is to employ similar estimates as in Lemma \ref{lm:pointwise}, one can derive
\begin{equation}\label{linearized:0001}
	\Lambda_0 \lesssim P(\bar M_0) ( C_\delta \mathcal{ E}_1 +  \delta \mathcal{\bar E}_1 + 1).
\end{equation}
Therefore, we have the following inequality
\begin{align*}
	& \mathcal E_1 \leq T(1+C_\delta) P_8(\bar\Lambda_0)(\mathcal E_1 + \mathcal{\bar E}_1 + 1) + \delta P_9(\bar M_0,\mathcal{\bar E}_1) ( \mathcal{ E}_1 +  \mathcal{\bar E}_1 + 1)^2\\
	& ~~~~~~~ + P_{10}(\bar M_0) \mathcal E_1^0,
\end{align*}
for some positive polynomials $ P_8 = P_8(\cdot), P_9= P_9(\cdot), P_{10}= P_{10}(\cdot) $. Then after choosing $ \delta $ so that $ \delta P_9(\bar M_0, \mathcal{\bar E}_1)(\mathcal E_1 + \mathcal{\bar E}_1 + 1 )^2 \leq 1 $, we have
\begin{equation*}
	\mathcal E_1 \leq T P_{11}(\bar M) (\mathcal E_1 + \bar M + 1 )^3 + 1 + C_3 \mathcal E_1^0,
\end{equation*}
for some positive polynomial $ P_{11} = P_{11}(\cdot) $ and some constant $ 0 < C_3 < \infty $. Hence, let
\begin{equation*}
	\tilde T_* := \dfrac{1}{P_{11}(\bar M)(\mathcal E_1 + \bar M + 1 )^3 }.
\end{equation*}
Then for $ T < \tilde T_* $ and $ \bar M > 2 +  C_3 \mathcal E_1^0 $, 
\begin{equation}\label{linearized:0002}
	\mathcal E_1 \leq 1 + 1 + C_3 \mathcal E_1^0 \leq \bar M.
\end{equation}
Also, from \eqref{linearized:0001}, for some positive polynomial $ P_{12} = P_{12}(\cdot) $,  choosing $ \delta \sqrt{\bar M} \simeq 1 $,
\begin{equation*}
	\Lambda_0 \leq P_{12}(2) (C_{\delta} \mathcal E_1 + \delta \bar M  + 1) \leq P_{12}(2) (\sqrt{\bar M} (2 + C_3 \mathcal E_1^0) + \sqrt{\bar M} + 1)\leq \bar M,
\end{equation*}
for $ \bar M $ sufficiently large. Also, applying the fundamental theory of calculus yileds, 
$$ M_0 \leq 1 + (\Lambda_0 + \Lambda_0^2 ) T \leq 2 ~ \text{for} ~ T < 1/({\bar M}+\bar M^2). $$
Therefore, we can refine $ \tilde T_* $ to be
\begin{equation}\label{linearized:0003}
	\tilde T_* = \dfrac{1}{P_{13}(\bar M)},
\end{equation}
for some positive polynomial $ P_{13} = P_{13}(\cdot) $. Finally, the elliptic estimates as in Section \ref{sec:elliptic} will show $ (r,v,\Theta) \in \mathfrak X $ and in particular $ \norm{(r,v,\Theta)}{\mathfrak M} < \tilde P(\bar M) $ for some positive polynomial $ \tilde P = \tilde P(\cdot) $. 

Therefore, after taking $ n \rightarrow \infty $, we get the solution $ ( r, v,\Theta) \in \mathfrak S \cap \mathfrak X $ of the linearized problem \eqref{eq:linearized}. It is easy to verify that the solution is unique.

\subsection{The Fixed Point}
Let $ \mathcal T(\bar r, \bar v,\bar \Theta ) : = (r, v, \Theta) $ be the unique solution to the linearized problem \eqref{eq:linearized}. Then from the last section, $ \mathcal T $ is mapping from $ \mathfrak S $ to $ \mathfrak S $. Also, it is furthermore clear that $ \mathcal T $ is weakly continuous in the $ \mathfrak M $ norm. Then the Tychonoff fixed point theory yields there is at least one solution to the nonlinear problem \eqref{eq:LagrangianCoordinates} in $ \mathfrak S $ for $ T < \tilde T_* $ (defined in \eqref{linearized:0003}). 

Then the a priori estimates in Section \ref{sec:elliptic} show that $ (r, v, \Theta) \in \mathfrak X $. Consequently, we get the well-posedness theory of strong solutions to \eqref{eq:LagrangianCoordinates}. On the other hand, the estimates in Section \ref{sec:classical} show that the strong solution becomes classical if the initial data satisfies $ \mathcal E_2^0 <\infty $. Thus we have the well-posedness theory of classical solutions to \eqref{eq:LagrangianCoordinates}.

\paragraph{Acknowledgements}
This work is part of the doctoral dissertation of the author under the supervision of Professor Zhouping Xin at the Institute of Mathematical Sciences of the Chinese University of Hong Kong, Hong Kong. The author would like to express great gratitude to Prof. Xin for his kindly support and professional guidance.

\bibliographystyle{plain}


\end{document}